\documentstyle[12pt,bezier]{book} 

\newtheorem{Theorem}{Theorem}[section]
\newtheorem{Definition}[Theorem]{Definition}
\newtheorem{Proposition}[Theorem]{Proposition}
\newtheorem{Lemma}[Theorem]{Lemma}

\title{The Method of Characteristics Revisited\\ 
A Viability Approach \\ \mbox{} \\A Mini-Course\footnote{This is the first
preliminary draft of lecture notes of a mini-course at University of
California at Davis in 1999. Any critical comment and correction will be
gratefully acknowledged.}}
\author{Jean-Pierre Aubin}

\begin{document}
\setlength{\unitlength}{1mm}

\pagenumbering{roman}
 \maketitle

{\sc Jean-Pierre Aubin}

\mbox{}

{\sc Universit\'{e} de Paris-Dauphine}

     {\sc Centre de Recherche Viabilit', Jeux, Contr"le}

F-75775 Paris cx(16), France

\begin{center}
{\Large \bf Acknowledgments}  \\
\end{center}

\mbox{}

The author wishes to thank Roger Wets for giving him the opportunity to
present theses lectures at University of California at Davis.

He is grateful to H'lŠne Frankowska for  her contributions about
Hamilton-Jacobi equations and systems of first-order partial differential
equations on which the main results presented here are based and daily
discussions on these topics.

\vspace{ 85 mm}

Typesetting: \LaTeX

 \newpage
\mbox{}  
\begin{center}

>From the same author\\
\end{center}
 \small
\begin{itemize}

        \item    
        {\sc Approximation of Elliptic Boundary-Value Problems} (1972)
         Wiley
	 \item     {\sc Applied Abstract Analysis} (1977)  Wiley-Interscience 
        \item    {\sc Applied Functional Analysis}  (1979) 
        Wiley-Interscience, Second Edition, 1999
 (Version  Fran\c{c}aise: {\sc Analyse Fonctionnelle Appliqu\'{e}e, Tomes
1  \& 2.} (1987) Presses Universitaires de France)

\item     {\sc Mathematical Methods of Game and Economic Theory} (1979) 
         North-Holland Second Edition, 1982
        \item     {\sc M\'{e}thodes Explicites de l'Optimisation} (1982) 
        Dunod (English translation: {\sc Explicit Methods of Optimization},
Dunod, 1985) 
        \item    {\sc L'Analyse non lin\'{e}aire et ses motivations
\'{e}conomiques} (1983) 
         Masson (English translation: {\sc Optima and Equilibria}, Springer
Verlag, 1993, Second Edition, 1998)
\item     {\sc Differential Inclusions} [in collaboration with A. CELLINA],
(1984) 
Springer-Verlag 
\item     {\sc Applied Nonlinear Analysis} [in collaboration with I.
EKELAND],  (1984) Wiley-Interscience
\item    {\sc Exercices d'analyse non lin\'{e}aire}  (1987) 
Masson 
\item   
{\sc  Set-Valued   Analysis} [with H. FRANKOWSKA],  
(1990) Birkh\"auser
\item       {\sc Viability Theory}, (1991)
Birkh\"auser
\item  {\sc Initiation \`a l'Analyse Appliqu\'ee} (1994) Masson (English
Version to appear)
\item      {\sc Neural Networks and Qualitative Physics: 
A Viability Approach} (1996) Cambridge University Press
\item       {\sc Dynamic Economic Theory: A Viability Approach} (1997)
Springer-Verlag (Studies in Economic Theory). Second Edition to appear
\item {\sc Mutational and Morphological Analysis: Tools for Shape
Regulation and Morphogenesis} (1999), Birkh\"auser
\item  {\sc La mort du devin, l'\'emergence du d\'emiurge. Essai sur la
contingence et la viabilit\'e des syst\`emes} (in preparation)
\end{itemize}
\normalsize
\newpage
\tableofcontents
\newpage

\pagenumbering{arabic}
\markboth{Introduction}{Introduction}
\mbox{} 

\vspace{ 12 mm}
{\Huge \bf Introduction}

\vspace{ 6 mm}
 
This mini-course provides a presentation of the method of characteristics 
to  initial/boundary-value  problems  for  systems  of  first-order partial
differential equations and to Hamilton-Jacobi variational inequalities.  

These  results are indeed  useful  for the treatment of  hybrid  systems of
control  theory.  

We use the tools forged by set-valued analysis and viability  theory, which
happen  to be both  efficient and   versatile to cover  many problems. They
find here an unexpected relevance.

Indeed, since solutions to first-order systems 
\begin{displaymath}
\forall \;  j=1, \ldots ,p, \;  \;  \frac{ \partial }{ \partial t}u (t,x) +
\sum_{i=1}^{n}\frac{ \partial }{ \partial x_{i}}u_{i}
(t,x) f_{i}(t,x,u (t,x)) - g_{j}(t,x,u (t,x)) \; = \; 0
\end{displaymath}
may have ``shocks'',  i.e.,  may be set-valued maps (or multi-valued maps),
it seems to us natural to use the concept of ``graphical derivative''  of a
set-valued  map  from  set-valued   analysis  instead  of  ``distributional
derivative''  from  the theory  of  distributions to  give  a meaning  to a
concept of solution to such systems of partial differential equations. 

  The basic  concept useful in our framework is the concept of {\em capture
basin} of a subset $C$  under a differential equation,  which is the set of
points from which a solution  to  the differential equation reaches  $C$ in
finite time. 
        
        Then we shall prove that {\em the graph of the  solution  $ (t,x)
\leadsto U (t,x)$ to the above boundary value problem  is the capture basin
of the graph of the initial/boundary data  under the characteristic system
of differential equations
        \begin{displaymath} \left\{ \begin{array}{ll}
        i) &  \tau ' (t) \; = \; -1 \\
        ii) &  x' (t ) \; = \; -f ( \tau (t),x (t),y (t))) \\
        iii) & y' (t) \; = \; - g ( \tau  (t),x (t),y (t))
        \end{array} \right. \end{displaymath}
         and that, under adequate assumptions,  this solution is {\bf
unique} among the solutions with closed graph to this boundary value
problem}. 
        
        Such a solution is thus taken in a  generalized --- or weak ---
sense  (Frankowska solutions),  since such maps, even when they are
single-valued, are not differentiable in the usual sense.  But when they
are, they naturally coincide with the above concept of solution thanks to
the uniqueness property.

        Existence and uniqueness of the solution is obtained
\begin{enumerate}
        \item from a first characterization of the capture basin of $C$ as
the union of reachable sets $ \vartheta _{-f} (t,C) := \{\vartheta _{-f}
(t,c)\}_{c \in C}$, where $\vartheta _{-f} (t,c)$ denotes the value at time
$t$ of the solution to the differential equation $x'=-f (x)$ starting at
$c$, allowing explicit computations of the capture basins in specific
instances,
        \item from a second characterization of the capture basin of $C$ 
of a closed backward invariant subset $M $ which is a
repeller\footnote{This means that all solutions to the differential
equation $x'=f (x)$  starting from $M$ leave $M$ in finite time.} under a
differential equation $x'=f (x)$. It states that the capture basin is the
{\bf unique} closed subset $K$ such that
        \begin{enumerate}
        \item $C\; \subset  \; K \; \subset  \; M$,
        \item $K$ is {\em backward  invariant} in the sense that any
backward solution starting from $K$ is viable in $K$,
        \item $K \backslash C$ is {\em locally viable} in the sense that
for any $x \in K \backslash C$, there exist $T>0$   and a solution to the
differential equation  starting at $x$ viable in $K \backslash C$ on the
interval $[0,T]$,
        \end{enumerate}
        \item from the 1942 Nagumo Theorem stating that 
        \begin{enumerate}
        \item the  subset $K$ is backward invariant under a differential
inclusion $x' =f (x)$ if and only if, for every $x \in K$, $f (x) \in -
T_{K} (x)$,  or equivalently, if and only if, for every $x \in K $, $p \in
N_{K} (x)$, $  \langle p, f (x) \rangle \geq 0$,
        \item the subset $K \backslash C$ is locally viable if and only if,
for every $x \in K \backslash C$, $f (x)\in T_{K} (x)$, or equivalently, if
and only if, for every $x \in K \backslash C$, $p \in N_{K} (x)$, $ 
\langle p, f (x) \rangle \leq 0$,
        \end{enumerate}
        where the ``contingent cone'' $T_{K} (x)$ to a subset $K$ at a
point $x \in K$, introduced in the early thirties independently by
Bouligand and Severi, adapts to any subset the concept of tangent space to
manifolds, and where the ``normal cone'' $N_{K} (x):= T_{K} (x)^{-}$,
defined as the  polar cone to the contingent cone, adapts to any subset the
concept of normal space to manifolds,
        \item from the definition of the contingent derivative $DU (t,x,y)$
of the set-valued map $U: (t,x) \leadsto  U(t,x)$ at a point $ (t,x,y)$ of
the graph of $U$ as the set-valued map from $ {\bf R} \times X$ to $Y$ the
graph of which is the contingent cone to the graph of $U$ at the point $
(t,x,y)$.
        \end{enumerate}
        
        The above results --- which are interesting by themselves for other
mathematical models of evolutionary economics, population dynamics,
epidemiology ---  can  be applied to many other problems. Dealing with
subsets, they can be applied to graphs of single-valued maps as well as
set-valued maps, to epigraphs and hypographs of (extended) real-valued
functions, to graph of ``impulse'' maps (which take empty values except in
a discrete sets, useful in the study of hybrid systems or inventory
management), etc.
        Since these results are also valid for underlying differential
inclusions, we are able to treat control problems for such boundary-value
problems for systems of first-order partial differential equations.
        
We also illustrate the strategy  of  using the properties  of the viability
kernel    --- the subset of elements from which  starts at least a solution
to the differential equation viable in $K$ --- and the capture basin   of a
subset advocated by H'lŠne Frankowska\footnote{H\'el\`ene Frankowska proved
that the epigraph of the value function  of an optimal control  problem ---
assumed  to  be only lower semicontinuous  ---  is  invariant and  backward
viable   under a (natural) auxiliary system.   Furthermore,    when   it is
continuous,  she  proved  that  its   epigraph    is    viable    and   its
hypograph     invariant     (\cite[Frankowska]{f89hjc,f89hje,HJB92}).    By
duality,   she proved that the latter   property is equivalent  to the fact
that    the value   function  is  a viscosity solution   of  the associated
Hamilton-Jacobi   equation    in the sense    of   M.  Crandall  and  P.-L.
Lions.        See       also      \cite[Barron      \&     Jensen]{bj90hj},
\cite[Barles]{ba93hj}    and   \cite[Bardi   \& Capuzzo-Dolcetta]{bcd98vis}
for   more  details.  Such  concepts  have been extended  to  solutions  of
systems   of   first-order  partial differential equation without  boundary
conditions      by      H\'el\`ene      Frankowska     and    the    author
(see                             \cite[Aubin                             \&
Frankowska]{af90cm1,af90hyp,af91hyp,af90cm,af94cpdi}   and   chapter~8   of
\cite[Aubin]{avt}).  See also \cite[Aubin  \& Da Prato]{adp90hyp,adp92hyp}.
This point of view is used here in the case   of  boundary value problems.}
for characterizing  the value functions  of  some  variational  problems or
stopping  time  problems  as  ``contingent  solutions''  and/or ``viscosity
solutions to
Hamilton-Jacobi ``differential variational inequalities''
        \begin{displaymath} \left\{ \begin{array}{ll}
        i) &  {\bf u} (x) \; \leq  \; {\bf u} ^{\top } (x)\\
         ii)  &   \displaystyle{\left\langle    \frac{ \partial }{ \partial
x}{\bf u}
^{\top } (x),f (x)\right\rangle +l (x,f (x))  +a{\bf u} ^{\top } (x)\; \leq
\;
0}\\
         iii)  &  \displaystyle{({\bf u}  (x)-{\bf u} ^{\top  } (x)) \left(
\left\langle  
\frac{ \partial  }{ \partial x}{\bf  u} ^{\top } (x),f  (x)\right\rangle +l
(x,f
(x))  +a{\bf u} ^{\top } (x) \right) \; = \; 0}
        \end{array} \right. \end{displaymath} 
and
        \begin{displaymath} \left\{ \begin{array}{ll}
        i) &  0 \; \leq  \; {\bf u} ^{\bot }(x) \; \leq  \; {\bf u}  (x)\\
         ii)  &   \displaystyle{\left\langle    \frac{ \partial }{ \partial
x}{\bf u}
^{\bot } (x),f (x)\right\rangle +l (x,f (x))  +a{\bf u} ^{\bot } (x)\; \geq
\;
0}\\
         iii)  &  \displaystyle{({\bf  u} (x)-{\bf u}  ^{\bot } (x)) \left(
\left\langle  
\frac{ \partial }{  \partial x}{\bf u}  ^{\bot }  (x),f (x)\right\rangle +l
(x,f
(x))  +a{\bf u} ^{\bot } (x) \right) \; = \; 0}
        \end{array} \right. \end{displaymath} 
where
\begin{enumerate}
\item $f:{\bf   R}^{n}  \mapsto {\bf   R}^{n}$  defines the dynamics of the
differential equation $x'= f(x)$, 
\item  $l : (x,p) \in {\bf  R}^{n} \times {\bf  R}^{n}  \mapsto l (x,p) \in
{\bf R}_{+}$
is a nonnegative ``Lagrangian'',

\item ${\bf u}:{\bf  R}^{n} \mapsto {\bf R}_{+} \cup \{+\infty \}$
is an extended nonnegative function (regarded as an obstacle, as in
unilateral mechanics).
\end{enumerate}

We shall observe that the epigraphs of  the solutions are  respectively the
viability kernel and the capture basin of the epigraph of $  {\bf u}$ under
the map $g: {\bf  R}^{n}
\times {\bf R} \mapsto {\bf  R}^{n} \times  {\bf R}$ defined by
\begin{displaymath}
 g (x,y) \; := \;  ( f (x), -ay-l (x,f (x)))
\end{displaymath}

This allows  us  to  compute  the  ``solutions''  to  these Hamilton-Jacobi
equations and to check that they are respectively defined by
\begin{displaymath}
{\bf u}^{\top} (x) \; := \;   \sup_{t \geq 0} \left(   e^{a  t}{\bf u} (x
(t)) + \int_{0}^{t}e^{a   \tau }l ( x( \tau ),x' ( \tau )) d \tau \right)
\end{displaymath}
and 
\begin{displaymath}
{\bf u}^{\bot} (x) \; := \;   \inf_{t \geq 0} \left(   e^{a  t}{\bf u}
(x (t)) + \int_{0}^{t}e^{a   \tau }l ( x( \tau ),x' ( \tau )) d \tau
\right) 
\end{displaymath} 
where $x ( \cdot )$ is the solution to the differential equation $x'=f (x)$
starting at $x$.

Using  the fact that the contingent  cone  to the epigraph   of an extended
lower  semicontinuous  function  is  the   epigraph  of   its  ``contingent
epiderivative'',  the  characterizations of  viability  kernel  and capture
basins  in  terms  of  tangential  conditions  allows  to  interpret  these
``solutions''  as  ``contingent  solutions''  to the  above Hamilton-Jacobi
variational  inequalities.  Using  the characterization in terms  of normal
cones,  we obtain  the  equivalent  interpretation  of  these  solutions as
viscosity  solutions whenever the  solution is continuous  instead of being
merely lower semicontinuous.

This this  is  the reason  why this  mini  course begins   with the  Nagumo
Theorem  on  closed  subsets viable  and/or invariant  under a differential
equation.  Next,  the notions of the  ``capture basin'' of a closed subset,
as well as its ``viability kernel''  are characterized and  linked with the
Nagumo  Theorem.  After recalling the concepts  of contingent epiderivative
and subdifferential of  an extended function,  these theorems are  used for
solving some Hamilton-Jacobi  variational inequalities.  They are also used
for solving initial/boundary  value  problems  for  systems  of first-order
partial differential equations.

For  the sake  of  simplicity,  we restrict ourselves  in this introductory
mini-course to the simple  case when the characteristic  system  is made of
differential equations.  However,  the very same methods can be  adapted to
the case when the characteristic system is made of differential inclusions.
Actually, the methods and the view points developed in these notes are more
important  than  the examples  presented here,  since  they  can  be easily
generalized and efficiently used for solving many more difficult problems.

\chapter{The Nagumo Theorem}

\vspace{ 6 mm}
{\Huge \bf Introduction}

\vspace{ 3 mm}
This chapter is  a   presentation  of  the  basic  Nagumo  Theorem  and its
corollaries   and   consequence   in  the  simple   framework  of  ordinary
differential equations $x' = f(x)$. These results --- which are interesting
by  themselves    as  mathematical  metaphors  of  evolutionary  economics,
population dynamics,  epidemiology,  biological  evolution   when  they are
extended to  differential  inclusions  ---   can be  applied to  many other
problems, such as control problems and, as we shall illustrate here, can be
used  as  tools  for  solving  other  mathematical  problems.  Dealing with
subsets,  they can be  applied to graphs of  single-valued maps  as well as
set-valued maps,  to epigraphs  and  hypographs  of  (extended) real-valued
functions,  and  to  be used as  versatile and efficient tools  for solving
systems of first-order partial differential equations.

A  function $[0,T] \ni t \rightarrow x(t)$ is said to be (locally) {\em
viable} in a given subset $K$ on $[0,T]$ if, for any $t \in [0,T]$, the
state $x(t)$ remains in $K$.

Therefore, if a continuous map $f: {\bf  R}^{n} \mapsto {\bf  R}^{n} $
describes the dynamics of the system
\begin{displaymath}
 \forall \; t \geq 0, \;  \; x' (t) \; = \; f (x (t))
\end{displaymath}
we shall say that $K$ is {\em viable under\/} $f$ if starting from any
initial point of $K$, {\em at least one solution} to the differential
equation is viable in  $K$.

The  Nagumo Theorem characterizes such a viability property for any locally
compact subset $K$ in terms of contingent and normal cones.

Hence, we begin by recalling the concept of contingent cones and normal
cones to a arbitrary subset of a finite dimensional vector space. For that
purpose, we start with the notion of upper convergence of sets introduced
by Painlev\'e: The upper limit of a sequence of subsets $K_{n} \subset X$
is the set of cluster points of sequences $x_{n} \in K_{n}$. With that
concept, we can define  {\bf the contingent cone} $T_{K} (x)$ to a subset
$K$ at $x$ is the upper limit of the ``difference quotients'' $
\frac{K-x}{h}$. The normal cone $N_{K} (x)$ is next defined as the polar
cone $T_{K} (x)^{-}$ to the contingent cone.

The  Nagumo Theorem characterizes such a viability property for any locally
compact subset $K$ by stating that $K$ is viable under $f$ if and only if
\begin{displaymath}
 \forall \; x \in K, \;  \; f (x) \; \in  \; T_{K} (x)
\end{displaymath}
or, equivalently, if and only if
\begin{displaymath}
 \forall \; x \in K, \; \forall \; p \in N_{K} (x)\;  \;  \langle p,f (x)
\rangle  \; \leq   \; 0
\end{displaymath}
 
Since open subsets, closed subsets, intersection of open and closed subsets
of $ {\bf  R}^{n} $ are locally compact, we shall be able to specify this
theorem in each of theses cases.

We then state and prove these viability theorems. Many proofs of viability
theorems are now available: We chose the most elementary one (which is not
the shortest) because it is the prototype of the extensions of the
viability theorems. It is just a modification of the Euler method of
approximating a solution by piecewise linear functions (polygonal lines) in
order to force the solution to remain viable in $K$.

\section{Viability \& Invariance Properties}

 \begin{Definition}[Viable functions]
Let $K   $  be a subset of a finite dimensional vector-space  ${\bf  R}^{n}
$. We shall
  say that a function $x(\cdot   )   $  from $[0,T] $ to ${\bf  R}^{n} $ is
{\em viable 
in $K$ on $[0,T]$}\index{viable function}
  if for all $t \in [0,T],\;\; x(t) \in K $.
  \end{Definition}

  Let  us describe the   dynamics of the system  by a   map $f$  from ${\bf
R}^{n}$  to
${\bf  R}^{n}  $. We consider the initial value problem (or Cauchy problem)
associated
with the differential equation

 \begin{equation} \label{01A11} 
  \forall \;  t \in [0,T],\; \;x'(t)\;=\;f(x(t  )  )    \end{equation}
 satisfying the initial condition $x(0) = x_{0}$.

 \begin{Definition} [Viability \& Invariance Properties]
\label{02A11} 
Let  $K \subset$ ${\bf R}^{n} $.  We shall say that $K$ {\em is
locally viable
under $f$\/}  if for any initial state $x_{0 }$ of $K$, {\em there exist }
$T >0$ and a {\em viable }solution on $[0,T]$ to differential equation
(\ref{01A11}) starting at $x_{0 }$.\index{viability property} It is said to
be (globally) viable  under $f$   if we can always take $T = \infty$.

 The subset $K  $ is said to be locally {\em invariant under $f$}  if for
any initial state $ x_{0 }$ of $K  $ and for  {\em  all } solutions $x (
\cdot )$ to differential  equation (\ref{01A11}) (a priori defined on ${\bf
R}^{n} $)
starting from $x_{0}$, there exists $T>0$ such that $x ( \cdot )$ is {\em
viable in $K$} on $[0,T]$. \index{ invariance property} It is said to be
(globally) invariant under $f$   if we can always take $T = \infty$ for all
solutions.

A subset $K$ is a {\sf repeller} \index{repeller} if   from any initial
element $x_{0} \in K$, all solutions to the differential equation
(\ref{01A11}) starting at $x_{0} \in K$ leave $K$ in finite time.

  \end{Definition}

 {\bf Remark} --- \hspace{ 2 mm}
 We should emphasize that the concept of invariance {\em depends upon the
  behavior of  $f   $  on the domain ${\bf   R}^{n}  $ outside $K$}. But we
observe that
viability property depends only on the behavior of $f $ on $K $.
   $\; \; \Box$  

So, the viability property requires only the existence of at least one 
viable solution whereas the invariance property demands that all solutions,
if any,  are viable. 
        
        Observe also that whenever there exists a unique solution to
differential equation $x'=f (x)$ starting from any initial state $x_{0}$,
then viability and invariance properties of a closed subset $K$ are
naturally equivalent.

\vspace{ 5 mm}

  We shall begin by characterizing the subsets $K  $ which are viable under
$f$. The idea is simple, intuitive and makes good sense:
 {\em A subset $K $ is viable under $f$ if at each state $x $  of $K  $,
the velocity $f(x  )  $ is ``tangent'' to $K$ at $x$, so to speak, for
bringing back a solution to  the differential equation inside $K $}.
        
        \section{Contingent and Normal Cones} 
        \subsubsection{Limits of Sets} 
        
        Limits of sets have been introduced by  Paul Painlev\'e in 1902 
without the concept of topology. They have been  popularized by 
Kuratowski  in his famous book {\sc Topologie} and thus, often called {\em
Kuratowski lower and upper limits\/} of sequences of sets. They are defined
without the concept of a topology on the power space. 
        
        \begin{Definition} \index{limit of sets} \label{04A391}
        Let $(K_{n})_{n \in {\bf N}}$ be a sequence of subsets of a finite
dimensional vector space $X$. We say that the subset
        \begin{displaymath}
        \mbox{\rm Limsup}_{n \rightarrow \infty}K_{n}   \;  := \; \left\{ x
\in X \;\; | \;\; \liminf_{n \rightarrow \infty}d(x,K_{n}) = 0 \right\}
        \end{displaymath}
         is the {\em  upper limit} or {\em  outer limit}\footnote{The terms
outer and inner limits of sets have been proposed by R.T Rockafellar and 
R. Wets.} \index{upper limit of sets}\index{outer limit of sets}  of the
sequence $K_{n}$ and that the subset
        \begin{displaymath}
        \mbox{\rm Liminf}_{n \rightarrow \infty}K_{n}   \; 
          := \; \left\{ x \in X  \;\; | \;\; \mbox{\rm lim}_{n \rightarrow
\infty}d(x,K_{n}) = 0 \right\}  
        \end{displaymath}
        is its {\em  lower limit} or {\em  inner limit}. \index{lower limit
of sets} \index{inner limit of sets}
         A subset $K$ is said to be the {\em  limit\/}  or the {\em set
limit\/} of the sequence $K_{n}$ if
        \begin{displaymath}
        K   \; =   \; \mbox{\rm Liminf}_{n \rightarrow \infty }K_{n} \;  
=  \; \mbox{\rm Limsup}_{n \rightarrow \infty } K_{n}   \;  =: \;  
\mbox{\rm Lim}_{n \rightarrow \infty }K_{n} 
        \end{displaymath}
        \end{Definition}
        Lower and upper limits are obviously {\em closed}.
         We also see at once that
        \begin{displaymath}
        \mbox{\rm Liminf}_{n \rightarrow \infty }K_{n} \;  \subset \; 
\mbox{\rm Limsup}_{n \rightarrow \infty }K_{n}
        \end{displaymath}
        and that the  upper limits and  lower limits of the subsets $K_{n}$
and of their closures $\overline{K}_{n}$  do coincide, since $d(x,K_{n}) =
d(x,\overline{K}_{n})$.

        Any decreasing  sequence of subsets $K_{n}$ has a limit:
        \begin{displaymath}
         \mbox{\rm if} \;  \; K_{n} \subset K_{m} \;  \; \mbox{when} \;  \;
n \geq m, \;  \; \mbox{then} \;  \; \mbox{\rm Lim}_{n \rightarrow \infty
}K_{n} \; =  \; \bigcap_{n \geq 0}\overline{K_{n}}
        \end{displaymath}
        An upper limit may be empty (no subsequence  of elements $x_{n} \in
K_{n}$ has a cluster point.)
        
        Concerning sequences of singletons $\{x_{n}\}$, the set limit, when
it exists, is either empty (the sequence of elements $x_{n}$ is not
converging), or is a singleton made of the limit of the sequence.
        
        It is easy to check that:
        \begin{Proposition} \label{02A391}
        If $(K_{n})_{n \in {\bf N}}$ is a sequence of subsets of a finite
dimensional vector space, then
         $\mbox{\rm Liminf}_{n \rightarrow \infty}K_{n}$ is the {\em set of
limits of sequences} $x_{n} \in K_{n}$  and $\mbox{\rm Limsup}_{n
\rightarrow \infty}K_{n}$ is the set {\em of cluster points of sequences}
$x_{n} \in K_{n}$, i.e., of limits of subsequences $x_{n'} \in K_{n'}$.
\end{Proposition}

        \subsubsection{Contingent Cone}
        We reformulate the definition of contingent direction to
a subset of a finite dimensional vector space  introduced
independently by Georges Bouligand and Francesco Severi
\footnote{The concepts of {\em semitangenti\/} and of {\em corde
improprie\/} to a set at a point of its closure had been
introduced by the Italian geometer Francesco Severi (1879-1961)
and are equivalent to the concepts of {\em contingentes\/} and
{\em paratingentes\/} introduced independently by the French
mathematician Georges Bouligand, slightly later. Severi explains
for the second time that he had discovered these concepts
developed by  Bouligand in {\em `` suo interessante libro
recente''} and comments: {\em ``All'egregio geometra \`e
evidentemente sfuggito che le sue ricerche in proposito sono
state iniziate un po' pi\`u tardi delle mie ... Ma non gli muovo
rimprovero per questo, perch\'e neppur io riesco a seguire  con
cura minuziosa la bibliografica e leggo pi\`u volontieri una
memoria o un libro dopo aver pensato per conto mio
all'argomento.''} I am grateful to M. Bardi for this information
about Severi.} in the 30's:

        \begin{Definition}
         When $K \subset X$ is a subset of a normed vector space $X$ and
when $x \in K$, the set $T_{K} (x)$ 
        \begin{displaymath}
         T_{K}(x) \; := \; \mbox{\rm Limsup}_{h \rightarrow 0+}
\frac{K-x}{h}
        \end{displaymath}
        of {\em contingent directions\/} to $K$ at $x$ is a closed cone,
called the {\em contingent cone} or simply, the {\em tangent cone}, to $K$
at $x$.
         \end{Definition}

Therefore, a direction $v \in X$ is {\em contingent} to $K$ at $x$ if and
only if
\begin{displaymath}
 \liminf_{h \rightarrow 0+} \frac{d (x+hv,K)}{h} \; = \; 0
\end{displaymath}
or, equivalently, if and only if there exists a sequence of elements $h_{n}
>0$ converging to $0$ and a sequence of $v_{n} \in X$ converging to $v$
such that
\begin{displaymath}
 \forall \; n \geq 0, \;  \; x + h_{n}v_{n} \; \in  \; K
\end{displaymath}

The lemma below shows right away why these cones will play a crucial role:
they appear naturally whenever we wish to differentiate viable functions.
\begin{Lemma} \label{derviablem}
Let $x(\cdot)$ be a differentiable {\em viable} function from $[0,T]$ to
$K$. Then $$\forall \;  t \in [0,T[,\;  \;  \; x'(t) \; \in \;
T_{K}(x(t))$$
 \end{Lemma}
{\bf Proof} --- \hspace{ 2 mm}
        Let us consider a function $x(\cdot) $ viable in $K$. It is easy to
check that $x'(0)$ belongs
to the contingent cone $T_{K}(x_{0})$ because $x(h)$ belongs to $K$ and
consequently, 
        \begin{displaymath}
        \frac{ d_{K}(x_{0}+ hx'(0))}{h} \; \leq  \; \frac{\|x(0)+h x'(0)
- x(h)\|}{h} \;\mbox{\rm converges to}\; 0
        \end{displaymath}
        Hence $x' (0)$ belongs to  the contingent cone to $K$ at $x_{0}$.
$\; \; \Box$

        For convex subsets $K$, the  contingent cone coincides with the
closed cone spanned by $K-x$:
        \begin{Proposition} \index{tangent cones to convex sets}
\label{01A955} 
        Let us assume that $K$ is convex. Then the contingent cone
$T_{K}(x)$ to $K$ at $x$ is convex and
        \begin{displaymath}
         T_{K}(x) \; =  \; \overline{\bigcup_{h>0}^{} \frac{K-x}{h}}  
        \end{displaymath}
        We shall say in this case that it is  the   {\em tangent cone} to
the convex subset $K$ at $x$. 
        \end{Proposition}

        {\bf Proof} --- \hspace{ 2mm} 
        We begin by stating the following consequence of convexity: If $0 <
h_{1} \leq h_{2}$, then
        \begin{displaymath}
         \frac{K-x}{h_{2}} \; \subset  \;  \frac{K-x}{h_{1}}
        \end{displaymath}
        because $x + h_{1}v = \frac{h_{1}}{h_{2}} (x+h_{2}v)+ \left(
1-\frac{h_{1}}{h_{2}} \right)x$ belongs to $K$ whenever $x+h_{2}v$ belongs
to $K$. The sequence of the subsets $\frac{K-x}{h}$ being increasing, our
proposition ensues.
        $\; \; \Box$

 \begin{Definition} [Viability Domain]
 Let $K  $ be a subset of ${\bf  R}^{n} $. We shall say that $K $ is a {\em
 viability domain}\index{ viability domain for single-valued maps}  of the
map $f :{\bf  R}^{n}  \mapsto {\bf  R}^{n} $ if
 \begin{equation}
\forall \;  x \in K , \; \; f(x) \; \in \; T_{K}(x)
\end{equation}
  \end{Definition} 

We recall that a subset $K \subset {\bf  R}^{n}$ is {\em locally compact\/}
if there
exists $r>0$ such that the ball $B_{K}(x_{0},r) := K \cap (x_{0}+rB)$ is
{\em compact}. Closed subsets, open subsets and intersections of closed and
open subsets of a finite dimensional vector space are locally compact.

        \subsubsection{Normals}
\begin{Definition} \label{normalcondef}
        The polar cone $N_{K} (x) :=\left( T_{K} (x)   \right)^{-}$  is
called the {\em normal cone}  to $K$ at $x$. 
 \end{Definition}

It is also called the {\em Bouligand normal cone}, or the {\em contingent
normal cone}, or also, the {\em sub-normal cone} and more recently, the
{\em regular normal cone} by R.T. Rockafellar and R. Wets. In this book,
only (regular) normals are used, so that we shall drop the adjective
``regular''.
        $\; \; \Box$ 
        
        \mbox{}
        
        \begin{Lemma} \label{localsep}
        Let $K \subset X$ be a closed subset. Let $y \notin K$ and $x \in 
\Pi _{K} (y)$ a best approximation of $y$ by elements of $K$: $ \|y-x\| =d
(y,K)$. Then
        \begin{displaymath}
      \forall \; y \notin K, \; \forall \; x \in  \Pi_{K} (y), \;  \;\; y-x
\; \in \; N_{K} (x)
     \end{displaymath}
        \end{Lemma}
        
        {\bf Proof} --- \hspace{ 2 mm}
         Since $x \in  \Pi _{K} (y)$ minimizes the distance $ z \mapsto
\frac{1}{2}\|y-z\|^{2}$ on $K$, we deduce that
        \begin{displaymath}
        \forall \; v \in T_{K} (x), \;  \; 0 \; \leq  \;   \left\langle 
x-y , v \right\rangle  
        \end{displaymath}
        so that $y-x$ belongs to $T_{K} (x)^{-}=:N_{K} (x)$. $\; \; \Box$ 

\mbox{}

When $K$ is convex, we deduce that $N_{K} (x)$ is the polar cone to $ K-x$
because the tangent cone  is spanned by $K-x$:

 \begin{Theorem} \label{03A955}
        Let $K $ be a closed convex subset of a  finite dimensional vector
space. Then 
        \begin{displaymath}
 p \in N_{K} (x)  \;\; \mbox{if and only if} \;\; \forall \; y \in K, \; 
\; \langle p, y \rangle \; \leq  \; \langle p, x \rangle  
\end{displaymath}
        and the graph of the set-valued map $x \leadsto N_{K} (x)$  is
closed in $X \times X^{\star}$.
        \end{Theorem}
        {\bf Proof} --- \hspace{ 2 mm}
        Let us consider sequences of elements $x_{n} \in K$  converging to
$x$ and $p_{n} \in N_{K}(x_{n})$ converging to $p$. Then inequalities
        \begin{displaymath}
        \forall \; y \in K,  \;\;  \langle  p_{n},y\rangle   \;  \leq  
\;   \langle p_{n},x_{n} \rangle 
        \end{displaymath} 
        imply by passing to the limit inequalities
        \begin{displaymath}
        \forall \; y \in K,  \;   \;   \langle p,y \rangle   \;  \leq  
\;   \langle p,x \rangle 
        \end{displaymath} 
        which state that $p$ belongs to $N_{K}(x)$. Hence the graph is
closed. $\; \; \Box$

In the general case, we provide now the following characterization of the
normal cone:
        
        \begin{Proposition} \label{mormconechar}
        Let $K$ be a subset of  a finite dimensional vector-space $X$. Then
$p \in N_{K} (x) $ if and only 
        \begin{equation} \left\{ \begin{array}{l}  \label{boh} 
        \forall \; \varepsilon >0, \; \exists \; \eta >0 \;\; \mbox{\rm
such that} \;\;\forall  \;   y \in K\cap B(x,\eta ),\\
         \\
         \langle p,y-x \rangle  \;  \leq \;  \varepsilon \|y- x\|
        \end{array} \right. \end{equation}
        \end{Proposition}
        
        {\bf Proof} --- \hspace{ 2 mm}
        Let  $p$ satisfy above property (\ref{boh}) and  $v \in  T_{K}
(x)$. Then there exist $h_{n}$ converging to $0$ and $v_{n}$ converging to
$v$ such that $$y \; := \;  x + h_{n} v_{n} \;  \in \;  K \cap B(x,\eta )$$
for $n$ large enough. Consequently, inequalities $\langle p,v_{n} \rangle 
\leq \varepsilon $ imply by taking the limit that $\langle p,v \rangle 
\leq \varepsilon $ for all $\varepsilon >0$. Hence $\langle p,v \rangle 
\leq 0$, so that any element $p$ satisfying the above property belongs to
the polar cone of $T_{K} (x)$.
        
        Conversely, assume that  $p$ violates property  (\ref{boh}): There
exist $\varepsilon >0$ and a sequence of elements $x_{n} \in K$ converging
to $x$ such that $$\left\langle p,x_{n}-x  \right\rangle \;  > \; 
\varepsilon \|x_{n} - x\|$$ We set $h_{n} := \|x_{n}- x\|$, which converges
to  $0$, and $v_{n} := (x_{n}-x)/h_{n}$. These elements belonging to the
unit sphere, a subsequence (again denoted) $v_{n}$ converges to some $v$.
By definition, this limit belongs to $T_{K} (x)$, so that $\langle p,v
\rangle  \leq 0$. But our choice implies that $\langle p,v_{n} \rangle 
>\varepsilon $, so that $\langle p,v \rangle  \geq \varepsilon $, a
contradiction.
        $\; \; \Box$

We provide a useful characterization by duality of viability domains in
terms of normal cones:
\begin{Theorem} \label{3russianfrank}
Let $K$  be a locally compact subset  of a finite dimensional  vector space
${\bf  R}^{n}$ and $f:K \mapsto {\bf  R}^{n}$ be a continuous single-valued
map. Then
\begin{equation}
\forall \;  x \in K , \; \; f(x) \; \in \; T_{K}(x)
\end{equation}
if and only if
\begin{equation}
\forall \;  x \in K , \; \; f(x) \; \in \; \overline{ \mbox{\rm
co}}(T_{K}(x))
\end{equation}
or, equivalently, in terms of normal cone, if and only if
\begin{displaymath}
 \forall \; x \in K, \;  \forall \; p \in N_{K} (x), \;  \;   \langle  p,f
(x)  \rangle \; \leq \; 0
\end{displaymath}
\end{Theorem}
 
{\bf Proof} --- \hspace{ 2 mm}
Since the normal cone $N_{K} (x)$ is the polar cone to $T_{K} (x)$, and
thus, to $\overline{ \mbox{\rm co}}(T_{K}(x))$, then $\overline{ \mbox{\rm
co}}(T_{K}(x))$ is the polar cone to $N_{K} (x)$, so that the two last
statements are equivalent by polarity. Since the first statement implies
the second one, it remains to prove that if for any $x \in K$, $f (x)$
belongs to $\overline{ \mbox{\rm co}}(T_{K}(x))$, then $f (x)$ is actually
contingent to $K$ at $x$ for any $x \in K$.
This follows from the following 
\begin{Lemma} \label{01A142}
Let  $K  \subset {\bf   R}^{n}$  be a  locally compact   subset of a finite
dimensional vector space  ${\bf  R}^{n}$ and $f:K
\mapsto {\bf  R}^{n}$ be a continuous single-valued map.
Assume that there exists $ \alpha >0$ such that 
\begin{displaymath}
  \forall  \;  x \in K \cap  B (x_{0},  \alpha ),  \;   \;  f  (x)\; \in \;
\overline{
\mbox{\rm co}}(T_{K}(x))
\end{displaymath}
Then, $f (x)$ is contingent to $K$ at elements $x$ of a
neighborhood of $x_{0}$. Actually, for any $ \varepsilon >0$, there exists
$ \eta (x_{0}, \varepsilon )\in ]0, \alpha ]$ such that
\begin{equation} \label{01A1421}
 \forall \; x \in B (x_{0}, \eta (x_{0}, \varepsilon )), \; \forall \; h
\leq \eta (x_{0}, \varepsilon ), \;  \; d \left( f (x), \frac{ \Pi _{K}
(x+hf (x))-x}{h} \right)  \; \leq  \; \varepsilon 
\end{equation}
where $  \Pi_{K} (y)$  denotes the set of best approximations $ \bar{x} \in
K$ of $y$, i.e., the solutions to $ \| \bar{x}-y\|=d (y,K)$.
\end{Lemma}
{\bf Proof} --- \hspace{ 2 mm}
It is sufficient to check the Lemma when $f (x) \ne 0$. Let us set
\begin{displaymath}
 g (t) \; := \; \frac{1}{2}d (x+tf (x),K)^{2} \; = \; \|x+tf
(x)-x_{t}\|^{2}
\end{displaymath}
where $x_{t} \in \Pi_{K} (x+tf (x))$  is a best approximation of $x+tf (x)$
by elements of $K$.
We take $ \alpha $ small enough for $K \cap B (x_{0}, \alpha )$ to be
compact.
We observe that there exists $ \beta \in ]0, \alpha ]$ such that for all $x
\in B (x_{0}, \beta )$, $ \|f (x)\| \leq 2 \|f (x_{0})\|$ because $f$ is
continuous at $x_{0}$. Furthermore, 
\begin{displaymath}
 \|x+tf (x)-x_{t}\| \; = \; d (x+tf (x),K) \; \leq  \; t \|f (x)\| \; \leq 
\; 2 t\|f (x_{0})\|
\end{displaymath}
because $x$  belongs to $K \cap   B (x_{0}, \beta )$, so that $ \|x-x_{t}\|
\leq
2t \|f (x)\|\leq 4t \|f (x_{0})\|$ converges to $0$ with $t$.

On the other hand, for every $v_{t} \in T_{K} (x_{t})$, there exists a
sequence of $h_{n}>0$ converging to $0$ and $v_{t}^{n}$ converging to
$v_{t}$ such that $x_{t}+h_{n}v_{t}^{n}$ belongs to $K$. Therefore,
\begin{displaymath}
 g (t+h_{n})-g (t) \; \leq  \; \frac{1}{2} \left( \|x+tf (x)-x_{t}+ h_{n}(f
(x)-v_{t}^{n})\|^{2}- \|x+tf (x)-x_{t}\|^{2} \right)
\end{displaymath}
and thus, dividing by $h_{n}>0$ and letting $h_{n}$ converge to $0$, that
\begin{displaymath}
\forall \; v_{t} \in T_{K} (x_{t}), \;  \; g' (t) \; \leq  \;
\left\langle   x+tf (x)-x_{t}, f (x)-v _{t} \right\rangle 
\end{displaymath}
Since it is true for any $v_{t} \in T_{K} (x_{t})$ and since the right-hand
side is affine with respect to $v_{t}$, we deduce that this inequality
remains true for any $v_{t} \in \overline{ \mbox{\rm co}}(T_{K} (x_{t}))$,
and thus, by assumption, for $f (x_{t})\in\overline{ \mbox{\rm co}}(T_{K}
(x_{t}))$:
\begin{displaymath}
 g' (t) \, \leq  \, \left\langle   x+tf (x)-x_{t}, f (x)-f
(x_{t})\right\rangle   \leq  2t \|f (x_{0})\| (\|f (x)-f (x_{0})\|
+\|f (x_{0})-f (x_{t})\|)
\end{displaymath}
For any $ \varepsilon $, let $ \gamma  >0$ such that $\|f (y)-f (x_{0})\|
\leq  \frac{\varepsilon^{2}}{8 \|f (x_{0})\|}$ whenever $y \in B (x_{0}, 2
\gamma   )   $. Since  $ \|x_{0}-x_{t}\| \leq \|x- x_{0}\|+4t \|f (x_{0})\|
\leq 2\gamma $ whenever $ \|x-x_{0}\| \leq  \gamma $ and $t \leq \frac{
\gamma }{4 \|f (x_{0})\|}$, then, setting $ \eta (x_{0}, \varepsilon ):=
\min \left( \gamma , \frac{ \gamma }{4 \|f (x_{0})\|} \right)$, we obtain
\begin{equation}
 \forall \; x \in K \cap B (x_{0}, \eta (x_{0}, \varepsilon )), \; \forall
\; t \in ]0, \eta (x_{0}, \varepsilon )], \;  \; g' (t) \; \leq  \;
\frac{t}{2} \varepsilon^{2}
\end{equation}

 Therefore, after integration from $0$ to $h$, we obtain
\begin{displaymath}
 \forall \; h \in [0, \eta (x_{0}, \varepsilon )], \;  \; g (h) -g (0) \;
\leq  \; h^{2} \varepsilon^{2} 
\end{displaymath}
Observing that $g (0)=0$, we derive the conclusion of the Lemma: $\forall
\; x_{h} \in  \Pi _{K} (x +hf (x))$,
\begin{displaymath}
d\left(  f (x), \frac{ \Pi _{K} (x+hf (x))-x}{h} \right) \; \leq  \; \frac{
\|x+h f (x)-x_{h}\|}{h} \;  =  \;  \frac{d  (x+hf (x),K)}{h}    \; \leq  \;
\varepsilon \; \; \Box
\end{displaymath}

\section{Statement of the Viability Theorems}
 Nagumo was the first one to prove the viability theorem 
 for ordinary differential equations in 1942. This theorem was apparently
 forgotten, for it was rediscovered many times during the next twenty
years.

We shall prove it when the subset $K \subset {\bf   R}^{n}$ is {\em locally
compact\/}:

 \begin{Theorem} [Nagumo] \index{Nagumo's Theorem}  \label{01A43}
 Let us assume that
\begin{equation} \left\{ \begin{array}{ll}
i) & K  \; \mbox{is locally compact}   \\
 ii) & f  \;\mbox{is continuous from $K $  to ${\bf  R}^{n} $ }
\end{array} \right. \end{equation}  
Then $K $  is locally viable under $f$ if and only if $K $ is
 a viability domain of $f$ in the sense that
\begin{displaymath}
 \forall \; x \in K, \;  \; f (x) \; \in  \; T_{K} (x)
\end{displaymath}
or, equivalently, in terms of normal cones, 
\begin{displaymath}
 \forall \; x \in K, \;  \forall \; p \in N_{K} (x), \;  \;   \langle  p,f
(x)  \rangle \; \leq \; 0
\end{displaymath}
  \end{Theorem}

 Since the contingent cone to an open subset is equal to the whole space,
an open subset is  a viability domain of any map. So, it is viable under
any continuous map  because any open subset of a finite
 dimensional vector space is locally compact. The Peano existence theorem
is then a consequence of Theorem~\ref{01A43}.

 \begin{Theorem}[Peano]\index{Peano's Theorem} \label{03A13}
   Let  $\Omega $  be an open subset  of a  finite dimensional vector space
${\bf  R}^{n}$
 and $f  :\Omega  \mapsto {\bf  R}^{n} $ be a continuous map.

Then, for every 
 $x_{0}  \in \Omega$, there exists $T > 0$ such that differential
 equation~(\ref{01A11}) has a solution on the interval $[0,T] $ starting at
$x_{0}$. 
 \end{Theorem}

If $C \subset K \subset {\bf  R}^{n}$ is a closed subset of a closed subset
$K$ of a
finite dimensional vector space, then $K \backslash C$ is locally compact,
because for any $x \in K \backslash C$, there exists $r>0$ such that $K
\cap B (x,r) \subset {\bf  R}^{n} \backslash C$. On the other hand,
\begin{displaymath}
 \forall \; x \in K \backslash C, \;  \; T_{K \backslash C} (x) \; = \;
T_{K} (x)
\end{displaymath}

Therefore, Theorem~\ref{01A43} implies 
 \begin{Theorem} \label{03A131}
   Let $C \subset K \subset {\bf   R}^{n}$  is a closed subset  of a closed
subset $K$ of
a finite dimensional vector space  ${\bf  R}^{n}$
 and $f  : K \backslash C \mapsto {\bf  R}^{n} $ be a continuous map.
Then $K \backslash C$ is locally viable under $f$ if and only if
\begin{displaymath}
 \forall \; x \in K \backslash C, \;  \; f (x) \; \in  \; T_{K} (x)
\end{displaymath}
or, equivalently, in terms of normal cones, 
\begin{displaymath}
 \forall \; x \in K \backslash C, \;  \forall \; p \in N_{K} (x), \;  \;  
\langle  p,f (x)  \rangle \; \leq \; 0
\end{displaymath}

 \end{Theorem}

  The interesting case from the viability point of view is the one when the
 viability subset is {\em closed}. This is possible because any closed
subset of a finite dimensional vector space is locally compact. However, in
this case, we derive from  Theorem~\ref{01A43} a more precise statement.

 \begin{Theorem}[Viability]\index{viability theorem for differential
equations} \label{03A43} 
 Let us consider a {\em closed}  subset $K $ of a finite dimensional
  vector space $  {\bf  R}^{n}$ and a {\em continuous }  map $f $ from $K $
to ${\bf  R}^{n} $.
Then $K$ is locally viable under $f$ if and only if
\begin{displaymath}
 \forall \; x \in K , \;  \; f (x) \; \in  \; T_{K} (x)
\end{displaymath}
If such is the case, then for every initial state $x_{0} \in K$, there
exist  a positive $T $ and a {\em viable} solution on $[0,T[ $ to
differential equation~(\ref{01A11}) starting at  $x_{0}$ such that 
\begin{equation} \left\{ \begin{array}{cll}
  \mbox{either}   & T = \infty & \\ 
 \mbox{or}  & T < \infty & \mbox{and}  \;\;  \limsup_{t\rightarrow T-}\|
x(t ) \|  = \infty
\end{array} \right. \end{equation} 
  \end{Theorem}
 
  Further adequate information --- a priori estimates on the growth of
 $f $ --- allows us to exclude the case when \mbox{$\limsup_{t\rightarrow
T-}\| x(t ) \|  = \infty $ }.

This is the case for instance when $f$ is
 bounded on $K $, and, in particular, when $K $  is bounded. 

More generally, we can take $T = \infty $  when $f $ enjoys  linear growth:
 \begin{Theorem} \label{04A43}
Let us consider a subset $K $  of a finite dimensional vector  space $ {\bf
R}^{n}$ and
a   map $f $  from $K $  to ${\bf  R}^{n} $. We assume that the map $f $ is
{\em
continuous} from $K$  to ${\bf  R}^{n} $, that it has {\em linear growth\/}
in the
sense that
\begin{displaymath}
\exists \; c > 0   \;\; \mbox{ such that} \;\;\forall \;  x \in K,\;\; 
\|f(x ) \| \; \leq \; c(\| x\|+1  )   
\end{displaymath}
If
\begin{displaymath}
 \forall \; x \in K , \;  \; f (x) \; \in  \; T_{K} (x)
\end{displaymath}
then $K$ is globally viable under $f$: For every  initial state $x_{0}  \in
K $, there exists a {\em viable solution on} $[0,\infty] $ to differential
equation~(\ref{01A11}) starting at $x_{0}$ and satisfying
\begin{displaymath}
 \forall \; t \geq 0, \;  \; \|x (t)\| \; \leq  \;  \|x_{0}\| e^{ ct}
+e^{ct}-1
\end{displaymath}
  \end{Theorem} 
        
        \section{Proofs of the Viability Theorems}
         We  shall  begin  by  proving  Theorem~\ref{01A43}.  The necessary
condition follows from Lemma~\ref{derviablem}.
        
For proving the sufficient condition,  we begin by constructing approximate
solutions by
modifying the classical the Euler method to take into account the viability
constraints, we
then deduce from available estimates that a subsequence of these solutions
converges uniformly to a limit, and finally check that this limit is a
viable solution to differential equation (\ref{01A11}.)

        1. \hspace{ 2 mm} --- \hspace{ 2 mm} {\bf Construction of
Approximate Solutions} 
        
        Since $K$ is locally compact, there exists $r>0$ such that the ball
$B_{K}(x_{0},r) := K \cap (x_{0}+rB)$ is {\em compact}. When $C$ is a
subset, we set 
        \begin{displaymath}
        \|C\| \; := \; \sup_{v \in C}\|v\| 
        \end{displaymath}
        and 
        \begin{displaymath}
        K_{0}  \; :=  \;K \cap B(x_{0},r), \;C  \;:=  \; B(f(K_{0}),1), \;
T \;  := \; \frac{r}{\|C\|}
        \end{displaymath}
        We observe that {\em $C$ is bounded} since {\em $K_{0}$ is
compact}.

        Let us consider the balls $B (x, \theta (x, \varepsilon) )$ defined
in Lemma~\ref{01A142} with $ \varepsilon := \frac{1}{m}$.
        The compact subset $K_{0}$ can be covered by  $q$ balls $B (x_{i},
\eta (x_{i},\frac{1}{m}))$.  Taking $  \theta :=  \min_{i=1, \ldots ,q}\eta
(x_{i},\frac{1}{m})
>0$, $J$ the smallest integer larger than or equal to $ \frac{T}{ \theta }$
and setting $h:= \frac{T}{J} \leq  \theta $, we infer that
        \begin{displaymath}
         \forall \; x \in K_{0}, \;  \; d \left( f (x), \frac{ \Pi _{K}
(x+hf (x))-x}{h} \right)  \; \leq  \; \frac{1}{m}
        \end{displaymath}
        Starting from $x_{0}$, instead of defining recursively the sequence
of elements $y_{j+1}:=y_{j}+hf (y_{j})$ as in the classical Euler method, 
we define recursively a sequence of elements
\begin{displaymath}
 x_{j+1}  \; :=  \; x_{j}+hu_{j} \; \in  \;\Pi _{K} (x_{j}+hf (x_{j}))
\end{displaymath}
where $f (x_{j})$ is replaced by $u_{j}$:
\begin{displaymath}
 u_{j} \; \in  \;  \frac{\Pi _{K} (x_{j}+hf
(x_{j}))-x_{j}}{h}  \;  \subset   \;  C  \;  \mbox{\rm  satisfies}  \;  \|f
(x_{j})-u_{j}\| \; \leq  \; \frac{1}{m}
\end{displaymath}
for keeping the elements $x_{j} \in K$.
        
The elements $x_{j}$ belong to $K_{0}$, since they belong to $K$ and 
\begin{displaymath} 
\|x_{j}-x_{0}\| \;\leq  \;
\sum_{i=0}^{i=j-1}\|x_{i+1}-x_{i}\|  \;  \leq  \; jh\|C\|  \; \leq 
\;T\|C\| = r
        \end{displaymath}
        whenever $j \leq J$.
        We interpolate the sequence of elements
        $x_{j}$ at the nodes $jh$
         by the piecewise linear functions $x_{m}(t)$ defined on each
interval $[jh, (j+1)h[$ by
        \begin{displaymath}
         \forall \;  t \in [jh, (j+1)h[,  \;\; x_{m}(t)  \;:=  \; x_{j} +
(t - jh)u_{j}
        \end{displaymath}
        We observe that this sequence satisfies the following estimates
        \begin{equation} \left\{ \begin{array}{ll} \label{02A142}
        i) & \forall \;  t \in [0,T], \;\; x_{m}(t) \in \; B (K_{0},
\varepsilon _{m}) \\ 
        & \\
        ii) & \forall \;  t \in [0,T], \;\; \|x'_{m}(t)\| \leq \|C\|
        \end{array} \right. \end{equation}
        Let us fix $t \in [\tau_{m}^{j}, \tau_{m}^{j+1}[$. Since
$\|x_{m}(t) - x_{m}(\tau_{m}^{j})\|$  = $h_{j}\|u_{j}\| \leq \|C\|/m$, and
since $(x_{j},u_{j})$ belongs to $B(\mbox{\rm Graph}(f),\frac{1}{m})$ by
Lemma~\ref{01A142}, we deduce that these functions are approximate
solutions in the sense that
        \begin{equation} \left\{ \begin{array}{ll} \label{03A142}
        i) & \forall \;  t \in [0,T], \;\; x_{m}(t) \in
B(K_{0},\varepsilon_{m})\\
         & \\
        ii) & \forall \;  t \in [0,T], \; (x_{m}(t),x'_{m}(t)) \in B(
\mbox{\rm Graph}(f), \varepsilon_{m})
        \end{array} \right. \end{equation}
        where $\varepsilon_{m} := \frac{\|C\|+1}{m}$ converges to $0$.

        \vspace{ 5 mm}
        
        2. \hspace{ 2 mm} --- \hspace{ 2 mm} {\bf Convergence of the
Approximate Solutions} 
        
        Estimates (\ref{02A142}) imply that for all $t \in [0,T]$, the
sequence $x_{m}(t)$ remains in the compact subset $B(K_{0},1)$ 
         and that the sequence $x_{m}(\cdot)$ is {\em equicontinuous},
because the derivatives $x'_{m}(\cdot)$ are bounded. We then deduce from
Ascoli's Theorem\footnote{ Let us recall that a subset ${\cal H}$ of
continuous functions of ${\cal C}(0,T;{\bf  R}^{n})$ is {\em
equicontinuous}\index{equicontinuous}  if and only if
        \begin{displaymath}
        \forall \;  t \in [0,T],  \; \forall \;  \varepsilon > 0,  \;
\exists \; \eta  \; :=  \; \eta({\cal H},t, \varepsilon )  \; |  \; 
        \forall \;  s \in [t-\eta,t+\eta], \; \sup_{x(\cdot) \in {\cal
H}}\|x(t)-x(s)\| \leq \varepsilon
        \end{displaymath}
        Locally Lipschitz functions with the  same Lipschitz constant form
an equicontinuous set of functions. In particular, a subset of
differentiable functions satisfying
        \begin{displaymath}
        \sup_{t \in [0,T]} \|x'(t)\| \leq c < + \infty
        \end{displaymath} is equicontinuous.
        
        {\em Ascoli's Theorem}\index{Ascoli's Theorem}  states that a
subset ${\cal H}$ of functions  is {\em relatively compact} in ${\cal
C}(0,T;{\bf  R}^{n})$ if and only if it is equicontinuous and satisfies
        \begin{displaymath}
        \forall \;  t \in [0,T],  \; {\cal H}(t) \; :=  \;
\{x(t)\}_{x(\cdot) \in {\cal H}} \; \mbox{\rm is compact.}
        \end{displaymath}} that it remains in a compact subset of the
Banach space ${\cal  C}(0,T;{\bf   R}^{n})$,  and thus,  that a subsequence
(again
denoted) $x_{m}(\cdot)$ converges uniformly to some function $x(\cdot)$. 
        
        \vspace{ 5 mm}
        
        3. \hspace{ 2 mm} --- \hspace{ 2 mm} {\bf The Limit is a Solution}
        
        Condition (\ref{03A142})i) implies that
        \begin{displaymath}
        \forall \;  t \in [0,T],  \;\; x(t) \in K_{0}
        \end{displaymath}
        i.e., that $x(\cdot)$ is viable.
        
        Property (\ref{03A142})ii) implies that for almost every $t \in
[0,T]$,  there exist $u_{m}$ and $v_{m}$ converging to $0$ such that
        \begin{displaymath}
         x_{m}' (t) \; = \; f (x_{m} (t)-u_{m}) +v_{m}
        \end{displaymath}
        We thus deduce that for almost $t \geq 0$, $x_{m}' (t)$ converges
to $f(x (t))$. On the other hand, 
        \begin{displaymath}
         x_{m} (t) -x_{m} (s) \; =  \int_{s}^{t}x'_{m} ( \tau )d \tau 
        \end{displaymath}
        implies that $x'_{m} (t)$ converges almost every where to $x' (t)$.
We thus infer that $x ( \cdot )$ is a solution to the differential
equation.
        $ \; \; \Box $

        \vspace{ 5 mm}
        
        {\bf Proof of Theorem~\ref{03A43}}  --- \hspace{ 2 mm}
         First, $K$ is locally compact since it is closed and the dimension
of ${\bf  R}^{n}$ is finite.
        
        Second, we claim that starting from any $x_{0}$, there exists a
maximal solution. Indeed, denote by ${\cal S}_{[0,T[}( x_{0})$ the set of
solutions to the differential equation defined on $[0,T[$.
        
        We introduce the set
        of pairs $ \{(T,x(\cdot ))\}_{T>0, \; x(\cdot ) \in {\cal
S}_{[0,T[}(x_{0})}$  on which we consider the order relation $ \prec $
defined by
        \begin{displaymath}
         (T,x(\cdot ))  \prec  (S, y(\cdot ))  \; \mbox{if and only if} \;
T \; \leq  \; S  \;\; \& \;\;\forall \; t \in [0,T[, \;  \; x(t)  =  y(t)
        \end{displaymath}
        
        Since every totally ordered subset has obviously a majorant, Zorn's
Lemma implies that any solution $y(\cdot ) \in {\cal S}_{[0,S[}(x_{0})$
defined on some interval $[0,S[$ can be extended to a solution $x(\cdot )
\in {\cal S}_{[0,T[}(x_{0})$  defined on a maximal interval $[0,T[$.
        
        Third, we have to prove that if $T$ is finite, we cannot have
        \begin{displaymath}
        c  \; := \; \limsup_{t \rightarrow T-}\|x(t)\| \; < \; +\infty 
\end{displaymath}
        Indeed, if $c <+\infty$, there would exist a constant $\eta \in
]0,T[$ such that
        \begin{displaymath}
        \forall \; t  \in [T-\eta,T[, \;\;\; \|x(t)\|  \;\leq  \; c+1
        \end{displaymath}
        Since $f$ is continuous images on the compact subset $K \cap
(c+1)B$, we infer that there exists a constant $ \rho $ such that 
for all $  s  \in [T-\eta,T[$, $ \|f(x (s))\| \leq  \rho $.

        Therefore, for all $\tau, \sigma \in [T-\eta,T[$, we obtain:
        \begin{displaymath}
         \|x(\tau)-x(\sigma)\| \;  \leq \;  \int_{\sigma}^{\tau}\|f(x(s))\|
ds
\; \leq \; \rho |\tau - \sigma|
        \end{displaymath}
        Hence the Cauchy criterion implies that $x(t)$ has a limit when $t
\rightarrow T-$. We denote by $x(T)$ this limit, which belongs to $K$
because it is closed.
        Equalities
        \begin{displaymath}
        x(T_{k})  \; = \; x_{0} + \int_{0}^{T_{k}} f(x(\tau))d\tau 
        \end{displaymath}
 imply that by letting $k \rightarrow
\infty$, we obtain:
        \begin{displaymath}
        x(T) \; =  \; x_{0} + \int_{0}^{T}f(x(\tau))d\tau  
        \end{displaymath}
        This means that we can extend the solution up to $T$ and even
beyond, since  Theorem~\ref{01A43} allows us to find a viable solution
starting at $x(T)$ on some interval $[T,S]$ where $S>T$.
        Hence $c$ cannot be finite.
         $  \;  \; \Box$
        
        \vspace{ 5 mm}
        
        {\bf Proof of Theorem~\ref{04A43}} --- \hspace{ 2 mm} Since the
growth of $f$ is linear, \begin{displaymath}
        \exists \; c \geq 0,   \;\; \mbox{ such that} \;\; \forall \; x \in
{\bf  R}^{n}, \;\; \|f(x)\| \leq c(\|x\|+1)
        \end{displaymath}
        Therefore, any solution to differential equation (\ref{01A11})
satisfies the estimate:
        \begin{displaymath}
        \|x'(t)\|  \;  \leq  \;  c(\|x(t)\|+1)
        \end{displaymath}
        The function $t \rightarrow \|x(t)\|$ being locally Lipschitz, it
is almost everywhere differentiable. Therefore, for any $t$ where $x(t)$ is
different from $0$ and differentiable, we have
        \begin{displaymath}
        \frac{d}{dt}\|x(t)\|  \;  = \;  \left\langle  
\frac{x(t)}{\|x(t)\|},x'(t) \right\rangle   \;  \leq  \;  \|x'(t)\|
        \end{displaymath}
        These two inequalities imply the estimates:
         \begin{equation} \label{xA445}
         \|x(t)\|  \; \leq \; ( \|x_{0}\|+1)e^{ct}    \;  \; \& \;  \; 
        \|x'(t)\|  \; \leq  \;  c( \|x_{0}\|+1)e^{ct}
         \end{equation}
        Hence, for any $T > 0$, we infer that \[ \limsup_{t \rightarrow
T-}\|x(t)\| \; < \; + \infty \] Theorem~\ref{03A43} implies that we can
extend the solution on the interval $[0,\infty[$.
         $ \;  \;  \Box$ 
        
        \section{The Solution Map}

        \begin{Definition}[Solution Map]  \label{02A44}
        We denote by ${\cal S}_{f}(x_{0})$ \label{01A44} the  set of
solutions to  differential equation (\ref{01A11}) and call the set-valued
map ${\cal S}_{f}: x \leadsto {\cal S}_{f} (x)$  the {\em solution
map\/}\index{solution map} of $f$ (or of  differential inclusion 
(\ref{01A11}).)
         \end{Definition}

        \begin{Theorem} \index{solution map}  \label{03A44}
Let  us  consider a finite dimensional vector space $ {\bf  R}^{n}$ and
   a {\em  continuous  map $f  :{\bf   R}^{n} \mapsto {\bf  R}^{n}$\/} with
linear growth. 
Then  the graph of the restriction of ${\cal S}_{f}|_{L}$ to any compact
subset  $L$   is compact in ${\bf   R}^{n} \times {\cal C}(0,  \infty ;{\bf
R}^{n})$ where the space
${\cal C}(0,\infty ;{\bf  R}^{n})$ is supplied with the compact convergence
topology.
 \end{Theorem}
        
        {\bf Proof}  --- \hspace{ 2 mm}
        We shall show that the graph of the restriction ${\cal S}_{f}|_{L}$
of the solution map ${\cal S}_{f}$ to a compact subset $ L $  is compact. 
                    Let     us     choose     a     sequence    of elements
$(x_{0_{n}},x_{n}(\cdot))$ of
the graph of the solution map ${\cal S}_{f}$. They satisfy:
        \begin{displaymath} 
         x'_{n}(t)   =f (x_{n}(t))  \;  \; \& \;  \;  x_{n}(0)   = 
x_{0_{n}} \in L
         \end{displaymath} 
         A subsequence (again denoted) $x_{0_{n}}$ converges to some $x_{0}
\in L$ because $L$ is compact. By Theorem~\ref{04A43}, 
        \begin{displaymath}
        \forall \; n \geq 0, \;  \;  \|x_{n}(t)\| \; \leq \; (
\|x_{0_{n}}\|+1) e^{ct}  \; \; \& \; \; \|x'_{n}(t)\| \; \leq \; c(
\|x_{0_{n}}\|+1) e^{ct}
        \end{displaymath}
        Thus, by Ascoli's Theorem, the sequence $x_{n}(\cdot)$ is
relatively compact in ${\cal C}(0,\infty;{\bf  R}^{n})$. 
        We deduce from this,
	that a subsequence (again denoted) $x_{n}  ( \cdot)$ converges
	to a continuous function $x ( \cdot )$ uniformly  on compact
intervals. Therefore, passing to the limit in equalities
        \begin{displaymath}
         x_{m} (t) \; = \; x_{0_{n}} + \int_{0}^{t}f (x_{m} ( \tau ))d \tau
\end{displaymath}
        we deduce that $x ( \cdot )$ is a solution to the differential
equation starting at $x_{0}$. $ \;  \;  \Box$ 
        
        \section{Uniqueness Criteria}
        
        Whenever there exists a unique solution to differential equation
$x'=f (x)$ starting from any initial state $x_{0}$, then viability and
invariance properties of a closed subset $K$ are naturally equivalent. This
is one of the motivations for providing uniqueness criteria.
        
        \begin{Definition}
         We shall say that a map $f:{\bf   R}^{n}  \mapsto {\bf  R}^{n}$ is
{\em monotone\/}
\index{monotone} if there exists 
 there exists $ \mu  \in  {\bf R}$ such that 
\begin{equation} \label{monassg}
  \langle f(x_{1}) - f (x_{2}) ,x_{1}-x_{2} \rangle \; \leq  \; -\mu
\|x_{1}-x_{2}\|^{2}
\end{equation}

         \end{Definition}
        The interesting case is obtained when $   \mu >0$. When $f$ is
Lipschitz  with constant $ \lambda $, then it is monotone with $ \mu =-
\lambda $.
        
 \begin{Theorem} \label{05A13m}
Let us consider a subset $K $  of a finite dimensional vector  space $ {\bf
R}^{n}$ and
a  continuous and monotone map $f $ from ${\bf  R}^{n}$ to ${\bf  R}^{n} $.
The solution to
differential equation $x'=f (x)$ starting from $x_{0}$ is unique. If $x_{i}
( \cdot )$ are two solutions to the differential equation $x'=f (x)$,then
\begin{displaymath}
 \|x_{1} (t)-x_{2} (t)\| \; \leq  \; e^{- \mu t} \|x_{1} (0)-x_{2} (0)\|
\end{displaymath}
  \end{Theorem} 
        
        {\bf Proof} --- \hspace{ 2 mm}
        Indeed, integrating the two sides of  inequality
\begin{displaymath}
 \frac{d}{dt} \|x_{1} (t)-x_{2} (t)\|^{2} \; = \; 2  \langle f(x _{1}
(t))-f (x_{2} (t)),x_{1} (t)-x_{2} (t) \rangle  \; \leq  \; - 2 \mu \|x_{1}
(t)-x_{2} (t)\|^{2}
\end{displaymath}
yields 
\begin{displaymath}
 \|x_{1} (t)-x_{2} (t)\|^{2} \; \leq  \; e^{-2 \mu t} \|x_{1} (0)-x_{2}
(0)\|^{2}
\end{displaymath}

        \section{Backward Viability}
        
        \begin{Definition}
        The subset $K$ is {\sf locally backward viable} 
\index{locally backward viable} under $ f$ if for any $x\in K $,
for any $t>0$, there exist $s \in [0,t[$  and a solution $x (
\cdot )$ to differential  equation (\ref{01A11}) such that
\begin{displaymath}
 \forall \;  \tau \in [s,t], \; x( \tau ) \in K \;   \& \; x (t)=x
\end{displaymath}
It is (globally) backward viable if we can take $s=0$ in the above
statement, and locally (resp. globally) backward invariant if 
for any $x\in K $, for any $t>0$, for all solutions $x ( \cdot )$ to
differential  equation (\ref{01A11}), there exist $s \in [0,t[$ such that
\begin{displaymath}
 \forall \;  \tau \in [s,t], \; x( \tau ) \in K  \; \& \; x (t)=x
\end{displaymath}

         \end{Definition}
        We now compare the invariance of a subset and the backward 
invariance of its complement:
        \begin{Lemma}
         \label{bacwinvcomp}
        A subset $K$ is invariant under a map $ f$ if and only if its
complement $K^{c}:=  {\bf  R}^{n} \backslash K$ is backward invariant under
$ f $.
        \end{Lemma}
        {\bf Proof} --- \hspace{ 2 mm}
        To say that $K$ is not invariant under $f$ amounts to saying that
        there exist a solution $x ( \cdot )$ to differential equation
(\ref{01A11}) and $T>0$ such that
        \begin{displaymath}
          x (0)  \;  \in   \;  K  \;  \&  \;  x (T) \; \in  \; {\bf  R}^{n}
\backslash K
        \end{displaymath}
         and to  say  that  ${\bf   R}^{n}  \backslash  K$  is not backward
invariant amounts
to saying that there exist a solution $y( \cdot )$ to differential equation
(\ref{01A11}), $T >0$ and $S \in [0,T[$ such that
        \begin{displaymath}
          y  (S)  \;  \in   \;  K  \;  \&  \; y (T) \; \in  \; {\bf  R}^{n}
\backslash K
        \end{displaymath}
        It is obvious that the first statement implies the second one by
taking $y ( \cdot )=x ( \cdot )$ and $S=0$. Conversely, the second
statement  implies the first one by taking $x (t):= y (t+S)$ and replacing
$T$ by $T-S >0$ since $x (0)=y (S)$ belongs to $K$ and $x (T-S)=y (T)$
belongs to ${\bf  R}^{n} \backslash K$. $\; \; \Box$

        \mbox{}
        
        It is also useful to relate backward viability and
invariance under $ f$ to viability and invariance under $ - f$:
        \begin{Lemma} \label{backinvab}
        Let us assume that $f$ is continuous with linear growth.
        Then $K$ is locally backward viable (resp. invariant)
under $ f$ if and only if $f$ is locally viable (resp. invariant)
under $ - f $.
        \end{Lemma}
        {\bf Proof} --- \hspace{ 2 mm}
        Let us check this statement for local viability.
        Assume that $K$ is locally backward viable and infer that
$K$ is locally invariant under $ - f$. Indeed, let $x \in K$.
Then, for any $T >0$, there exists $S \in [0,T[$ and  a solution
$ x( \cdot ) $ to differential equation (\ref{01A11}) viable in
$K$ on the interval $[S,T]$ and satisfying $x (T)=x$. Let $y (
\cdot ) $ be a solution to the differential equation $y'=-f (y)$ 
starting at $y (0)=x (S)$. Then the  function $ z ( \cdot )$
defined by
        \begin{displaymath}
         z (t) \; = \left\{ \begin{array}{lll}
        x (T-t) & \mbox{\rm if} & t \in [0,T-S] \\
        y (t+T-S) & \mbox{\rm if} & t \geq T-S
        \end{array} \right.
        \end{displaymath}
        is a solution to the differential equation $z'=-f (z)$  starting at
$z (0)=x (T)=x$ and viable in $K$ on the interval $[0,T-S]$.
        
        Conversely, assume that $K$ is locally viable under $ -f
$ and check that $K$ is locally backward invariant. Let $x \in
K$, $T >0$ and one solution $x ( \cdot )$ to differential
equation $x'=-f (x)$ viable in $K$ on $[0,R]$ where $R>0$. Let be
any  solution $y ( \cdot )$ to $y' (t)=f (y (t))$ starting at $x
$ and  set
        \begin{displaymath}
         z (t) \; = \left\{ \begin{array}{lll}
        x (T-t) & \mbox{\rm if} & t \in [0,T] \\
        y (t-T) & \mbox{\rm if} & t \geq T
        \end{array} \right.
        \end{displaymath}
        Hence $ z ( \cdot )$ to differential equation (\ref{01A11})
satisfying $x (T) =x \in K$ and viable in $K$ on the interval $[S,T]$ where
$S := \max (T-R,0)$.
        $\; \; \Box$ 
        
        \section{Time-Dependent Differential Equations}

 \begin{Theorem} \label{04A43t}
Let us consider a subset $K $  of a finite dimensional vector space  $ {\bf
R}^{n}$ and
a   map $f $  from $  {\bf  R}_{+} \times K $ to ${\bf  R}^{n} $. We assume
that the map
$f $ is {\em continuous} from $ {\bf  R}_{+}  \times K$ to ${\bf  R}^{n} $,
that it
has {\em uniform linear growth\/} in the sense that
\begin{displaymath}
\exists \; c > 0   \;\; \mbox{ such that} \;\;\forall \;  t \geq 0, \; x
\in K,\;\;  \|f(t,x ) \| \; \leq \; c(\| x\|+1  )   
\end{displaymath}
If
\begin{displaymath}
 \forall \; t \geq 0, \forall \;  x \in K , \;  \; f (t,x) \; \in  \; T_{K}
(x)
\end{displaymath}
then $K$ is globally viable under $f$: for every  initial state $x_{0}  \in
K $, there exists a {\em viable solution on} $[0,\infty] $ to differential
equation
\begin{displaymath}
 x' (t) \; := \; f (t,x (t))
\end{displaymath}
starting at $x_{0}$ and satisfying
\begin{displaymath}
 \forall \; t \geq 0, \;  \; \|x (t)\| \; \leq  \;  \|x_{0}\| e^{ ct}
+e^{ct}-1
\end{displaymath}
Assume moreover that $f$ is {\em uniformly monotone\/} \index{uniformly
monotone} in the sense that\footnote{The interesting case is obtained when
$   \mu >0$. When $f$ is uniformly Lipschitz  with constant $ \lambda $,
then it is uniformly monotone with $ \mu =- \lambda $.} there exists $ \mu 
\in  {\bf R}$ such that 
\begin{equation} \label{monassgt}
  \langle f(t,x_{1}) - f (t,x_{2}) ,x_{1}-x_{2} \rangle \; \leq  \; -\mu
\|x_{1}-x_{2}\|^{2}
\end{equation}
If $x_{1}$ and $x_{2}$ are two initial states, then the solutions $x_{i} (
\cdot )$ starting from $x_{i}$, $ (i=1,2)$, satisfy 
\begin{displaymath}
 \|x_{1} (t)-x_{2} (t)\| \; \leq  \; e^{- \mu t} \|x_{1} (0)-x_{2} (0)\|
\end{displaymath}
  \end{Theorem} 
        
        {\bf Proof} --- \hspace{ 2 mm}
        We deduce the first statement from the standard trick which amounts
to observing that a solution $x ( \cdot )$ to the time-dependent
differential equation $x'=f (t,x)$ starting at time $0$ from the initial
state $x_{0}$ if and only if $ ( \tau  ( \cdot ),x ( \cdot ))$ is a
solution to the system of differential equations
        \begin{displaymath} \left\{ \begin{array}{ll}
        i) &  \tau ' (t) \; = \; 1 \\
        ii) & x' (t) \; = \; f ( \tau (t),x (t))
        \end{array} \right. \end{displaymath}
        starting at time $0$ from $ (0,x_{0})$. The solution $x ( \cdot )$
is viable in $K$ under $f$ if and only if $ ( \tau  ( \cdot ),x ( \cdot ))$
is viable in $ {\bf  R}_{+}  \times K$. By the Nagumo Theorem, this is
equivalent to require that $ (1,f (t,x)) \in T_{ {\bf  R}_{+} \times K}
(t,x)$, i.e., that $f (t,x)$ belongs to $T_{K} (x)$.
        
        The proof of the second statement is same than the one of
Theorem~\ref{05A13m}.
        $\; \; \Box$

\chapter{Viability Kernels and Capture Basins}

\vspace{ 22 mm}
{\Huge \bf Introduction}

\vspace{ 13 mm}

When a closed subset $K$ is not viable under a dynamical economy, then two
questions arise naturally: 
\begin{enumerate}
\item find solutions starting from $K$ which {\em remain viable in $K$ as
long as possible, hopefully, forever},
\item and starting outside of $K$, find solutions which {\em return to $K$
as soon as possible, hopefully, in finite time}
\end{enumerate} 

Studying these questions leads to the concepts of  
\begin{enumerate}
\item {\em viability kernel\/} of a subset $K$ under a dynamical system, as
the set of elements of $K$ from which starts  a solution viable in $K$,
\item {\em capture basin\/} of $C$, which is the set  of points of $K$ from
which a solution reaches $C$ in finite time,
\item when $C \subset K$, {\em viable-capture basin\/}, which is the subset
of points of $K$ from which starts a solution reaching $C$ before leaving
$K$.
\end{enumerate}

We  shall provide characterizations of these concepts and derive their
properties we shall use later for solving some Hamilton-Jacobi equations
and boundary-value problems for systems of first-order partial differential
equations. For instance, if $K$ is backward invariant and a repeller, the
capture basin of $C$ is the unique closed subset $D$ satisfying
\begin{displaymath} \left\{ \begin{array}{ll}
i) & C \; \subset  \; D \; \subset  \; K\\
ii) & \forall \; x \in D \backslash C, \; \;f(x) \; \in  \;T_{D} (x)
\\
iii) &  \forall \; x \in D, \;  \; -f(x) \; \in  \; T_{D} (x)\\
\end{array} \right. \end{displaymath}
or, equivalently, by duality, the ``normal conditions''
\begin{displaymath} \left\{ \begin{array}{ll}
i) & C \; \subset  \; D \; \subset  \; K\\
ii) & \forall \; x \in D \backslash C, \;  \forall \; p \in N_{D} (x),
\;  \;  \langle p,f (x) \rangle  \; = \; 0 \\
iii) &  \forall \; x \in D, \; \forall \; p \in N_{D} (x), \;  \; 
\langle p,f (x) \rangle  \; \geq  \; 0\\
\end{array} \right. \end{displaymath}

\section{Reachable, Viability and Capture Tubes} 

\begin{Definition}
Let  $f:{\bf   R}^{n} \mapsto  {\bf   R}^{n}$  be a map and $C \subset {\bf
R}^{n}$ be a subset. 
        The {\sf reachable map} \index{reachable map}
        $\vartheta_{f} ( \cdot ,x)$ is defined by
\begin{displaymath}
\forall \; x \in {\bf  R}^{n}, \;  \; \vartheta_{f} (t,x) \; := \;  \left\{
x(t)
\right\}_{ x(\cdot ) \in{\cal
S}_{f}(x)}
\end{displaymath}
        We associate with it the {\sf reachable tube} $ t \leadsto
\vartheta_{f} ( t,C)$ defined by \index{reachable tube}
\begin{displaymath}
 \vartheta_{f} (t,C) \; := \;  \left\{ x(t) \right\}_{ x(\cdot ) \in{\cal
S}_{f}(C)}
\end{displaymath}
        
 \end{Definition}
We derive the following properties:

\begin{Proposition}
The reachable map $t\leadsto \vartheta_{f}(t,x)$ enjoys the semi-group
property:
$\forall t,s\geq 0, \;  \;\vartheta_{f}(t+s,x)=\vartheta_{f}
(t,\vartheta_{f} (s,x))$.

Furthermore,
\begin{displaymath}
 ( \vartheta _{f } (t, \cdot ))^{-1} \; := \; \vartheta _{-f } (t, \cdot )
\end{displaymath}
Therefore, the subset $ \vartheta _{-f} (t,C)$ is the subset of elements $x
\in {\bf  R}^{n}$ which reach the subset at the prescribed time $t$.

If $f$  is continuous  with linear growth  and $K \subset  {\bf  R}^{n}$ is
closed, the
graph of the reachable map $ t \leadsto \vartheta_{f} ( t,K)$ is
closed.
\end{Proposition} 

{\bf Proof} --- \hspace{ 2 mm}
The semi-group property is obvious.
Let us prove the second one: If $ y \in  \vartheta_{f } (t,x)$, there
exists a solution $x ( \cdot )$ to the differential equation $x'=  f(x)$
starting at $x$ such that $y = x (t)$. We set $y (s) := x (t-s)$ if $ s \in
[0,t]$ and we choose any solution $y ( \cdot )$ to the differential
equation $y' \in - f(y)$ starting at $x$ at time $t$ for $s \geq t$. Then
such a function $y ( \cdot )$ is a solution to the differential equation
$y' \in  - f(y)$ starting at $y$ and satisfying $y (t)=x$. This shows that
$x \in \vartheta _{-f } (t,y)$.

The last statement is a consequence of Theorem~\ref{04A43}.
$\; \; \Box$ 

        \begin{Definition}
         Let  $f:{\bf  R}^{n} \mapsto {\bf  R}^{n}$ be a map and $C \subset
{\bf  R}^{n}$ be any subset. 
        \begin{enumerate}
        \item The subset $\mbox{\rm Viab}_{f}(C,T)$ of initial states
$x_{0} \in C$ such that one solution $x(\cdot )$ to differential equation
$x'= f(x)$ starting at $x_{0}$ is viable in $C$ for all $t \in  [0,T]$ is
called the {\em $T$-viability kernel\/} \index{viability kernel} and the
subset $\mbox{\rm Viab}_{f}(C) := \mbox{\rm Viab}_{f}(C, \infty )$ is
called the {\sf viability kernel} \index{viability kernel} of $C$ under
$f$. A subset $C$ is a repeller if its viability kernel is empty.
\index{repeller}

        \item The subset $ \mbox{\rm Capt}_{f}(C,T)$ of initial states
$x_{0} \in {\bf   R}^{n}$  such  that  $C$  is reached  before  $T$  by one
solution $x(\cdot
)$ to differential equation $x'= f(x)$ starting at $x_{0}$ is called the
{\em $T$-capture basin\/} \index{capture basin}  and
\begin{displaymath}
 \mbox{\rm Capt}_{f}(C) \; := \;  \bigcup_{T > 0}^{} \mbox{\rm
Capt}_{f}(C,T)
\end{displaymath}
is said to be the {\sf capture basin} \index{capture basin} of $C$. 
\item When $C \subset K$, the {\sf viable-capture basin}
\index{viable-capture basin} $ \mbox{\rm Capt}_{f}^{K}(C)$
of  $C$ in $K$ by $f$  is the set of initial states $x_{0} \in
{\bf   R}^{n}$  from which starts at least one solution to the $x'=  f (x)$
viable in
$K$ until it reaches $C$ in finite time.
        
        \end{enumerate}
 \end{Definition}

{\bf Remarks} --- \hspace{ 2 mm}
We observe that if $T_{1} \leq T_{2}$,

 \begin{displaymath} \left\{ \begin{array}{ccccc}
 \mbox{\rm Viab}_{f}(C)    & \subset  & \mbox{\rm Viab}_{f}(C,T_{2})  &
\subset  & \mbox{\rm Viab}_{f}(C,T_{1})     \\
 \cap & & \cap & & \cap \\
 \mbox{\rm Capt}_{f}(C,0)    & = & C  & = & \mbox{\rm Viab}_{f}(C,0)     \\
 \cap & & \cap & & \cap \\
 \mbox{\rm Capt}_{f}(C)    & \supset  & \mbox{\rm Capt}_{f}(C,T_{2})  &
\supset  & \mbox{\rm Capt}_{f}(C,T_{1})     \\
 \end{array} \right. \end{displaymath} 

One can write
\begin{displaymath}
 \mbox{\rm Capt}_{f}(C,T) \; = \; \bigcup_{t \in [0,T]}^{} \vartheta_{-f}
(t,C)
\end{displaymath}
        
We point out the following obvious properties:
\begin{Lemma}
Let  $C  \subset  K$  be  closed  subsets.  The  capture  basin  $\mbox{\rm
Capt}_{f}(C)$ is smallest backward invariant containing $C$ and
$\mbox{\rm Capt}_{f}(C) \backslash C$ is locally viable. The capture basin
of any union of subsets $C_{i} \; (i \in I)$ is the union of the capture
basins of the $C_{i}$.
When $C \subset K$ where $K$ is assumed to be backward invariant, then the
viable-capture basin satisfies
\begin{displaymath}
\mbox{\rm Capt}_{f}^{K}(C) \; = \;  \mbox{\rm Capt}_{f}(C) \; \subset  \; K
\end{displaymath}

The complement $\mbox{\rm Capt}_{f}^{K}(C) \backslash C$ of $C$ in the
viable-capture basin  $\mbox{\rm Capt}_{f}^{K}(C)$ is locally viable.

\end{Lemma}
{\bf Proof} --- \hspace{ 2 mm}
Indeed,  whenever $K$ is backward invariant, each backward reachable set 
$\vartheta_{-f} (t,C) $ is contained in $K$, so that
\begin{displaymath}
  \mbox{\rm Capt}_{f}(C) \; = \; \bigcup_{t \geq 0}^{} \vartheta_{-f} (t,C)
\; \subset  \; K
\end{displaymath}
Since the intersection of
backward invariant subsets is backward invariant, the capture basin is
contained in the smallest backward invariant subset containing $C$. The
semi group property implies that the capture basin, which is the union of
the backward reachable subsets, is backward invariant.

When $C \subset K$ where $K$ is assumed to be backward invariant, then the
capture basin 
\begin{displaymath}
 \bigcup_{t \geq 0}^{} \vartheta_{-f} (t,C)
\end{displaymath}
is contained in $K$.

If $x$ belongs to $\mbox{\rm Capt}_{f}^{K}(C) \backslash C$, then there
exists a solution to the differential equation $x'=f (x)$ starting from $x$
which reaches $C$ before leaving $K$, and thus, which is viable in
$\mbox{\rm Capt}_{f}^{K}(C) \backslash C$ on some nonempty interval. $\; \;
\Box$ 

\mbox{}
        
        \begin{Proposition} \label{largestvaibinv}
The viability kernel $\mbox{\rm Viab}_{f}(C)$  of $C$ under $f$ is the
largest subset of $C$ viable  under $f$.

Furthermore,  $C  \backslash  \mbox{\rm  Viab}_{f}(C)$  is  a  repeller and
$\mbox{\rm Viab}_{f}(C) \backslash \partial C$ is locally
backward invariant.
        \end{Proposition} 
        {\bf Proof} --- \hspace{ 2 mm}
Every subset $L \subset C$ viable  under $f$ is obviously contained in the
viability kernel $\mbox{\rm Viab}_{f}(C)$  of $C$ under $f$.
        
        On the other hand, if  $x(\cdot)$ is a solution to the differential
equation $x'= f(x)$ viable in $C$, then for all $t >0$, the function
$y(\cdot)$ defined by $y(\tau) := x(t+ \tau)$ is also a solution to the
differential equation, starting at $x(t)$, viable  in $C$.
        
        Therefore, for  any element $x_{0} \in \mbox{\rm Viab}_{f}(C)$,
there exists  a viable solution $x(\cdot)$ to the differential equation
starting from $x_{0}$, and thus, for all $t \geq 0$,  $ x(t) \in \mbox{\rm
Viab}_{f}(C)$, so that it is viable under $f$.

Let us assume that $\mbox{\rm Viab}_{f}(C) \backslash \partial C$ is not
locally backward invariant: There would exist $x \in  \mbox{\rm
Viab}_{f}(C) \backslash \partial C$, $T >0$ and a solution $x ( \cdot )$ to
the differential equation $x'=f (x)$ satisfying $x (T)=x$ such that for all
$S <T$, there  exist $S ' \in [S,T]$ such that $x (S')$  belong to the
union of $ \partial C$ and of the complement of the viability kernel
$\mbox{\rm Viab}_{f}(C)$. Since $x (T)=x$ does not belong to the boundary $
\partial C$ of $C$, we know that for $S$ close enough to $T$, $x ([S,T])
\cap  \partial C = \emptyset $. Hence $x (S')$ does not belong to the
boundary of $C$, so that it belongs to complement of the viability kernel
$\mbox{\rm Viab}_{f}(C)$, and thus, the solution $x ( \cdot )$
starting  from  $x  (S')  \in \mbox{\rm Viab}_{f}(C)$  at time $S' $ should
leave $C$ in finite time, a
contradiction.
$\; \; \Box$

\section{Hitting and Exit Times}
\begin{Definition}
     We say that 
the {\sf hitting functional\/} (or {\sf minimal time functional})
\index{minimal time functional}  associating with $x (\cdot) $ its {\sf
hitting time\/} \index{hitting time} $\omega_{C}(x(\cdot ))$ is defined
by   \index{hitting  functional}
\begin{displaymath}
\omega_{C}(x(\cdot )) \; := \; \inf \left\{ t \in [0,+\infty [ \; | \; x(t)
\in C  \right\}
\end{displaymath}
and the function  $\omega_{C}^{f ^{\flat }}: C \mapsto {\bf R}_{+} \cup 
\{+\infty \}$  defined by 
\begin{displaymath}
\omega_{C}^{f ^{\flat }}(x) \; := \; \inf_{ x(\cdot ) \in {\cal
S}_{f}(x)}\omega_{C}(x(\cdot )) 
\end{displaymath}
is  called   the  (lower)  {\em hitting   function\/} or  {\em minimal time
function\/}. \index{minimal time function}
In the same way, the {\sf exit functional} is defined by
\begin{displaymath}
\tau _{C} (x ( \cdot )) \; := \; \omega_{{\bf  R}^{n} \backslash C}(x(\cdot
)) \; :=
\; \inf \left\{ t \in [0,+\infty [ \; | \; x(t) \notin C  \right\}
\end{displaymath}
and  the function  $ \tau _{C}^{f ^{\sharp }}:C \mapsto {\bf R}_{+} \cup 
\{+\infty \}$  defined by  
\begin{displaymath}
\tau _{C}^{f ^{\sharp }}(x) \; := \;  \sup _{ x(\cdot ) \in {\cal
S}_{f}(x)}\tau _{C}(x(\cdot ))
\end{displaymath}
is called the (upper)  {\em exit  function\/}. \index{exit  function}

Let $C \subset K \subset  {\bf  R}^{n}$ be  two closed subsets. 
We also introduce the function
\begin{displaymath}
  \gamma _{ (K,C)}^{f ^{\flat }} (x) \; := \; \inf_{ x ( \cdot ) \in 
{\cal S}_{f} (x)} ( \omega _{C} (x ( \cdot ))- \tau _{K} (x ( \cdot
)))
\end{displaymath}
 \end{Definition}

We observe that if $C=K$, $\omega _{C} (x ( \cdot )) =0$ for all
solutions $x ( \cdot ) \in {\cal S}_{f} (x)$ starting from $K$ and
thus, that
\begin{displaymath}
  \gamma _{ (K,K)}^{f ^{\flat }} (x) \; = \; - \tau _{K}^{f ^{\sharp
}} (x)
\end{displaymath}

To say that $K \backslash C$ is a repeller amounts to saying that for
every solution $ x ( \cdot ) \in  {\cal S}_{f} (x)$ starting from $x
\in K \backslash C$, $\min ( \omega _{C} (x ( \cdot ), \tau _{K} (x (
\cdot )))) < + \infty $ and to say that $K$ is a  repeller under $f$
amounts to saying that the exit function $ \tau _{K}^{f ^{\sharp }}$
is finite on $K$.

\begin{Lemma} \label{Picard02}
When $K$ is closed, the exit functional  $ \tau _{K}$ is upper
semicontinuous  and   the   hitting   functional   $\omega_{K}$   is  lower
semicontinuous 
when ${\cal  C}(0,\infty ;{\bf   R}^{n})$   is  supplied  with  the compact
convergence
topology. 
 \end{Lemma}

{\bf Proof} --- \hspace{ 2 mm}
For proving that the exit functional $ \tau_{K}$ is upper semicontinuous,
we shall check that the subsets $ \{x ( \cdot ) \; | \; \tau _{K} ( x(\cdot
)) < T \}$ are open for the pointwise convergence, and thus, the compact
convergence. Let $x_{0} ( \cdot )$ belong to such a set when it is not
empty. Since $x_{0}(T)$ does not belong $K$ which is closed, there exist $
\alpha >0$ such that $B (x_{0}(T), \alpha ) \cap K = \emptyset $. Then the
set of continuous functions $x ( \cdot )$ such that $x (T) \in
\stackrel{\circ}{B} (x_{0}(T), \alpha )$ is open and satisfy $\tau _{K} (
x(\cdot )) < T$.

For proving that the hitting  functional is lower semicontinuous,  we shall
check
that the subsets $ \{x ( \cdot ) \; | \; \omega _{K} ( x(\cdot )) \leq  T
\}$ are closed. Let $x_{n} ( \cdot )$ satisfying $\omega _{K} ( x_{n}(\cdot
))
\leq  T $ converge to $x ( \cdot )$ uniformly on compact intervals. For any
$ \varepsilon >0$, one can find $t_{n} \leq T+ \varepsilon $ such that
$x_{n} (t_{n})$ belongs to $K$. A subsequence (again denoted by) $t_{n}$
converges to some $t$. Since $x_{n} ( \cdot )$ converges uniformly to $x (
\cdot )$ on $[t - \varepsilon ,T + \varepsilon ]$, we deduce that $x (t)$
is the limit of $x_{n} (t_{n}) \in K$ and thus, that $x (t)$ belongs to the
closed subset $K$ and thus, that $ \omega _{K} (x ( \cdot )) \leq T +
\varepsilon $. Letting $ \varepsilon $ converge to $0$, we infer that $
\omega  _{K} (x ( \cdot )) \leq T$. $\; \; \Box$ 

\mbox{}

We deduce the following properties of these hitting and exit functions:

\begin{Proposition} \label{Picard03se}
Let $f:{\bf  R}^{n} \mapsto {\bf  R}^{n}$ be a   continuous map with linear
growth and $C \subset
K \subset {\bf   R}^{n}$  be  two closed  subsets.  Assume  that  $K$  is a
repeller.

The function $\gamma _{ (K,C)}^{f ^{\flat }}$ is lower semicontinuous
and for any $x \in K$, there exists a solution $x_{ (K,C)} ( \cdot )
\in  {\cal S}_{f} (x)$ satisfying 
\begin{displaymath}
  \gamma _{ (K,C)}^{f ^{\flat }} (x) \; := \;  \omega _{C} (x_{
(K,C)}( \cdot ))- \tau _{K} (x _{ (K,C)}( \cdot )))
\end{displaymath}

In the same way, the  hitting function $\omega_{K}^{f ^{\flat }}$ is lower
semicontinuous and the  exit function $\tau _{K}^{f ^{\sharp }}$ is upper
semicontinuous. Furthermore, for any $x \in  \mbox{\rm Dom}(\omega_{K}^{f
^{\flat }})$, there exists one solution $x ^{\flat } ( \cdot ) \in  {\cal
S}_{f} (x)$ which hits $K$ as soon as possible
\begin{displaymath}
 \omega_{K}^{f ^{\flat }} (x) \; = \;  \omega_{K}  (x ^{\flat } ( \cdot ))
\end{displaymath}
and for any $x \in  \mbox{\rm Dom}(\tau _{K}^{f ^{\sharp }})$, there exists
one solution $x ^{\sharp} ( \cdot ) \in  {\cal S}_{f} (x)$ which remains
viable in $K$ as long as possible:
\begin{displaymath}
\tau _{K}^{f ^{\sharp }} (x) \; = \;  \tau _{K}  (x ^{\sharp} ( \cdot ))
\end{displaymath}
\end{Proposition} 

{\bf Proof} --- \hspace{ 2 mm}
 Since the function  $ x ( \cdot ) \mapsto \omega _{C} (x ( \cdot ))-
\tau _{K} (x ( \cdot ))$ is lower semicontinuous on $ {\cal C} (0,
\infty ,{\bf  R}^{n})$ supplied with the compact convergence by
Theorem~\ref{Picard02},  we deduce first from Theorem~\ref{03A44} 
that the infimum is reached by a solution $x_{ (K,C)} ( \cdot ) \in 
{\cal S}_{f} (x)$ because the set $ {\cal S}_{f} (x)$ is compact and
second, that this function $\gamma _{ (K,C)}^{f ^{\flat }}$ is lower
semicontinuous, by checking that the subsets $ \{x \in K \; | \;
\gamma _{ (K,C)}^{f ^{\flat }} (x) \leq T\}$ are closed. Indeed, let
us consider a sequence of elements $x_{n}$ of such a subset converging
to $x$. There exist solutions $x_{n} ( \cdot ) \in  {\cal S}_{f} (x)$
such that 
\begin{displaymath}
\gamma _{ (K,C)} (x_{n} ( \cdot )) \; \leq  \; \gamma _{ (K,C)}^{f
^{\flat }} (x_{n}) + \frac{1}{n} \; \leq  \; T+ \frac{1}{n}
\end{displaymath}
        On the other hand, since $x_{n}$ belongs to the compact ball $B (x,
1)$, Theorem~\ref{03A44} implies that a subsequence (again denoted by)
$x_{n} ( \cdot )$ converges to some solution $x ( \cdot ) \in  {\cal S}_{f}
(x)$ uniformly on compact intervals. Since the functional $\gamma _{
(K,C)}$ is lower semicontinuous, we infer that
\begin{displaymath}
\gamma _{ (K,C)} ^{f ^{\flat }} (x) \; \leq  \;  \gamma _{ (K,C)}(x (
\cdot )) \; \leq  \; \liminf_{n \rightarrow + \infty } \gamma _{
(K,C)} (x_{n} ( \cdot )) \; \leq  \; T
\end{displaymath}
        
        In particular, taking $C:=K$, we observe that $ \gamma _{ (K,K)}^{f
^{\flat }} =- \tau _{K}^{f ^{\sharp }}$, and thus, we deduce the upper
semicontinuity of the exit function.
The same  proof shows that the hitting function $ \omega _{C}^{f
^{\flat }}$ is lower semicontinuous.
$\; \; \Box$

Viability kernels and capture  basins can be characterized in terms of exit
and  hitting functionals:

\begin{Theorem} \label{characviabkern}
If $f:{\bf   R}^{n} \mapsto {\bf   R}^{n}$ is continuous with linear growth
and $C \subset K
\subset {\bf  R}^{n}$ are closed subsets, then
\begin{displaymath} \left\{ \begin{array}{ll}
  \mbox{\rm Capt}_{f}(K,T) & = \; \left\{ x \in {\bf  R}^{n} \; | \; \omega
_{K}^{f
^{\flat }}(x) \; \leq \;   T\right\} \; \; \& \; \;\mbox{\rm Capt}_{f}(K,T)
\; = \; \mbox{\rm Dom}(\omega _{K}^{f ^{\flat }})\\
 & \\
   \mbox{\rm Capt}_{f}^{K}(C,T)  &  = \; \left\{ x \in {\bf  R}^{n} \; | \;
\gamma _{
(K,C)}^{f ^{\flat }}(x) \; \leq \;   0\right\} \\
 & \\
  \mbox{\rm Viab}_{f}(K,T)  &  :=  \; \left\{ x \in {\bf  R}^{n}\; | \;\tau
_{K} ^{\sharp
}(x) \; \geq  \; T\right\} \\
\end{array} \right. \end{displaymath}
In particular, the $T$-viability kernels $ \mbox{\rm Viab}_{f}(K,T)$ of a
closed  subset  $K  \subset {\bf   R}^{n}$,   the $T$-capture basins of $K$
under $f$ and
the viable-capture basin $\mbox{\rm Capt}_{f}^{K}(C,T) $ are closed.
\end{Theorem} 

{\bf Proof} --- \hspace{ 2 mm}
The subset of initial states $x \in {\bf  R}^{n} $ such that $K$ is reached
before $T$
by  a solution $x(\cdot )$ to  the differential equation $x'=f (x)$
starting at $x$  is obviously contained in the subset  $\left\{  x \in {\bf
R}^{n} \; |
\; \omega_{K} ^{f ^{\flat }}(x) \; \leq \;   T\right\}$.

 Conversely, consider an element $x $ satisfying $ \omega_{K} ^{f ^{\flat
}}(x) \; \leq \;   T$. Hence the solution $x ^{\flat}( \cdot ) \in  {\cal
S}_{f} (x)$ such that $ \omega _{K} (x ^{\flat } ) = \omega _{K}^{f ^{\flat
}} (x) \leq T$ belongs to the $T$-capture basin. 

Now, to say that a solution $x ( \cdot ) \in  {\cal S}_{f} (x)$ is
viable in $K$ until it reaches the target $C$ means that $ \omega 
_{C} (x ( \cdot )) \leq  \tau _{K} (x ( \cdot ))$. Therefore,  $x$
belongs to the viable-capture basin $ \mbox{\rm Capt}_{f}^{K} (C)$ if
and only if $\gamma _{ (K,C)}^{f ^{\flat }} (x) \leq 0$.

The proof of the characterization of the $T$-viability kernel
$\mbox{\rm Viab}_{f}(K,T)$ as upper sections of the exit time function
$\tau  _{K} ^{\sharp }$ is analogous.

The topological properties then follow from the semicontinuity
properties of the above functions stated in
Proposition~\ref{Picard03se}.
$\; \; \Box$

\section{Characterization of the Viability Kernel}
We deduce at once the following consequence:
\begin{Theorem} \label{viablbascharthm}
Let $f:{\bf   R}^{n} \mapsto {\bf  R}^{n}$ be a  continuous map with linear
growth and $K \subset
{\bf   R}^{n}$ be a closed subset. Then the viability kernel is the largest
closed
subset $D \subset K$ viable under $f$, or, equivalently, the largest closed
subset of $K$ satisfying 
\begin{equation} \left\{ \begin{array}{ll}
i) &  D \; \subset  \; K\\
ii) &  \forall \; x \in D, \;  \; f (x) \; \in  \; T_{D} (x)
\end{array} \right. \end{equation}
or, equivalently, in terms of normal cones, the largest closed subset of
$K$ satisfying 
\begin{equation} \left\{ \begin{array}{ll}
i) &  D \; \subset  \; K\\
ii) &  \forall \; x \in D, \;  \forall \; p \in N_{D} (x)\;  \langle p,f
(x) \rangle \; \leq  \; 0
\end{array} \right. \end{equation}
Furthermore, the viability kernel satisfies the following properties
$$ \forall \; x \in \mbox{\rm Viab}_{f} (K) \backslash  \partial K, \;  \;
f (x) \; \in  \; T_{\mbox{\rm Viab}_{f} (K)} (x) \cap - T_{\mbox{\rm
Viab}_{f} (K)} (x)
$$
or, equivalently, in terms of normal cones,
\begin{displaymath}
 \forall \; x \in \mbox{\rm Viab}_{f} (K) \backslash  \partial K, \;
\forall \; p \in N_{\mbox{\rm Viab}_{f} (K)} (x), \;  \;  \langle p,f (x)
\rangle \; = \; 0
\end{displaymath}
\end{Theorem} 
{\bf Proof} --- \hspace{ 2 mm}
The first property follows from the Nagumo Theorem characterizing viable
subsets in terms of tangential and/or normal conditions. The second
property translates the fact that $\mbox{\rm Viab}_{f} (K) \backslash 
\partial K$  is locally backward invariant, and thus, locally backward
viable.
$\; \; \Box$ 

\begin{Proposition} \label{absfinitime}
Let $f:{\bf   R}^{n} \mapsto {\bf  R}^{n}$ be a  continuous map with linear
growth and $K \subset
{\bf   R}^{n}$  be a closed  subset.  If $M\subset  {\bf  R}^{n} \backslash
\mbox{\rm Viab}_{f}(K)$
is  compact,  then, for every $x \in M$ and every  solution $x(\cdot ) \in
{\cal S}_{f}(x)$, there exists $t \in  [0, \sup_{x \in M} \tau _{K}^{f
^{\sharp }} (x)]$ such that $x(t) \notin K$. 
\end{Proposition}  

{\bf Proof} --- \hspace{ 2 mm}
Indeed,    $M$   being  compact   and   the   exit  function   being  upper
semicontinuous, then
$\sup_{x \in M} \tau _{K}^{f ^{\sharp }} (x)$ is finite because, for each
$x \in M$, $ \tau _{K}^{f ^{\sharp }} (x)$ is finite. $\; \; \Box$ 
        
        \mbox{}
        
        In particular,
        \begin{Proposition}
         Let us  assume  that  $K$  is a compact and  that  $f:{\bf  R}^{n}
\mapsto {\bf  R}^{n}$ is
continuous with linear growth. 
Then  either the viability kernel of $K$ is not empty or  $K$ is a 
repeller, \index{repeller} and in this case, 
$\overline{T} := \sup_{x \in K} \tau _{K}^{f ^{\sharp }} (x)$ is finite and
satisfies
\begin{displaymath}
         \mbox{\rm Viab}_{f}(K,\overline{T}) \; \ne  \; \emptyset  \; \; \&
\; \; \forall \; T>\overline{T}, \;  \;  \mbox{\rm Viab}_{f}(K,T) \; = \;
\emptyset  
        \end{displaymath}
        \end{Proposition} 
        {\bf Proof} --- \hspace{ 2 mm}
        When $K$ is a repeller, the exit function is finite. Being compact,
$\overline{T}:=\sup_{x \in K} \tau _{K}^{f ^{\sharp }} (x)$ is thus finite
and achieves its maximum at  some $  \bar{x}$. By Theorem~\ref{Picard03se},
there exists a solution $
\bar{x} ( \cdot ) \in  {\cal S}_{f} ( \bar{x})$ such that  $ \tau _{K} (
\bar{x} ( \cdot )) = \tau _{K}^{f ^{\sharp }} ( \bar{x}) =\overline{T}$.
        $\; \; \Box$ 
        
        \mbox{}
        
        In other words, when $K$ is a compact repeller, there exists a
smallest nonempty $T$-viability kernel of $K$, the ``viability core'', so
to speak, because it is the subset of initial states from which one
solution which enjoys the longest ``life expectation'' $\overline{T} $ in
$K$.  The viability kernel, when it is nonempty, is the viability core with
infinite life expectation.

\section{Characterization of Viable-Capture Basins}

Let $C \subset K$ be  a closed subset of a closed subset $K \subset
{\bf  R}^{n}$.

\begin{Theorem} \label{viablcaptbascharthm}
Let us assume that $f$ is continuous with linear growth and that $K$
is a closed repeller under $f$. Then  the viable-capture basin $
\mbox{\rm Capt}_{f}^{K} (C)$ is the largest  closed subset $D$
satisfying 
\begin{displaymath} \left\{ \begin{array}{ll}
i) &  C \; \subset  \; D \; \subset  \; K\\
ii) &  D \backslash C \; \mbox{\rm is locally viable under $f$}
\end{array} \right. \end{displaymath} 
or, equivalently, is the largest  closed subset $D$ satisfying 
\begin{equation} \left\{ \begin{array}{ll} \label{wonderfulthmeq27}
i) &  C \; \subset  \; D \; \subset  \; K\\
ii) &  \forall \; x\in D \backslash C, \; \; f(x) \in  T_{D} (x) 
\end{array} \right. \end{equation}
or, equivalently, in terms of normal cones,  is the largest closed
subset $D$ satisfying
\begin{equation} \left\{ \begin{array}{ll} \label{wonderfulthmeq28}
i) &  C \; \subset  \; D \; \subset  \; K\\
ii) &  \forall \; x \in D \backslash C, \;  \forall \; p \in N_{D}
(x), \;  \;   \langle  p,f (x)  \rangle \; \leq \; 0
\end{array} \right. \end{equation}
\end{Theorem} 

{\bf Proof of Theorem~\ref{viablcaptbascharthm} } --- \hspace{ 2 mm}
     Assume that a closed subset $D$ such that $C \subset D \subset K$ is a
repeller under $f$ such that $D \backslash C$ is locally viable under $f$
and let us check that it is contained in $ \mbox{\rm Capt}_{f}^{K} (C)$.
Since  $ C \subset \mbox{\rm Capt}_{f}^{K} (C)$, let $x$ belong to $D
\backslash C$ and show that it belongs to  $ \mbox{\rm Capt}_{f}^{K} (C)$.
Since $K$ is a repeller, all solutions starting from $x$ leave $D
\backslash C$ in finite time.
      At least one of them,  the solution $x ^{\sharp }  ( \cdot )\in 
{\cal S}_{f} (x)$ which maximizes $ \tau _{D} (x ( \cdot ))$:
\begin{displaymath}
\tau _{D}^{f ^{\sharp }}(x) \; := \; \sup_{ x(\cdot ) \in {\cal
S}_{f}(x)}\tau _{D}^{f ^{\sharp }}(x(\cdot )) \; = \; \tau _{D}^{f ^{\sharp
}}(x ^{\sharp } (x))
\end{displaymath}
leaves $D \backslash C$ through $C$.
This solution exists by Proposition~\ref{Picard03se} since $D$ is closed
and $f$ is continuous with linear growth. Then we claim that $ x ^{\sharp
}:=x ^{\sharp } (\tau _{D}^{f ^{\sharp }}(x))$ belongs to $C$. If not,  $D
\backslash C$ being locally viable, one could associate with $ x ^{\sharp
}\in D \backslash C$ a solution $y ( \cdot ) \in  {\cal S}_{f} (x ^{\sharp
} )$ and $T>0$ such that $y ( \tau ) \in D \backslash C$ for all $ \tau \in
[0,T]$. Concatenating  this solution to $x ^{\sharp } ( \cdot )$, we obtain
a  solution viable in $D$ on an interval $[0, \tau _{D}^{f ^{\sharp
}}(x)+T]$, which contradicts the definition of $x ^{\sharp } ( \cdot )$.
Furthermore, $x ^{\sharp } ( \cdot )$ is viable in $K$ since $D \subset K$.
This implies that $D \subset \mbox{\rm Capt}_{f}^{K} (C)$.

The viable-capture basin being a closed subset such that $\mbox{\rm
Capt}_{f}^{K} (C) \backslash C$ is locally viable, we conclude  that it is
the largest closed subset $D$ of $K$ containing $C$ such that $D \backslash
C$ is locally viable under $f$.

Since $f$ is continuous and since $D \backslash C$ is locally compact,
the Viability Theorem~\ref{03A131} states that $D \backslash C$ is
locally viable if and only if (\ref{wonderfulthmeq27})ii) or
(\ref{wonderfulthmeq28})ii) holds true.
$\; \; \Box$

\section{Characterization of Capture Basins}

\begin{Theorem} \label{wonderfulthmthmbis}
Let us assume that the closed subset $K \subset {\bf  R}^{n}$ is a repeller
under
$f$ and backward invariant, that $C \subset K$ is closed  and that $f$
is continuous with linear growth.

Then the capture basin $ \mbox{\rm Capt}_{f} (C)$ 
is  the {\em unique\/} closed subset $D $ which satisfies
\begin{displaymath} \left\{ \begin{array}{ll}
i) & C \; \subset  \; D \; \subset  \; K\\
ii) & D\backslash C \;\; \mbox{\rm is locally viable under $f$}\\
iii) & D  \;\; \mbox{\rm is backward invariant under $f$}\\
\end{array} \right. \end{displaymath}
If we assume furthermore that $f$ is Lipschitz, it is 
the unique closed subset satisfying the ``tangential conditions''
\begin{displaymath} \left\{ \begin{array}{ll}
i) & C \; \subset  \; D \; \subset  \; K\\
ii) & \forall \; x \in D \backslash C, \; \;f(x) \; \in  \;T_{D} (x)
\\
iii) &  \forall \; x \in D, \;  \; -f(x) \; \in  \; T_{D} (x)\\
\end{array} \right. \end{displaymath}
or, equivalently, by duality, the ``normal conditions''
\begin{displaymath} \left\{ \begin{array}{ll}
i) & C \; \subset  \; D \; \subset  \; K\\
ii) & \forall \; x \in D \backslash C, \;  \forall \; p \in N_{D} (x),
\;  \;  \langle p,f (x) \rangle  \; = \; 0 \\
iii) &  \forall \; x \in D, \; \forall \; p \in N_{D} (x), \;  \; 
\langle p,f (x) \rangle  \; \geq  \; 0\\
\end{array} \right. \end{displaymath}
\end{Theorem} 
{\bf Proof} --- \hspace{ 2 mm}
Since $C \subset K$ and $K$ is backward invariant, we already observed
that 
\begin{displaymath}
\mbox{\rm Capt}_{f}^{K}(C) \; = \;  \mbox{\rm Capt}_{f}(C) \; \subset  \; K
\end{displaymath}
        
        Since $K$ is a repeller, the capture basin $\mbox{\rm
Capt}_{f}(C)$ is closed: For that purpose, let $x_{n} \in \mbox{\rm
Capt} _{f}(C)$ converge to $x$ and infer that $x$ belongs to the capture
basin of $C$. There exist $t_{n}$ and $c_{n}:= \vartheta _{f}
(t_{n},x_{n})$ which belongs to $C$. Since  $K$ is a repeller, the exit
function $ \tau _{K}^{f ^{\sharp }}$ is finite, and, being upper
semicontinuous thanks to Lemma~\ref{Picard02}, $T:= \sup_{x_{n}\in B
(x,1)}\tau _{K}^{f ^{\sharp }} (x_{n}) <+ \infty $. Since $C \subset K$, we
infer that $t_{n} \leq T <+ \infty $. Hence we can extract a subsequence 
(again denoted by) $t_{n}$ converging to some $t \in [0, T]$, so that, the
reachable map being continuous,  $c_{n}$ converges to $c= \vartheta _{f}
(t,x)$ which belongs to $C$. Hence $x $ belongs to $ \mbox{\rm Capt}
_{f}(C)$, which is then closed.

Therefore, the capture basin $\mbox{\rm Capt}_{f}(C)$ is a closed subset,
backward invariant and locally viable. 

If $D$ is closed and backward invariant, we infer that $\mbox{\rm
Capt}_{f}(C) \subset D$. By Theorem~\ref{viablcaptbascharthm}, if $D
\backslash C$ is locally viable, we know that 
\begin{displaymath}
\mbox{\rm Capt}_{f}(C) \; = \;  \mbox{\rm Capt}_{f}^{K}(C) \; \subset  \; D
\end{displaymath}
Hence,  the capture basin $\mbox{\rm Capt}_{f}(C)$ is the unique closed
subset, backward invariant and locally viable between $C$ and $K$. 
$\; \; \Box$

\section{Stability Properties}

        Let us consider now a sequence of closed subsets $K_{n}$ viable
under a  map $f$. 
        {\em Is the  upper limit   of these closed subsets still viable
under $f$?\/} The answer is positive.
        
        \begin{Theorem} \index{stability of viability domain} 
\label{01A452}
         Let  us  assume  that  $f:{\bf   R}^{n}  \mapsto  {\bf  R}^{n}$ is
continuous with linear
growth.  Then the  upper limit of a sequence of subsets viable under  $f$
is still viable under $f$ and the  upper limit of the viability kernels of
$K_{n}$ is contained in the viability kernel of the upper limit:
        \begin{equation} \label{liminfviabsup}
\mbox{\rm Limsup}_{n \rightarrow + \infty } \mbox{\rm Viab}_{f}(C_{n})
 \; \subset  \;  
\mbox{\rm Viab} _{f}(\mbox{\rm Limsup}_{n \rightarrow + \infty }C_{n})
\end{equation}
         \end{Theorem}
        
        In particular, the intersection of a decreasing family of closed
viability domains is a closed viability domain.

        {\bf Proof} --- \hspace{ 2 mm}
        We shall prove that the  upper limit $K^{\sharp}$ of a sequence of
subsets $K_{n}$ viable under $f$ is still viable under $f$.

        Let $x$ belong to $K^{\sharp}$. It is the limit of a subsequence
$x_{n'} \in K_{n'}$. Since the subsets $K_{n}$ are viable under $f$, there
exist solutions $y_{n'}(\cdot)$ to differential equation $x' =f(x)$
starting at $x_{n'}$ and viable in $K_{n'}$. Theorem~\ref{03A44}  implies
that a subsequence (again denoted) $y_{n'}(\cdot)$ converges uniformly on
compact intervals to a solution $y(\cdot)$ to  differential equation $x'
=f(x)$ starting at $x$. Since $y_{n'}(t)$ belongs to $K_{n'}$ for all $n'$,
we deduce that $y(t)$ does belong to $K^{\sharp}$ for all $t > 0$.
$\; \; \Box$ 

\mbox{}

\begin{Theorem} \label{convgraivenvthm}
Let us consider a sequence of closed subsets $C_{n}$. 
\begin{enumerate}
\item If the map $f$ is continuous with linear growth and if the subsets
$C_{n}$ are contained in a closed repeller $K$,  then 
\begin{equation} \label{liminfenvsup}
\mbox{\rm Limsup}_{n \rightarrow + \infty } \mbox{\rm Capt}_{f}^{K}(C_{n})
 \; \subset  \;  
\mbox{\rm Capt} _{f}^{K}(\mbox{\rm Limsup}_{n \rightarrow + \infty
}C_{n})\end{equation}

\item If the  map $f$ is furthermore monotone and $K$ is backward
invariant, then 
\begin{equation} \label{liminfenvinv}
 \mbox{\rm Capt}_{f}(\mbox{\rm Liminf}_{n \rightarrow + \infty }C_{n}) \;
\subset  \; \mbox{\rm Liminf}_{n \rightarrow + \infty } \mbox{\rm
Capt}_{f}(C_{n})
\end{equation}
\end{enumerate}
\end{Theorem}

{\bf Proof} --- \hspace{ 2 mm}
 For proving inclusion (\ref{liminfenvsup}), we  consider the limit
$x:=\lim_{n \rightarrow + \infty }x_{n}$ of elements $x_{n} $ of $\mbox{\rm
Capt}_{f}^{K}(C_{n})$. Let us consider solutions $x_{n} ( \cdot )$
satisfying
\begin{displaymath}
 t_{n} \; := \; \omega _{C_{n}} (x_{n} ( \cdot )) \; \leq  \; \tau _{K}
(x_{n} ( \cdot )) \; \leq  \; T \; := \; \sup_{y \in B (x,1) \cap K} \tau
_{K}^{f ^{\sharp }} (y)
\end{displaymath}
where $T$ is finite because $K$ is a repeller.
        
        Therefore, a subsequence (again denoted by) $t_{n}$
converges to some $t \leq T$ and another subsequence (again
denoted by) $x_{n} ( \cdot )$ converges uniformly on compact
intervals to some solution $x ( \cdot )$ starting from $x$. Since
$x_{n} (t_{n})$ belongs to $C_{n}$ and converges to $x (t)$, we
infer that $x (t)$ belongs to  the upper limit $C ^{\sharp }$ of
the $C_{n}$.  Hence 
\begin{displaymath}
  \omega_{C ^{\sharp }} (x) \; \leq  \; \lim_{n \rightarrow + \infty }t_{n}
\; \leq  \; \limsup_{n \rightarrow + \infty } \tau _{K}
(x_{n} ( \cdot )) \; \leq  \; \tau _{K}
(x ( \cdot ))
\end{displaymath}
and thus, $x$ belongs to  the viable-capture basin of
$C ^{\sharp }$.

We now prove (\ref{liminfenvinv}).
Let $C ^{\flat }$ denote the lower limit of the subsets $C_{n}$ and let us
consider $x $ an element of $\mbox{\rm Capt}_{f}(C ^{\flat })$, a solution
$x ( \cdot )$ starting from $x$ and reaching $C ^{\flat }$ at time $T$ at
$c:=x (T)$. Hence $y (t):=x (T-t)$ is a solution to the backward
differential equation $y'=-f (y)$ starting at $c$ and satisfying $y (T)=x$.
Since $c= \lim_{n \rightarrow + \infty }c_{n}$ where $c_{n} \in C_{n}$, 
Theorem~\ref{05A13m} states that the solutions $y_{n} ( \cdot )$ to the
differential equation $y'=-f (y)$ starting from $c_{n}$ satisfy
\begin{displaymath}
 \|y (t)-y_{n} (t)\| \; \leq  \; e^{ - \mu t} \|c-c_{n}\| 
\end{displaymath}
Then $x_{n} := y_{n} (T) \in \mbox{\rm Capt}_{f}(C_{n})$ converges to $x=y
(T) \in \mbox{\rm Capt}_{f}(C ^{\flat })$.  $\; \; \Box$

\chapter{Epiderivatives and Subdifferentials}

 For reasons motivated by optimization theory, Lyapunov
stability, control theory, Hamilton-Jacobi equations  and
mathematical morphology,  the order relation on $ {\bf R}$ is
involved. This leads us to associate with an  extended functions
$ {\bf u}: X \mapsto {\bf R} \cup \{+ \infty \}$ its epigraph
instead of its graph. It actually happens that the properties of
the extended functions $ {\bf u}: X \mapsto {\bf R} \cup \{+
\infty \}$ are actually properties of their {\em epigraphs}. This
``epigraphical point of view''  is the key to ``variational
analysis'' and to our treatment of Hamilton-Jacobi inequalities.
        
        In particular, one can define the following concepts:
\begin{enumerate}
\item The epigraph of the {\em lower epilimit} of a sequence of extended
functions $ {\bf u}_{n}:X \mapsto {\bf R} \cup \{+\infty \}$ is the upper
limit of the epigraphs of the $f_{n}$,
\item {\em The contingent epiderivative} $D _{\uparrow } {\bf u} (x)$ at
$x$  is the lower epigraphical limit of the difference quotients $ \nabla
_{h} {\bf u}(x)$, so that the epigraph of the contingent epiderivative is
the contingent cone to the epigraph of $ {\bf u}$.
\end{enumerate}
        
        By duality, the normal cones to the epigraph of an extended
function yields the concept of {\bf subdifferential}, which is in
particular used in the concepts of viscosity solutions to Hamilton-Jacobi
equations.
        
\section{Extended Functions and their Epigraphs}

        A function ${\bf v}: X \mapsto {\bf R} \cup \{+ \infty \}$ is
called an {\em extended  (real-valued) function\/}\index{extended 
(real-valued) function}. Its 
         {\em domain\/}\index{domain (of an extended function)}  is the set
of points at which ${\bf v}$ is finite:
        \begin{displaymath}
        \mbox{\rm Dom}({\bf v})  \;  :=  \;  \{x \in X \;\; | \;\; {\bf
v}(x) < +\infty \} 
        \end{displaymath}
        A function is said  to be \index{proper} {\em nontrivial\/}
\index{nontrivial (function)} if its domain is
        not empty. Any  function ${\bf v}$ defined on a subset $K \subset
X$ can be regarded as the extended function ${\bf v}_{K}$ equal to ${\bf
v}$ on $K$ and to $+ \infty $ outside of $K$, whose domain is $K$.
        
        Since the order relation on the real numbers is involved in the
definition of the Lyapunov property (as well as in minimization problems
and other dynamical inequalities), we no longer characterize a real-valued
function by its graph, but rather by its {\em epigraph\/} index{epigraph} 
        \begin{displaymath}
        {\cal E}p({\bf v}) \; := \; \{(x,\lambda ) \in X \times {\bf R} \;
| \; {\bf v}(x) \leq \lambda \}
        \end{displaymath}
        The {\em hypograph\/}\index{hypograph} of  a function
${\bf v}: X \mapsto {\bf R} \cup \{-\infty \}$ is defined in a
symmetric way by
        \begin{displaymath}
        {\cal H}yp({\bf v}) \; := \; \{(x,\lambda ) \in X \times {\bf R} \;
| \; {\bf v}(x) \geq \lambda \}  \; = \; - {\cal E}p(-{\bf v})
        \end{displaymath}
        
        {\em The  graph of a real-valued (finite) function is then the
intersection of its epigraph and its hypograph}.
        
        We also remark that some properties of a function are actually
properties of their epigraphs. For instance, {\em an extended function
${\bf v}$ is convex (resp. positively homogeneous) if and only if its
epigraph is convex (resp. a cone).\/} The epigraph of ${\bf v}$ is closed
if and only if ${\bf v}$ is lower semicontinuous:
        \begin{displaymath}
        \forall \; x \in X, \;  \; {\bf v}(x) \;   = \; \liminf_{y
\rightarrow x}{\bf v}(y) 
        \end{displaymath}
        
         We recall the convention $\inf (\emptyset) := +\infty  $.
        
        \begin{Lemma}
        Consider a function ${\bf v} :X \mapsto {\bf R} \cup \{\pm \infty
\} $. Its epigraph is closed if and only if
        \begin{displaymath}
        \forall \; x  \; \in \;  X, \;  \; {\bf v} (x)  \; = \; 
\liminf_{x' \rightarrow x} {\bf v} (x') 
        \end{displaymath}
         Assume that the epigraph of ${\bf v} $ is a closed cone. Then the
following conditions are equivalent:
        \begin{displaymath} \left\{ \begin{array}{ll}
        i) & \forall \; x  \; \in \;  X,  \;  \; {\bf v} (x) \;  > \;
-\infty  \\
        ii) & {\bf v} (0)  \; = \;  0 \\
        iii) & (0,-1) \; \notin  \; {\cal E}p({\bf v} )
        \end{array} \right. \end{displaymath}
        \end{Lemma}
        
        {\bf Proof} --- \hspace{ 2 mm}
        Assume that the epigraph of ${\bf v} $ is closed and pick $x \in
X$. There exists a sequence of elements $x_{n}$ converging to $x$ such that
$$\lim_{n \rightarrow \infty } {\bf v} (x_{n}) \;  =  \;  \liminf_{x'
\rightarrow x} {\bf v} (x')$$
        Hence, for any $\lambda > \liminf_{x' \rightarrow x} {\bf v} (x')$,
there exist $N$ such that, for all $n \geq N$, ${\bf v} (x_{n}) \leq
\lambda $, i.e., such that $(x_{n},\lambda ) \in {\cal E}p({\bf v} )$. By
taking the limit, we infer that ${\bf v} (x) \leq \lambda $, and thus, that
${\bf v} (x) \leq \liminf_{x' \rightarrow x} {\bf v} (x')$. The converse
statement is obvious.
        
        Suppose next that the epigraph of ${\bf v} $ is a cone. Then it 
contains $(0,0)$ and ${\bf v}(0) \leq 0$. The statements $ii)$ and $iii)$
are clearly equivalent.
        
        If $i)$ holds true and ${\bf v}(0)<0$, then 
        $$(0,-1) \;  = \;  \frac{1}{-{\bf v}(0)}(0,{\bf v}(0)) $$
        belongs to the epigraph of ${\bf v}$, as well as all $(0,-\lambda
)$, and (by letting $\lambda \rightarrow +\infty $) we deduce that ${\bf
v}(0)=-\infty $, so that $i)$ implies $ii)$.
        
        To end the proof, assume that ${\bf v} (0) =0$ and that for some
$x$, ${\bf v} (x) =-\infty $. Then, for any $\varepsilon >0$, the pair $(x,
-1/\varepsilon )$ belongs to the epigraph of ${\bf v} $, as well as the
pairs $(\varepsilon x,-1)$. By letting $\varepsilon $ converge to $0$, we
infer that $(0,-1)$ belongs also to the epigraph, since it is closed. Hence
${\bf v}(0)<0$, a contradiction.
        $\; \; \Box$
        
        \mbox{}
        
        {\em Indicators $\psi _{K}$ of subsets $K$\/}\index{indicators }
are cost functions defined by
        \begin{displaymath}
        \psi _{K}(x) := 0 \;\;\mbox{\rm if}\;\; x \in K \;\;\mbox{\rm
and}\;\; + \infty \;\;\mbox{\rm if not}
        \end{displaymath}
        which characterize subsets (as {\em characteristic functions\/} do
for other purposes) provide important examples of extended functions.
        It can be regarded as a {\em membership cost}\footnote{Functions
${\bf v}:X \mapsto [0, + \infty ]$ can be regarded as some kind of {\em
fuzzy sets}, called {\em toll sets}. \index{toll sets}} to $K$: it costs
nothing to belong to $K$, and $+ \infty $ to step outside of $K$.
\index{membership cost}
        
        Since 
        \begin{displaymath}
         {\cal E}p ( \psi K) \; = \; K \times {\bf  R}_{+} 
        \end{displaymath}
        we deduce that the indicator $\psi _{K}$ is  lower semicontinuous
if and only if $K$ is closed and that $\psi _{K}$ is convex if and only if
$K$ is convex. One can regard the sum ${\bf v} + \psi _{K}$ as the
restriction of ${\bf v}$ to $K$.
        
         We recall the convention $\inf (\emptyset) := +\infty  $.

 \section{Epilimits}
        \begin{Definition}
        The epigraph of the {\em lower epilimit} $\mbox{\rm
lim}_{\uparrow}^{\sharp}\mbox{}_{n \rightarrow \infty }{\bf u}_{n}$
        of a sequence of extended functions $ {\bf u}_{n} :X \mapsto {\bf
R} \cup \{+\infty \}$ is  the upper limit of the epigraphs:
        \begin{displaymath}
        {\cal E}p(\mbox{\rm lim}_{\uparrow}^{\sharp}\mbox{}_{n \rightarrow
\infty }{\bf u}_{n})   \;  :=  \;  \mbox{\rm Limsup}_{n \rightarrow \infty
} {\cal E}p({\bf u}_{n})  
        \end{displaymath}
        The function $\mbox{\rm lim}_{\uparrow}^{\flat}\mbox{}_{n
\rightarrow \infty }{\bf u}_{n}$ whose epigraph is the  lower limit of the
epigraphs of the functions $ {\bf u}_{n}$
        \begin{displaymath} 
        {\cal E}p( \mbox{\rm lim}_{\uparrow}^{\flat}\mbox{}_{n \rightarrow
\infty }{\bf u}_{n})   \;  :=  \;  \mbox{\rm Liminf}_{n \rightarrow \infty
} {\cal E}p({\bf u}_{n})  
        \end{displaymath}
        is the {\em  upper epilimit}\index{upper epilimit} of the functions
$ {\bf u}_{n}$
         \end{Definition}

        One can check that
        $$\mbox{\rm lim}_{\uparrow}^{\sharp}\mbox{}_{n \rightarrow \infty
}{\bf u}_{n}(x_{0}) 
         \;  =  \; \liminf_{n \rightarrow \infty, x \rightarrow x_{0}}{\bf
u}_{n}(x) $$

        \section{Contingent Epiderivatives}
        
        When $ {\bf u}$ is an extended function, we associate with it its
epigraph  and the contingent cones to this epigraph. This leads to the
concept of epiderivatives of extended functions.
        
        \begin{Definition}
        
        Let  $ {\bf u}: X \mapsto {\bf R} \cup \{\pm \infty \}$ be a
nontrivial extended  function and $x$ belong to its domain.
        
        We associate with it the {\em differential quotients}
        \begin{displaymath}
        u \; \leadsto \; \nabla _{h}{\bf u}(x)(u) \; := \; \frac{{\bf
u}(x+hu)-{\bf u}(x)}{h}
        \end{displaymath}
        
         The contingent epiderivative $D _{\uparrow }{\bf u}(x)$ of $ {\bf
u}$ at $ x \in  \mbox{\rm Dom}( {\bf u})$ is the lower epilimit of its
differential quotients:
        \begin{displaymath}
         D _{\uparrow }{\bf u}(x) \; = \; \mbox{\rm
lim}_{\uparrow}^{\sharp}\mbox{}_{h \rightarrow 0+}\nabla _{h}{\bf u}(x)
        \end{displaymath}
         We shall say that the function $ {\bf u}$ is {\em contingently
epidifferentiable\/}\index{contingently epidifferentiable}  at $x$ if for
any $u \in X$, $D_{\uparrow }{\bf u}(x)(u) >-\infty $ (or, equivalently, if
$D_{\uparrow }{\bf u}(x)(0)=0$).
         \end{Definition}

        \begin{Proposition} \label{cgepider}
        Let  $ {\bf u}: X \mapsto {\bf R} \cup \{\pm \infty \}$ be a
nontrivial extended  function and $x$ belong to its domain. Then the
contingent epiderivative $D _{\uparrow }{\bf u}(x)$ satisfies
        \begin{displaymath}
        \forall \; u \in X,  \;  \; D_{\uparrow }{\bf u}(x)(u) \; = \; 
\liminf_{h \rightarrow 0+, u'\rightarrow u}\frac{{\bf u}(x+hu')-{\bf
u}(x)}{h}
        \end{displaymath}
        and the epigraph of  the contingent epiderivative $D _{\uparrow
}{\bf u}( \cdot )$ is equal to the contingent cone to the epigraph of $
{\bf u}$ at $(x,{\bf u}(x))$ is \begin{displaymath} 
        {\cal E}p (D_{\uparrow}{\bf u}(x))  \; = \;   T_{{\cal E}p ( {\bf
u})}(x,{\bf u}(x))  
        \end{displaymath}
        \end{Proposition} 
        
        {\bf Proof} --- \hspace{ 2 mm}
        The first statement is obvious. For proving the second one,  we
recall that the contingent cone $$T_{{\cal E}p ( {\bf u})}(x,{\bf u}(x)) \;
= \; \mbox{\rm Limsup}_{h \rightarrow 0+} \frac{{\cal E}p ( {\bf u})
-(x,{\bf u}(x))}{h}$$ is the upper limit of the differential quotients $
\displaystyle{ \frac{{\cal E}p ( {\bf u})-(x,{\bf u}(x))}{h}}$ when $h
\rightarrow 0+$. It is enough to observe that 
        \begin{displaymath}
        {\cal E}p(D _{\uparrow }{\bf u}(x)) :=  T_{{\cal E}p ( {\bf
u})}(x,y)
        \; \; \& \; \; {\cal E}p(\nabla _{h}{\bf u}(x)) =  \frac{{\cal E}p
( {\bf u}) -(x,{\bf u}(x))}{h}
        \end{displaymath}
        to conclude.
        $\; \; \Box$

        Consequently,  {\em the epigraph of the contingent epiderivative at
$x$ is a closed cone. It is then lower semicontinuous and positively
homogeneous whenever $ {\bf u}$ is contingently epidifferentiable at
$x$.\/}
        
        \mbox{}
        
        We observe that the contingent epiderivative of the indicator
function $ \psi _{K}$ at $x \in K$ is the indicator of the contingent cone
to $K$ at $x$:
        \begin{displaymath}
         D _{\uparrow } \psi _{K} (x) \; = \;  \psi _{T_{K} (x)}
        \end{displaymath}
        making precise the intuition stating that the contingent cone
$T_{K} (x)$ plays the role of a ``derivative of a set'', as the limit of
differential quotients $ \displaystyle{ \frac{K-x}{h}}$ of sets.

        The hypoderivatives of an extended function are  defined in a
analogous way: The contingent hypoderivative $D _{\downarrow }{\bf u}(x)$
of $ {\bf u}$ at $ x \in  \mbox{\rm Dom}( {\bf u})$ is the upper hypolimit
of its differential quotients:
        \begin{displaymath}
         D _{ _{\downarrow } }{\bf u}(x) \; = \; \mbox{\rm lim}_{
_{\downarrow } }^{\sharp}\mbox{}_{h \rightarrow 0+}\nabla _{h}{\bf u}(x)
        \end{displaymath}
        We observe that it is equal to
        \begin{displaymath}
        \forall \; u \in X,  \;  \; D_{\downarrow }{\bf u}(x)(u) \; = \; 
\limsup_{h \rightarrow 0+, u'\rightarrow u}\frac{{\bf u}(x+hu')-{\bf
u}(x)}{h}
        \end{displaymath}
        and that {\em the hypograph of the {\em contingent
hypoderivative\/}\index{contingent hypoderivative} $D_{\downarrow }{\bf
u}(x) $ of $ {\bf u}$ at $x$ is the contingent cone to the hypograph of $
{\bf u}$ at $(x,{\bf u}(x))$\/}:
        \begin{displaymath} 
        {\cal E}p (D_{\downarrow}{\bf u}(x))  =  T_{{\cal H}yp( {\bf
u})}(x,{\bf u}(x))  
        \end{displaymath}
        
        \begin{Definition}
        We shall say that $ {\bf u}: X \mapsto W$ is {\em differentiable
from the right\/} \index{differentiable from the right} at $x$ if the
contingent epiderivative and hypoderivative coincide:
        \begin{displaymath}
         \forall \; v \in X, \;  \; D _{\uparrow }{\bf u}(x) (v) \; = \; D
_{\downarrow } {\bf u}(x) (v)
        \end{displaymath}
         \end{Definition}
        
        \begin{Lemma}
        Let $K \subset X$ be a closed subset and $ {\bf u}:X \mapsto  {\bf
R} \cup \{+\infty \}$ be an extended function. We denote by $ {\bf
u}|_{K}:=f + \psi _{K}$ the restriction to $ {\bf u}$ at $K$. Inequality
        \begin{displaymath}
         D _{\uparrow }{\bf u}(x) |_{T_{K}} (x) \; \leq  \; D _{\uparrow
}{\bf u}|_{K} (x)
        \end{displaymath}
        always holds true. It is an equality when $ {\bf u}$ is 
differentiable from the right:
         {\em the contingent derivative of the restriction of $ {\bf u}$ to
$K$ is the restriction of the derivative to the contingent cone}.
        \end{Lemma}
        {\bf Proof} --- \hspace{ 2 mm}
        Indeed, let $x \in K \cap \mbox{\rm Dom}( {\bf u})$. If $ u $
belongs to $T_{K} (x)$, there exist $h_{n} \rightarrow 0+$, $ \varepsilon
_{n} \rightarrow 0+$ and $x_{n}:=x+h_{n}u_{n} \in K$ such that
\begin{displaymath}
         D _{\uparrow }{\bf u}(x) ( u ) \; \leq  \; \liminf_{n \rightarrow
+ \infty } \frac{{\bf u}(x_{n})-{\bf u}(x)}{h_{n}} \; = \;\liminf_{n
\rightarrow + \infty } \frac{{\bf u}|_{K}(x_{n})-{\bf u}|_{K}(x)}{h_{n}}
        \end{displaymath}
        which implies the inequality. If $ {\bf u}$ is  differentiable from
the right, the differential quotient converges to the common value $ D
_{\uparrow }{\bf u}(x) = D _{\uparrow }{\bf u}|_{K} (x)= D _{\downarrow }
{\bf u}|_{K} (x)$.
        $\; \; \Box$ 
        
         \vspace{ 3 mm}
        
        \mbox{}
        
        For locally Lipschitz functions, the contingent epiderivatives are
finite:
        
        \begin{Proposition} \label{04A9627}
        If $ {\bf u} :X \mapsto {\bf R} \cup \{+\infty \}$ is Lipschitz
around $x \in \mbox{\rm Int}( \mbox{\rm Dom}( {\bf u}))$, then the
contingent epiderivative $D _{\uparrow }{\bf u}(x)$ is Lipschitz: there
exists $ \lambda >0$ such that
        \begin{displaymath}
        \forall \; u \in X, \;  \;  D _{\uparrow }{\bf u}(x) (u) \; = \;
\liminf_{h \rightarrow 0+} \frac{{\bf u}(x+hu)-{\bf u}(x)}{h} \; \leq  \;
\lambda \|u\|
        \end{displaymath}
        \end{Proposition} 
        {\bf Proof} --- \hspace{ 2 mm}
        Since $ {\bf u}$ is Lipschitz on some ball $B (x, \eta )$, the
above inequality follows immediately from
        \begin{displaymath}
         \forall \; u \in \eta B, \;  \; \frac{{\bf u}(x+hu)-{\bf u}(x)}{h}
\; \leq  \; \frac{{\bf u}(x+hu')-{\bf u}(x)}{h} + \lambda  (\|u\|+
\|u'-u\|)
        \end{displaymath}
        by taking the liminf when $h \rightarrow 0+$ and $u' \rightarrow
u$. $\; \; \Box$ 
        
        \mbox{}
        
        For convex functions, we obtain:
        \begin{Proposition} \label{04A962}
        When the function $ {\bf u}:X\mapsto {\bf R}\cup \{+\infty \}$ is
convex, the contingent epiderivative is equal to
        \begin{displaymath}
        D_{\uparrow }{\bf u}(x)(u)   \;  =  \;  \liminf_{u' \rightarrow
u}\left( \inf_{h>0}\frac{{\bf u}(x+hu')-{\bf u}(x)}{h}   \right)  
        \end{displaymath}
        \end{Proposition} 
        {\bf Proof} --- \hspace{ 2 mm}
        Indeed, Proposition~\ref{01A955} implies that if $0< h_{1} \leq
h_{2}$,
        \begin{displaymath}
         {\cal E}p(\nabla _{h_{2}}{\bf u}(x)) \; = \;   \frac{{\cal E}p (
{\bf u}) -(x,{\bf u}(x))}{h_{2}} \; \subset  \;  {\cal E}p(\nabla
_{h_{1}}{\bf u}(x))
        \end{displaymath}
        i.e.,
        \begin{displaymath}
        \forall \; u \in X, \;  \;\nabla _{h_{1}}{\bf u}(x)(u) \; \leq  \; 
\nabla _{h_{2}}{\bf u}(x)(u)
        \end{displaymath}
        Therefore, 
        \begin{displaymath}
         \forall \; u\in X, \;  \; D{\bf u}(x) (u) \; := \; \lim_{h
\rightarrow 0+} \frac{{\bf u}(x+hu)-{\bf u}(x)}{h} \; = \;
\inf_{h>0}\frac{{\bf u}(x+hu)-{\bf u}(x)}{h}
        \end{displaymath}
        and this function $D{\bf u}(x)$ is convex with respect to $u$.
Since the epigraph of $D{\bf u}(x)$ is the increasing union of the
epigraphs of the differential quotients $\nabla _{h}{\bf u}(x)$,  we infer
that
        \begin{displaymath}
         D _{\uparrow }{\bf u}(x) (u) \; := \; \liminf_{u' \rightarrow
u}D{\bf u}(x) (u')
        \end{displaymath}
        
        \mbox{}
        
        We recall the following important property of convex functions
defined on finite dimensional vector spaces:
        
        \begin{Theorem}
        An extended  convex function $ {\bf u}$ defined on a finite
dimensional vector-space is locally Lipschitz and subdifferentiable on the
interior of its domain. Therefore, when $x$ belongs to the interior of the
domain of $ {\bf u}$, there exists a constant $ \lambda _{x}$ such that
        \begin{displaymath}
        \forall \; u  \; \in \;  X, \;  \; D_{\uparrow }{\bf u}(x)(u)   \; 
=  \;   \inf_{h>0}\frac{{\bf u}(x+hu)-{\bf u}(x)}{h} \; \leq \; \lambda_{x}
\|u\|
        \end{displaymath}
        \end{Theorem}  
        The second statement follows from Proposition~\ref{04A9627}.
        $\; \; \Box$

        \section{Generalized Gradients}

         \begin{Definition}
Let $ {\bf u}:X \mapsto {\bf R} \cup \{+\infty \}$ be a nontrivial extended
function.        The continuous linear functionals $p \in X^{\star}$
satisfying
\begin{displaymath}
\forall \; v \in X, \;  \; \langle p,v \rangle  \;  \leq  \; D_{\uparrow
}{\bf u}(x)(v)
\end{displaymath}
are called the {\em (regular) subgradients} \index{subgradients} of $ {\bf
u}$ at $x$, which constitute the (possibly empty) closed convex subset
        \begin{displaymath}
\partial _{-}{\bf u}(x)    \; := \;      \{p \in X^{\star}\;|\;\forall \; v
\in X, \;\; \langle p,v \rangle   \;  \;
\leq  \; D_{\uparrow }{\bf u}(x)(v)\}
        \end{displaymath}
        called the {\em (regular) subdifferential\/}\index{
subdifferential}   of $ {\bf u}$ at $x_{0}$. 

        In a symmetric way, the {\em superdifferential\/}\index{
superdifferential}  $\partial
_{+}{\bf u}(x)$ of $ {\bf u}$ at $x$ is defined by 
$$\partial _{+}{\bf u}(x) \;  := \;  - \partial
_{-}(- {\bf u})(x)$$
 \end{Definition}
        
        Naturally, when $ {\bf u}$ is Fr\'{e}chet differentiable at $x$,
then
$$D_{\uparrow }{\bf u}(x)(v) \;  = \;  \langle f'(x),v \rangle   $$ so that
{\em the
subdifferential $\partial _{-}{\bf u}(x)$ is reduced to the
gradient $ {\bf u}'(x)$.}

We observe that
\begin{displaymath}
  \partial _{-} {\bf u}(x) + N_{K} (x) \; \subset  \; \partial ({\bf
u}|_{K}) (x)
\end{displaymath}

        If $ {\bf u}$ is differentiable at a point $x \in K$,
then the {\em subdifferential of the restriction is the sum of the 
gradient and the normal cone}:
        \begin{displaymath}
        \partial_{-}({\bf u}|_{K})(x)   \; =   \; {\bf u}'(x) + N_{K}(x) 
        \end{displaymath}
         
        We also note that the subdifferential of the indicator of a
subset is the normal cone:
        \begin{displaymath}
        \partial_{-}\psi _{K}(x)  \; =  \; N_{K}(x) 
        \end{displaymath}
that
     \begin{displaymath} \left\{ \begin{array}{ll}
     i) & (p,-1) \in  N_{ {\cal E}p ( {\bf u})} (x,{\bf u}(x)) \;\mbox{\rm
if and only if}\;p \in  \partial_{-}{\bf u}(x)\\
     ii) & (p,0) \in  N_{ {\cal E}p ( {\bf u})} (x,{\bf u}(x)) \;\mbox{\rm
if and only if}\;  p  \in \mbox{\rm Dom}(D _{\uparrow }{\bf u}(x))^{-}
     \end{array} \right. \end{displaymath}
     so that
We also deduce that
\begin{displaymath}
    N_{ {\cal E}p ( {\bf u})} (x,{\bf u}(x)) \; = \; 
  \{ \lambda (q,-1)\}_{q \in \partial_{-}{\bf u}(x), \; \lambda >0}
\bigcup_{}^{} \{ (q,0)\}_{q \in \mbox{\rm Dom}(D _{\uparrow }{\bf
u}(x))^{-}}
\end{displaymath}
 The subset $ \mbox{\rm Dom}(D _{\uparrow }{\bf u}(x))^{-} = \{0\}$
whenever the domain of the contingent epiderivative $D _{\uparrow }{\bf
u}(x)$ is dense in $X$. 
This happens when $ {\bf u}$ is locally Lipschitz and when  the dimension
of $X$ is finite:

        \begin{Proposition}
        Let $X$ be a finite dimensional vector space, $ {\bf u}:X \mapsto
{\bf R} \cup \{\pm \infty \}$ be a nontrivial extended function
and $x_{0} \in \mbox{\rm Dom}( {\bf u})$. Then the subdifferential
$\partial _{-}{\bf u}(x)$ is the set  of elements $p \in
X^{\star}$ satisfying 
        \begin{equation} \label{af-064-27} 
        \liminf_{x \rightarrow
x_{0}}\frac{{\bf u}(x)-{\bf u}(x_{0})-  \langle p,x-x_{0} \rangle
}{\|x-x_{0}\|}  \;  \geq \;  0
        \end{equation}
        is the \index{local subdifferential}{\em local subdifferential of $
{\bf u}$ at $x_{0}$\/}. 
        
        In a symmetric way, the {\em superdifferential\/}\index{local
superdifferential}  $\partial
_{+}{\bf u}(x_{0})$ of $ {\bf u}$ at $x_{0}$ is the subset of elements $p
\in X^{\star}$ satisfying 
        \begin{displaymath}
        \limsup_{x \rightarrow
x_{0}}\frac{{\bf u}(x)-{\bf u}(x_{0})-  \langle p,x-x_{0} \rangle
}{\|x-x_{0}\|}   \; \leq  \;  0
        \end{displaymath}
  \end{Proposition}
{\bf Proof} --- \hspace{ 2 mm}
This is an easy consequence of proposition~\ref{mormconechar}.
$\; \; \Box$

The equivalent formulation (\ref{af-064-27}) of the concept of
subdifferential has been introduced  by Crandall \& P.-L. Lions for
defining {\em viscosity solutions\/}\index{viscosity solutions}  to
Hamilton-Jacobi equations.

\mbox{}

\section{Moreau-Rockafellar Subdifferentials}
        
        When $ {\bf u}$ is convex, the generalized gradient coincides with
the subdifferential introduced by Moreau and Rockafellar for convex
functions in the early 60's:
        
        \begin{Definition}
        Consider a nontrivial function $ {\bf u}:X\mapsto {\bf R}\cup
\{+\infty \}$ and $x \in \mbox{\rm Dom}( {\bf u})$. The closed convex
subset $ \partial {\bf u}(x)$ defined by
        \begin{displaymath}
        \partial {\bf u}(x)  \; = \;  \{p  \; \in \;  X^{\star}\;| \;
\;\forall \; y  \; \in \;  X, \; \langle p,y-x \rangle    \;  \leq  \; {\bf
u}(y)-{\bf u}(x)\}
        \end{displaymath}
        (which may be empty) is called the {\em   Moreau-Rockafellar
subdifferential\/} \index{Moreau-Rockafellar subdifferential} of $ {\bf u}$
at $x$. We say that $ {\bf u}$ is {\em   subdifferentiable at $x$\/} if $
\partial {\bf u}(x) \ne \emptyset $.
         \end{Definition}
         
        Proposition~\ref{03A955} implies that in the convex case 
        \begin{Proposition}
        Let $ {\bf u}:X \mapsto {\bf  R}_{+} $ be a nontrivial extended
convex function. Then the subdifferential $ \partial_{-}{\bf u}(x)$
coincides with Moreau-Rockafellar subdifferential $ \partial {\bf u}(x)$.
        
        Furthermore, the graph of the subdifferential map $x \leadsto 
\partial {\bf u}(x)$ is closed.
        \end{Proposition} 
        
        Let us mention the following simple --- but useful --- remark:
        
        \begin{Proposition} \label{lopsidfrematruleprp}
        Assume that $ {\bf u}:={\bf v}+{\bf w}$ is the sum of a
differentiable function ${\bf v}$ and a convex function $ {\bf w}$. If $
\bar{x}$ minimizes $ {\bf u}$, then
        \begin{displaymath}
         -{\bf v}' ( \bar{x}) \; \in  \; \partial {\bf w} ( \bar{x})
        \end{displaymath}
        \end{Proposition} 
        {\bf Proof} --- \hspace{ 2 mm}
        Indeed, for $h>0$ small enough, $\bar{x} +h (y-\bar{x}) =
(1-h)\bar{x} +hy$ so that
        \begin{displaymath}
         0 \; \leq  \; \frac{{\bf u}(\bar{x} +h (y-\bar{x}))- {\bf u}
(\bar{x})}{h} \; \leq  \; \frac{{\bf u}(\bar{x} +h (y-\bar{x}))- {\bf u}
(\bar{x})}{h} + {\bf w} (y)-{\bf w} ( \bar{x})
        \end{displaymath}
        thanks to the convexity of ${\bf w}$. Letting $h$ converge to $0$
yields
        \begin{displaymath}
         0 \; \leq  \;  \langle {\bf v}' (\bar{x}), y-\bar{x} \rangle +
{\bf w} (y) - {\bf w} (\bar{x})
        \end{displaymath}
        so that $-{\bf v}' ( \bar{x}) $ belongs to $  \partial {\bf w} (
\bar{x})$.
        $\; \; \Box$ 
        
\chapter{Some Hamilton-Jacobi Equations}

\vspace{ 22 mm}
{\Huge \bf Introduction}

\vspace{ 13 mm}

Let us introduce 
\begin{enumerate}
\item a differential  equation  $x'=  f(x)$,  where $f:{\bf  R}^{n} \mapsto
{\bf  R}^{n}$ is continuous and has linear growth,
\item a nonnegative continuous ``Lagrangian''
 $$l : (x,p) \in {\bf  R}^{n} \times {\bf  R}^{n}  \mapsto l (x,p) \in {\bf
R}_{+}$$ 
\item an extended  nonnegative lower semicontinuous  function ${\bf u}:{\bf
R}^{n} \mapsto {\bf R}_{+} \cup
\{+\infty \}$
\end{enumerate}
We consider the problem
\begin{displaymath}
{\bf u}^{\top} (x)  \;  := \;  \alpha _{ (f,l)} ^{\top} ({\bf u}) (x) \; :=
\; \inf_{x (
\cdot )  \in {\cal S}_{f} (x)} \left( \sup_{t \geq 0} \left(   e^{a  t}{\bf
u} (x
(t)) + \int_{0}^{t}e^{a   \tau }l ( x( \tau ),x' ( \tau )) d \tau \right)
\right)
\end{displaymath}
and the ``stopping time'' problem
\begin{displaymath}
{\bf u}^{\bot} (x)  \; := \;  \alpha _{ (f,l)} ^{\bot } ({\bf u}) (x) \; :=
\;  \inf_{x
(  \cdot  )  \in  {\cal S}_{f}  (x)}  \left(  \inf_{t \geq 0} \left(   e^{a
t}{\bf u}
(x (t)) + \int_{0}^{t}e^{a   \tau }l ( x( \tau ),x' ( \tau )) d \tau
\right) \right)
\end{displaymath} 
whenever the graph of the function ${\bf u}$ is regarded as an obstacle, as
in
unilateral mechanics.
Taking $l \equiv 0$, we obtain the $a$-Lyapunov function
\begin{displaymath}
   \alpha _{ (f,0)} ^{\top} ({\bf  u})  (x)  \; := \; \inf_{x ( \cdot ) \in
{\cal
S}_{f} (x)}  \left( \sup_{t \geq 0}   e^{a  t}{\bf u} (x (t)) \right)
\end{displaymath}
as an example of the first problem and taking ${\bf u} \equiv 0$, we obtain
the
variational problem
\begin{displaymath}
 \alpha _{ (f,l)} ^{\top} (0) (x) \; := \; \inf_{x ( \cdot ) \in {\cal
S}_{f} (x)} \int_{0}^{ \infty }e^{a   \tau }l ( x( \tau ),x' ( \tau )) d
\tau 
\end{displaymath}

We shall prove that these functions are ``generalized'' solutions to
Hamilton-Jacobi ``differential variational inequalities''
        \begin{displaymath} \left\{ \begin{array}{ll}
        i) &  {\bf u} (x) \; \leq  \; {\bf u} ^{\top } (x)\\
         ii)  &   \displaystyle{\left\langle    \frac{ \partial }{ \partial
x}{\bf u}
^{\top } (x),f (x)\right\rangle +l (x,f (x))  +a{\bf u} ^{\top } (x)\; \leq
\;
0}\\
         iii)  &  \displaystyle{({\bf  u} (x)-{\bf u}  ^{\top } (x)) \left(
\left\langle  
\frac{ \partial }{  \partial x}{\bf u}  ^{\top }  (x),f (x)\right\rangle +l
(x,f
(x))  +a{\bf u} ^{\top } (x) \right) \; = \; 0}
        \end{array} \right. \end{displaymath} 
and
        \begin{displaymath} \left\{ \begin{array}{ll}
        i) &  0 \; \leq  \; {\bf u} ^{\bot }(x) \; \leq  \; {\bf u}  (x)\\
         ii)  &   \displaystyle{\left\langle    \frac{ \partial }{ \partial
x}{\bf u}
^{\bot } (x),f (x)\right\rangle +l (x,f (x))  +a{\bf u} ^{\bot } (x)\; \geq
\;
0}\\
         iii)  &  \displaystyle{({\bf u}  (x)-{\bf u} ^{\bot  } (x)) \left(
\left\langle  
\frac{ \partial  }{ \partial x}{\bf  u} ^{\bot } (x),f  (x)\right\rangle +l
(x,f
(x))  +a{\bf u} ^{\bot } (x) \right) \; = \; 0}
        \end{array} \right. \end{displaymath} 
respectively.

For this type of partial differential equations, the concept of
distributional derivatives happens to be much less adequate than other
concepts of  ``contingent epiderivatives'',  as it was shown  by H\'el\`ene
Frankowska, or  
subdifferentials, as they appear in the concept of ``viscosity
solutions'' introduced by Michael Crandall and Pierre-Louis Lions for a
general class of nonlinear Hamilton-Jacobi equations.

Indeed, for those Hamilton-Jacobi type equations derived from the calculus
of variations, control theory or differential games, one can  derive   from
the properties of the viability kernels and capture basins of auxiliary
systems that the various value functions involved in these problems are
solutions to such Hamilton-Jacobi equations or differential variational
inequalities in an adequate sense. 

We illustrate here this general approach for the two preceding problems only.
Namely, we shall prove that the epigraphs of the two ``value functions'' $
{\bf u}^{\top}$ and $ {\bf u}^{\bot}$ are respectively the viability kernel
and the
capture basin of the epigraph of the function ${\bf u}$  under  the map $g:
{\bf  R}^{n}
\times {\bf R} \mapsto {\bf  R}^{n} \times  {\bf R}$ defined by
\begin{displaymath}
 g (x,y) \; := \;  ( f (x), -ay-l (x,f (x)))
\end{displaymath}
which is a continuous map with linear growth whenever $f$ and the
Lagrangian $l$ are continuous with linear growth.
The basic observation allowing this transfer of properties is summarized in
the statement of the following:
\begin{Lemma} \label{tranprplemm}
If a closed subset $ {\cal V} \subset {\bf  R}^{n} \times {\bf  R}_{+} $ is
locally
viable under $g$,  so is the  epigraph of the associated  extended function
${\bf v}
:{\bf  R}^{n} \mapsto {\bf R}_{+} \cup \{+\infty \}$ defined by
\begin{displaymath}
 {\bf v} (x) \; := \; \inf_{ (x,y) \in  {\cal V}}y
\end{displaymath}
where we set as usual ${\bf v} (x)  :=+  \infty $ whenever the subset $ \{y
\; | \;
(x,y) \in  {\cal V}\}$ is empty. As a consequence, for any $x \in \mbox{\rm
Dom}({\bf v})$,  there exist a solution to the  differential equation $x'=f
(x)$
and $T>0$ satisfying 
\begin{equation} \label{dynineeq}
\forall \;  t \in [0,T], \;  \; e^{a  t} {\bf v} (x (t)) + \int_{0}^{t}e^{a
\tau
}l ( x( \tau ),x' ( \tau )) d \tau  \; \leq  \; {\bf v} (x)
\end{equation}
(where $T=+ \infty $ whenever $ {\cal V} $ is (globally) viable under $g$).
\end{Lemma}
{\bf Proof} --- \hspace{ 2 mm}
If $x$  belongs to the domain of ${\bf v}$, then $ (x,{\bf v} (x))$ belongs
to $  {\cal V}$  so that there exist a solution $   x (  \cdot ) \in  {\cal
S}_{ f} (x)$
and $T>0$ such that
  \begin{displaymath}
\forall \; t \in [0,T], \;  \;\left( x (t) , y (t):=e^{-at}{\bf v} (x)-
\int_{0}^{t}e^{-a(t - \tau }l ( x( \tau ),x' ( \tau )) d \tau \right) \;
\in  \; {\cal V}
\end{displaymath}
i.e.,  if and only if $y (t) \geq {\bf v} (x (t))$, which can be written in
the
form (\ref{dynineeq}).  If $y_{0} >{\bf v} (x)$,  then we  observe that for
all $t
\geq 0$,
  \begin{displaymath}
y_{0} (t) \; := \;  e^{-at}y_{0} -\int_{0}^{t}e^{-a(t - \tau }l ( x( \tau
),x' ( \tau )) d \tau  \; \geq  \; y (t)\; \geq  \; {\bf v} (x (t))
\end{displaymath}
and thus, that $ (x (t), y_{0} (t))$ is a solution to the differential
equation $ (x' (t),y' (t))=g (x (t),y (t))$ starting at $ (x,y_{0})$ and
viable in the epigraph of ${\bf v}$.
$\; \; \Box$ 

\section{Value Function}
  We  associate  with   ${\bf  u}:{\bf   R}^{n}  \mapsto  {\bf  R}_{+} \cup
\{+\infty \}$ the problem
\begin{displaymath}
 {\bf u}^{\top} (x) \; := \; \inf_{x ( \cdot ) \in {\cal S}_{f} (x)} \left(
\sup_{t \geq 0} \left(    e^{a  t}{\bf u} (x (t)) + \int_{0}^{t}e^{a   \tau
}l (
x( \tau ),x' ( \tau )) d \tau \right) \right)
\end{displaymath}
The function  $  {\bf u}^{\top}:=\alpha   _{  (f,l)} ^{\top }({\bf  u})$ is
called the {\sf
value function associated with ${\bf u}$}. \index{value function}

If $l=0$, the above problem can be written
\begin{displaymath}
  \alpha _{ (f,0)}  ^{\top} ({\bf u})  (x)  \;  := \; \inf_{x ( \cdot ) \in
{\cal
S}_{f} (x)} \left( \sup_{t \geq 0}    e^{a  t}{\bf u} (x (t))  \right)
\end{displaymath}
and if ${\bf u} \equiv 0$, the above problem boils down to
\begin{displaymath}
 \alpha _{ (f,l)} ^{\top} (0) (x) \; := \; \inf_{x ( \cdot ) \in {\cal
S}_{f} (x)} \int_{0}^{ \infty }e^{a   \tau }l ( x( \tau ),x' ( \tau )) d
\tau 
\end{displaymath}

Before investigating further these examples, we begin by characterizing the
epigraph of ${\bf u}^{\top}$:
\begin{Proposition} \label{optvalfctstoptimepbvkthm}
Let us assume that $f$ and $l$ are  continuous with linear growth and that
${\bf  u}  :  {\bf   R}^{n}  \mapsto {\bf   R}_{+}   \cup  \{+\infty \}$ is
nontrivial, non
negative and lower semicontinuous.  

Then the epigraph of $ {\bf u}^{\top}:=\alpha  _{ (f,l)} ^{\top }({\bf u})$
is the
viability  kernel $   \mbox{\rm Viab}_{g} (  {\cal E}p ({\bf u}))  $ of the
epigraph
of ${\bf u}$ under $g$.

Consequently,   the  function  $  {\bf u}^{\top}$  is  characterized as the
smallest of
the  lower semicontinuous functions  ${\bf v}:{\bf   R}^{n} \mapsto {\bf R}
\cup \{+\infty \}$
larger  than  or  equal to  ${\bf u}$  such that  for any $x  \in \mbox{\rm
Dom}({\bf v})$,
there exists a solution $x ( \cdot )$ to the differential equation $x'=f
(x)$ satisfying property (\ref{dynineeq}):
\begin{displaymath}
  \forall \;  t \geq  0,  \;   \; e^{a  t}{\bf v} (x (t))+ \int_{0}^{t}e^{a
\tau }l
(x ( \tau ),x' ( \tau ))d \tau  \; \leq  \; {\bf v} (x)
\end{displaymath}
\end{Proposition} 

{\bf Proof} --- \hspace{ 2 mm}
Indeed, to say that a pair $ (x,y)$ belongs to the viability kernel $ 
\mbox{\rm Viab}_{g} (  {\cal  E}p ({\bf u}))  $  means that there  exists a
solution
$  x ( \cdot ) \in  {\cal S}_{ f} (x)$ such that
  \begin{displaymath}
\forall \; t \geq 0, \;  \;\left( x (t) , e^{-at}y-\int_{0}^{t}e^{-a(t -
\tau }l (  x( \tau ),x' ( \tau )) d \tau \right) \; \in  \; {\cal E}p ({\bf
u})
\end{displaymath}
i.e., if and only if
\begin{displaymath}
\forall \;  t \geq 0,  \;   \;  e^{a  t} {\bf u} (x (t)) + \int_{0}^{t}e^{a
\tau }l
( x( \tau ),x' ( \tau )) d \tau  \; \leq  \; y
\end{displaymath}
This implies that
\begin{displaymath}
 {\bf u}^{\top} (x) \; := \;  \inf_{x ( \cdot ) \in {\cal S}_{f} (x)}
\left(  \sup_{t \geq 0} \left(    {\bf u} (x (t)) + \int_{0}^{t}l ( x( \tau
),x' (
\tau )) d \tau \right) \right) \; \leq \; y  
\end{displaymath}
and thus,  that $  \mbox{\rm Viab}_{g} ( {\cal E}p ({\bf u}))$ is contained
in $
{\cal E}p ( {\bf u}^{\top } )$.

Since the set $ {\cal S}_{f} (x)$ of solutions is compact in the space $
{\cal C} (0, \infty ;{\bf  R}^{n})$ thanks to Theorem~\ref{03A44} and since
the
function $x ( \cdot ) \mapsto \int_{0}^{t} l (x ( \tau ),f (x ( \tau )))$
is continuous on $ {\cal C} (0, \infty ;{\bf  R}^{n})$,  the infimum 
\begin{displaymath}
  {\bf u}^{\top }(x)  \; := \;    \sup_{t \geq 0} \left( e^{a   t}{\bf u} (
\bar{x} ( t))
+ \int_{0}^{ t} e^{a   \tau }l ( \bar{x}( \tau ), \bar{x}' ( \tau )) d
\tau  \right)
\end{displaymath}
is reached by  a solution $ \bar{x} ( \cdot ) \in  {\cal S}_{f} (x)$.
Consequently, the function
\begin{displaymath}
\left(   \bar{x} (t), e^{-at}{\bf u}^{\top } (x)-\int_{0}^{t}e^{-a(t - \tau
}l (
\bar{x}( \tau ), \bar{x}' ( \tau )) d \tau \right) \; \in  \; {\cal E}p
({\bf u}^{\top })
\end{displaymath}
is viable in the epigraph of ${\bf u} ^{\top }$. Therefore, 
$  (x,   {\bf  u}  ^{\top  } (x))$  belongs to the viability kernel  of the
epigraph of
${\bf u}$.

By Lemma~\ref{tranprplemm}, the epigraph of ${\bf u} ^{\top }$ being the
viability kernel, contains the epigraph of any  lower semicontinuous
functions ${\bf v}:{\bf   R}^{n} \mapsto {\bf R} \cup  \{+\infty \}$ larger
than or equal to
${\bf u}$ viable under $g$, i.e., satisfying property (\ref{dynineeq}).
Therefore, the  function $ {\bf u}^{\top }$ is the smallest of the lower
semicontinuous  functions  ${\bf  v}:{\bf   R}^{n}  \mapsto  {\bf  R}  \cup
\{+\infty \}$ larger
than or equal to ${\bf u}$ satisfying  property (\ref{dynineeq}).
$\; \; \Box$

\begin{Theorem} \label{optvalfctstoptimepbvkthm02}
We posit the assumptions of Theorem~\ref{optvalfctstoptimepbvkthm}.
Then the value function $ {\bf u}^{\top }$ is characterized as the smallest
of the nonnegative lower  semicontinuous  functions  ${\bf  v}:{\bf  R}^{n}
\mapsto {\bf R}
\cup \{+\infty \}$ satisfying for every $x$ 
\begin{displaymath} \left\{ \begin{array}{ll}
i) &  {\bf u} (x) \; \leq  \; {\bf v}(x)\\
ii) &   D _{\uparrow }{\bf v} (x) (f (x))+l (x,f (x))  +a{\bf v} (x)\; \leq
\; 0
\end{array} \right. \end{displaymath}
Furthermore, it satisfies the property
 \begin{displaymath} \left\{ \begin{array}{l}
      \forall \;  x   \;  \mbox{  such  that} \;  {\bf u} (x)  \; < \; {\bf
u}^{\top } (x),\\
D _{\uparrow }{\bf u}^{\top }({\bf u})   (x) (-f (x))-l (x,f (x))  - a {\bf
u}^{\top } (x)\;
\leq  \; 0
 \end{array} \right. \end{displaymath} 
\end{Theorem} 
{\bf Remark} --- \hspace{ 2 mm}
If  the function ${\bf u} ^{\top }$ is differentiable, then the contingent
epiderivative coincides with the usual derivatives, so that ${\bf u} ^{\top
}$ is
a solution to the linear Hamilton-Jacobi ``differential variational
inequalities''
\begin{displaymath} \left\{ \begin{array}{ll}
        i) &  {\bf u} (x) \; \leq  \; {\bf u} ^{\top } (x)\\
         ii)  &   \displaystyle{\left\langle    \frac{ \partial }{ \partial
x}{\bf u}
^{\top } (x),f (x)\right\rangle +l (x,f (x))  +a{\bf u} ^{\top } (x)\; \leq
\;
0}\\
         iii)  &  \displaystyle{({\bf  u} (x)-{\bf u}  ^{\top } (x)) \left(
\left\langle  
\frac{ \partial }{  \partial x}{\bf u}  ^{\top }  (x),f (x)\right\rangle +l
(x,f
(x))  +a{\bf u} ^{\top } (x) \right) \; = \; 0}
        \end{array} \right. \end{displaymath} 

\mbox{}

     {\bf Proof} --- \hspace{ 2 mm}
      By the Nagumo Theorem,  the epigraph of ${\bf v}$ is viable under $g$
if and
only if
     \begin{displaymath}
      \forall \;   (x,y) \in {\cal E}p ({\bf v}), \;  \; (f (x), -ay-l (x,f
(x))) 
\; \in  \; T_{ {\cal E}p ({\bf v})} (x,y)
     \end{displaymath}
     When $y={\bf v} (x)$, we deduce from the fact that the contingent cone
to
the epigraph of ${\bf v}$ at $ (x,{\bf v} (x))$
     \begin{displaymath}
       T_{  {\cal E}p ({\bf v})}  (x,{\bf v}  (x))  \;  := \;  {\cal E}p (D
_{\uparrow }{\bf v}
(x))
     \end{displaymath}
     that if the epigraph of ${\bf v}$ is viable under $g$, then $$ \forall
\; x
\in \mbox{\rm Dom} ({\bf v}),  \;    \;D _{\uparrow }{\bf v} (x)  (f (x))+l
(x,f (x))  +a{\bf v}
(x)\; \leq  \; 0$$
      Conversely,  this inequality implies that   $ (f (x), -a{\bf v} (x)-l
(x,f
(x)))  \in T_{ {\cal E}p ({\bf v})}  (x,{\bf v}  (x))$.  It also
implies that if
$y >{\bf v} (x)$,  $  (f (x),  -ay-l (x,f (x)))$ belongs to  $T_{ {\cal E}p
({\bf v})}
(x,y)$. Indeed, we know that there exist sequences $h_{n}>0$ converging to
$0$, $v_{n}$ converging to $f (x)$ and $ \varepsilon _{n}$ converging to
$0$ such that
     \begin{displaymath}
      (x+h_{n}v_{n}, {\bf v} (x) - h_{n}(a{\bf v} (x)+l (x,f (x)) + h_{n}
\varepsilon_{n}) \; \in  \; {\cal E}p ({\bf v})
     \end{displaymath}
     We thus deduce that
      \begin{displaymath} \left\{ \begin{array}{l}
      (x+h_{n}v_{n}, y - h_{n}(ay+l (x,f (x)) + h_{n} \varepsilon_{n}) \\ =
\;  (x+h_{n}v_{n},  {\bf v}  (x)  -  h_{n}(a{\bf v} (x)+l (x,f (x)) + h_{n}
\varepsilon_{n})
+  (0,  (1-h_{n})(y-{\bf  v}  (x))  )\\ \in   \; {\cal E}p ({\bf v})+ \{0\}
\times {\bf 
R}_{+} \; = \; {\cal E}p ({\bf v})     \\
      \end{array} \right. \end{displaymath} 
and thus, that $ (f (x), -ay-l (x,f (x)))$ belongs to  $T_{ {\cal E}p ({\bf
v})}
(x,y)$.

Finally,  Theorem~\ref{viablbascharthm}  states  that  ${\bf  u}  ^{\top }$
satisfies 
\begin{equation} \label{locbacinv}
      \forall \,   (x,{\bf u} ^{\top }(x)) \in {\cal E}p ({\bf u} ^{\top })
\backslash
\partial \left(  {\cal E}p ({\bf u})  \right),  \;  (-f (x),  a{\bf u}
^{\top } (x) +l
(x,f (x)))   \, \in  \, T_{ {\cal E}p ({\bf u} ^{\top })} (x,{\bf u} ^{\top
} (x))
     \end{equation}
     which, joined to the other properties,  can be translated as the
``Frankowska'' solution to the variational inequalities.
$\; \; \Box$ 
\section{Lyapunov Functions}
Consider a differential equation $x'= f(x)$ and a nontrivial nonnegative
lower semicontinuous extended function  ${\bf u}:{\bf   R}^{n} \mapsto {\bf
R}_{+} \cup
\{+\infty \}$.

The function $\alpha _{ (f,0)} ^{\top} ({\bf u}):{\bf  R}^{n}  \mapsto {\bf
R}_{+} \cup
\{+\infty \}$  defined by 
\begin{displaymath}
  \alpha _{  (f,0)} ^{\top} ({\bf  u})  (x)  \; := \; \inf_{x ( \cdot ) \in
{\cal
S}_{f} (x)}  \left( \sup_{t \geq 0}   e^{a  t}{\bf u} (x (t)) \right)
\end{displaymath}
        is said {\em to enjoy the  $a$-Lyapunov property\/} because for any
initial state $x_{0}$, there exists a solution to the differential equation
$x'=f (x)$  satisfying
        \begin{equation} \label{lyap03A82}
         \forall \;   t  \geq 0,   \;\;  {\bf u} (x(t))  \; \leq   \; e^{-a
t}\alpha
_{ (f,0)} ^{\top} ({\bf u}) (x_{0})
        \end{equation}
        Such inequalities allow us to deduce many properties on the
asymptotic behavior of ${\bf v}$ along the solutions to the differential
equation. 
        This may be quite useful when ${\bf u}$ is the distance function
$d_{M}(\cdot)$ to a subset. The domain of this Lyapunov function $\alpha _{
(f,0)} ^{\top} (d_{M})
        $ provides the $a$-basin of attraction of $M$, which is the set of
states   from which a solution $x ( \cdot )$ to the differential equation
converges exponentially to $M$:
        \begin{displaymath}
        \forall \; x_{0} \in \mbox{\rm Dom}(\alpha _{ (f,0)} ^{\top}
(d_{M})), \;  \; d_{M}(x(t)) \; \leq  \; e^{-a  t} \alpha _{ (f,0)} ^{\top}
(d_{M}) (x_{0})
        \end{displaymath}
        
        The main question we face is {\em to characterize this Lyapunov
function\/}. Ever since Lyapunov proposed in 1892 his second method for
studying the behavior of a solution around an equilibrium, finding Lyapunov
functions for such and such differential equation  has been a source of
numerous problems requiring most often  many  clever tricks. 
        
        We deduce from Theorems~\ref{optvalfctstoptimepbvkthm} and
\ref{optvalfctstoptimepbvkthm02} with $l=0$ the following characterization
of Lyapunov functions:
        
        \begin{Theorem} \label{optvalfctstoptimepbvkthm023}
Let us assume that $f$ is continuous with linear growth and that ${\bf u} :
{\bf  R}^{n}
\mapsto {\bf  R}_{+}  \cup \{+\infty \}$ is  nontrivial, non negative and
lower semicontinuous.  
Then   the epigraph  of  the Lyapunov  function  $  \alpha _{ (f,0)} ^{\top
}({\bf u})$
is the viability kernel $   \mbox{\rm Viab}_{g} ( {\cal E}p ({\bf u})) $ of
the
epigraph of ${\bf u}$ under $g$. 
          Therefore, the Lyapunov function $ \alpha _{ (f,0)} ^{\top }({\bf
u})$
is   the smallest of  the nonnegative lower semicontinuous  functions ${\bf
v}:{\bf  R}^{n}
\mapsto {\bf R} \cup \{+\infty \}$ enjoying the $a$-Lyapunov property,
i.e., such that from any $x_{0}\in \mbox{\rm Dom}({\bf v})$ starts at least
one
solution to the differential equation $x'=f (x)$ 
         satisfying 
        \begin{displaymath}
          \forall \;  t \geq 0,  \;  \; {\bf v} (x (t)) \; \leq  \; {\bf v}
(x_{0})e^{-at}
        \end{displaymath}
        or, equivalently, 
        \begin{displaymath}
         D _{\uparrow }{\bf v} (x) (f (x)) +a{\bf v} (x)\; \leq  \; 0
        \end{displaymath}
Furthermore,  if ${\bf u} (x) < \alpha _{ (f,0)} ^{\top }({\bf u}) (x)$, it
satisfies
\begin{displaymath}
D _{\uparrow }\alpha _{ (f,0)} ^{\top }({\bf u})   (x) (-f (x))  - a \alpha
_{
(f,0)} ^{\top }({\bf u}) (x)\; \leq  \; 0 
\end{displaymath}
\end{Theorem}

\section{Finite Length Solutions}
 
We define now $l (x,p):= \|p\|$, so that $l (x,f (x))= \|f (x)\|$, and take
$a:=0$ and ${\bf u} (x):=0$. Then 
\begin{displaymath}
  \alpha  _{ (f, \|f\|)} ^{\top }(0) (x_{0}) \; = \; \inf_{x ( \cdot ) \in 
{\cal S}_{f} (x)} \int_{0}^{+ \infty } \|x' ( \tau )\| d \tau 
\end{displaymath}
is the minimal length  of the trajectories of the solutions $x ( \cdot )$
to the differential equation $x'=f (x)$ starting from $x_{0}$.

Its epigraph is the viability kernel of ${\bf  R}^{n} \times {\bf  R}_{+} $
under the
system of differential equations $ (x',y') = (f (x), - \|f (x)\|)$. The
minimal length is  the smallest of the nonnegative lower semicontinuous
functions  ${\bf  v}:{\bf   R}^{n}  \mapsto  {\bf  R}  \cup  \{+\infty  \}$
satisfying for every $x$ 
\begin{displaymath}
 D _{\uparrow }{\bf v} (x) (f (x)) + \|f (x)\| \; \leq  \; 0
\end{displaymath}
and satisfies whenever the length $\alpha  _{ f, \|f\|} ^{\top }(0) (x)>0$
is strictly positive
\begin{displaymath}
 D _{\uparrow }\alpha  _{ f, \|f\|} ^{\top }(0)  (x) (-f (x))  - \|f (x)\|
\leq  \; 0
\end{displaymath}

\section{Stopping Time Problem}

We still consider  the   ${\bf u}:{\bf   R}^{n}  \mapsto  {\bf  R}_{+} \cup
\{+\infty \}$,
regarded as an ``obstacle'' in problems of unilateral mechanics. We
associate with it the stopping time problem
\begin{displaymath}
  {\bf u}^{\bot} (x) \; := \; \alpha _{ (f,l)} ^{\bot } ({\bf u}) (x) \; :=
\;  \inf_{x
(  \cdot  )  \in  {\cal S}_{f}  (x)}  \left(  \inf_{t \geq 0} \left(   e^{a
t}{\bf u}
(x (t)) + \int_{0}^{t}e^{a   \tau }l ( x( \tau ),x' ( \tau )) d \tau
\right) \right)
\end{displaymath}

We begin by characterizing its epigraph:
\begin{Proposition} \label{optvalfctstoptimepbthm}
Let us assume that $f$ and $l$ are continuous with linear growth and that
${\bf  u}  :  {\bf   R}^{n}  \mapsto {\bf   R}_{+}   \cup  \{+\infty \}$ is
nontrivial, non
negative and lower semicontinuous.  

Then the epigraph of $ {\bf u}^{\bot}:=\alpha  _{ (f,l)} ^{\bot }({\bf u})$
is the
capture  basin $   \mbox{\rm Capt}_{g} (  {\cal  E}p ({\bf  u}))  $  of the
epigraph of
${\bf u}$ under $g$.
\end{Proposition} 

{\bf Proof} --- \hspace{ 2 mm}
To say that a pair $ (x,y)$ belongs to the capture basin $  \mbox{\rm
Capt}_{g} (  {\cal E}p ({\bf u})) $ means that there exist a solution $ ( x
(
\cdot )) \in  {\cal S}_{ (f)} (x)$ and $t \geq 0$ such that
\begin{displaymath}
\left( x (t) , e^{-at}y-\int_{0}^{t}e^{-a(t - \tau }l ( x( \tau ),x' ( \tau
)) d \tau \right) \; \in  \; {\cal E}p ({\bf u})
\end{displaymath}
i.e., if and only if
\begin{displaymath}
e^{a  t} {\bf u} (x (t)) + \int_{0}^{t}e^{a   \tau }l ( x( \tau ),x' ( \tau
)) d
\tau  \; \leq  \; y
\end{displaymath}
This implies that
 \begin{displaymath} \left\{ \begin{array}{l}
  {\bf  u}^{\bot}(x)  \\  \displaystyle{:=  \;  \inf_{x ( \cdot ) \in {\cal
S}_{f} (x)}
\left( \inf_{t \geq 0} \left(  e^{a  t} {\bf u} (x (t)) + \int_{0}^{t} e^{a

\tau }l ( x( \tau ),x' ( \tau )) d \tau \right) \right) \; \leq \; y}    
\\
 \end{array} \right. \end{displaymath} 
and thus,  that $  \mbox{\rm Capt}_{g} ( {\cal E}p ({\bf u}))$ is contained
in $
{\cal E}p ( {\bf u}^{\bot } )$.

Since   the infimum 
\begin{displaymath}
  {\bf  u}^{\bot} (x)  \;  :=  \;      e^{a    \bar{t}}{\bf u}  ( \bar{x} (
\bar{t})) +
\int_{0}^{ \bar{t}}e^{a   \tau }l ( \bar{x}( \tau ), \bar{x}' ( \tau )) d
\tau 
\end{displaymath}
is reached by a solution $ \bar{x} ( \cdot ) \in  {\cal S}_{f} (x)$ at a
time $ \bar{t}$, this states that $ (x,  {\bf u}^{\bot}(x))$ belongs to the
capture basin of the epigraph of ${\bf u}$.
$\; \; \Box$ 

\mbox{}

In order to apply the properties of capture basins, we need to check that
$X \times {\bf  R}_{+} $  is a repeller under the auxiliary system $ g $.

\begin{Lemma} \label{awlem01}
Let us  assume that there exist real constants $ \gamma _{-}$ and
$ \delta ^{-}$ such that
\begin{displaymath} \left\{ \begin{array}{ll} \label{aweq02}
i) & \inf_{x \in X} \frac{  \langle x,f (x) \rangle }{ \|x\|} \;
\geq  \; \gamma_{-} ( \|x\|+1)\\
ii) & \inf_{x \in X}  l (x,f (x))\; \geq  \; \delta _{-} (
\|x\|+1)
\end{array} \right. \end{displaymath}
and if  $ a+ \gamma _{-} >0$, then $X \times {\bf  R}_{+} $ is a
repeller.
\end{Lemma}

{\bf Proof} --- \hspace{ 2 mm}
Whenever $l$ is nonnegative, the backward solutions  $ ( x ( \cdot
),y ( \cdot ))$ starting from $X \times  {\bf  R}_{+} $ are
viable in $X \times  {\bf  R}_{+}$ because $y' (t)= -ay (t)+l (x
(t)) \geq 0$.

Let $ (x ( \cdot ),y ( \cdot ))$ be the solution to the
differential equation $ (x',y') = g  (x,y)$ starting from $
(x_{0},y_{0})$. 

Therefore
\begin{displaymath}
 \frac{d}{dt} \|x (t)\| \; = \; \left\langle   x' (t), \frac{x
(t)}{ \|x (t)\|}\right\rangle \; = \; \left\langle   f(x (t))),
\frac{x (t)}{ \|x (t)\|}\right\rangle \; \geq  \; \gamma _{-} (
\|x (t)\|+1)
\end{displaymath}
so that
\begin{displaymath}
 \forall \; t \geq 0, \;   \; \|x (t) \|  \; \geq  \; e^{ \gamma
_{-}t} ( \|x_{0}\|+1) -1
\end{displaymath}
Furthermore, since
\begin{displaymath}
 l (x ( \tau ),x' ( \tau )) \; \geq  \; \delta_{-} ( \|x ( \tau
)\|+1) \; \geq  \; \delta_{-} ( \|x_{0}\|+1)e^{ \gamma _{-}t}
\end{displaymath}
and since
\begin{displaymath}
 e^{at} y (t) \; = \; y_{0} - \int_{0}^{t}e^{a \tau } l (x ( \tau
))d \tau 
\end{displaymath}
we infer that
\begin{displaymath}
 e^{at} y (t) \; \leq  \; y_{0} - \delta_{-} (
\|x_{0}\|+1)\int_{0}^{t}e^{ \gamma_{-}+ a)\tau } \; = \; y_{0} -
\frac{\delta_{-} ( \|x_{0}\|+1)}{ \gamma_{-}+a} \left(e^{
(\gamma_{-}+a)t}-1  \right)
\end{displaymath}
Consequently, if $  \gamma {-}+ a>0$ 
\begin{displaymath}
 e^{at}y (t) \; \leq  \; y_{0}+\frac{\delta_{-} ( \|x_{0}\|+1)}{
\gamma_{-}+a} - \frac{\delta_{-} ( \|x_{0}\|+1)}{ \gamma_{-}+a}
e^{ (\gamma _{-}-a)t}
\end{displaymath}
so that $y (t)$ becomes negative in finite time.
$\; \; \Box$ 

\mbox{}

\begin{Theorem} \label{optvalfctstoptimepbthm28}
We posit the assumptions of Proposition~\ref{optvalfctstoptimepbthm} and
we assume that
\begin{displaymath} 
 \forall \; x \in K, \;  \;\forall \; x ( \cdot ) \in  {\cal S}_{f} (x),
\;  \;  \int_{0}^{+ \infty }e^{a   \tau }l (x ( \tau ),x' ( \tau ))d \tau 
\; = \; + \infty 
\end{displaymath}

 Then  ${\bf u}^{\bot }$ is characterized as the {\bf unique} nonnegative
lower semicontinuous functions ${\bf v}:{\bf   R}^{n} \mapsto  {\bf R} \cup
\{+\infty \}$
such that  from  any $x  $  satisfying ${\bf v}  (x)<{\bf u} (x)$  starts a
solution $x
( \cdot ) \in  {\cal S}_{f} (x)$ satisfying, for some time $T>0$
\begin{displaymath}
 \forall \; t \in [0,T], \;  \; e^{a  t} {\bf v} (x (t)) +
\int_{0}^{t}e^{a   \tau }l ( x( \tau ),x' ( \tau )) d \tau  \; \leq 
\; {\bf v} (x)
\end{displaymath}
and that, for any $T>0$ and any $x_{T} \in \mbox{\rm Dom}({\bf u})$, all
solutions $x ( \cdot )$ to the differential equation $x'=f (x)$
arriving at $x_{T}$ at time $T$ satisfy
\begin{displaymath}
 \forall \; t \in [0,T], \;  \; e^{a  t} {\bf v} (x (t)) +
\int_{0}^{t}e^{a   \tau }l ( x( \tau ),x' ( \tau )) d \tau  \; \leq 
\; {\bf v} (x ( 0))
\end{displaymath}
The function ${\bf u} ^{\bot }$ is also the smallest of the lower
semicontinuous functions ${\bf v}$ satisfying
\begin{displaymath} \left\{ \begin{array}{ll}
i) & 0 \; \leq  \; {\bf v} (x) \; \leq  \; {\bf u} (x)\\
ii)  &  \mbox{\rm  if} \;  0  \leq {\bf  v}  (x)  <  {\bf u} (x), \; \;   D
_{\uparrow }{\bf v} (x) (f
(x))+l (x,f (x))  +a{\bf v} (x)\; \leq  \; 0
\end{array} \right. \end{displaymath}
If we assume furthermore that $f$ and $l$ are Lipschitz, then the  function
${\bf u} ^{\bot  }$  is the {\bf unique} solution  ${\bf v} \geq 0$  to the
system of
``differential inequalities'': for every $x \in \mbox{\rm Dom}({\bf v})$,
\begin{displaymath} \left\{ \begin{array}{ll}
i) & 0 \; \leq  \;{\bf v} (x) \; \leq  \; {\bf u} (x)\\
ii) & \mbox{\rm if} \; 0 \leq   {\bf v} (x) < {\bf u} (x) , \;  \;
 D _{\uparrow }{\bf v} (x) (f (x))+l (x) + a{\bf v} (x) \; \leq  \; 0\\
iii)  &   \mbox{\rm if} \; 0 \leq {\bf v} (x) \leq  {\bf u} (x) , \;  \;  D
_{\uparrow }{\bf v}
(x) (-f (x) -l (x),f (x)) - a {\bf v}(x)\; \leq  \;0\\
\end{array} \right. \end{displaymath}

Knowing the function  ${\bf u}  ^{\bot }$,  the stopping time is  the first
time $
\bar{t}\geq  0$  when  ${\bf v}  (x  (  \bar{t}))=   {\bf u} ^{\bot  }( x (
\bar{t}))$.

 \end{Theorem}

\mbox{}

{\bf Remark} --- \hspace{ 2 mm}
If   the function ${\bf u} ^{\bot }:=  \alpha  _{ (f,l)} ^{\bot }({\bf u})$
is
differentiable, then the contingent epiderivative coincides with the usual
derivatives, so that ${\bf u} ^{\bot }$ is a solution to the linear
Hamilton-Jacobi ``differential variational inequalities''
        \begin{displaymath} \left\{ \begin{array}{ll}
        i) &  0 \; \leq  \; {\bf u} ^{\bot }(x) \; \leq  \; {\bf u}  (x)\\
         ii)  &   \displaystyle{\left\langle    \frac{ \partial }{ \partial
x}{\bf u}
^{\bot } (x),f (x)\right\rangle +l (x,f (x))  +a{\bf u} ^{\bot } (x)\; \geq
\;
0}\\
         iii)  &  \displaystyle{({\bf u}  (x)-{\bf u} ^{\bot  } (x)) \left(
\left\langle  
\frac{ \partial  }{ \partial x}{\bf  u} ^{\bot } (x),f  (x)\right\rangle +l
(x,f
(x))  +a{\bf u} ^{\bot } (x) \right) \; = \; 0}
        \end{array} \right. \end{displaymath} 

\mbox{}

     {\bf Proof} --- \hspace{ 2 mm}
      Since the Lagrangian is nonnegative,  the closed subset ${\bf  R}^{n}
\times {\bf 
R}_{+} $ is backward invariant under $g$. It is  a repeller under $g$
whenever 
\begin{displaymath}
 \forall \; x \in {\bf  R}^{n}, \; \forall \; x ( \cdot ) \in  {\cal S}_{f}
(x), \; \;
\int_{0}^{+ \infty }e^{a   \tau }l (x ( \tau ),x' ( \tau ))d \tau  \; = \;
+ \infty 
\end{displaymath}
      By Theorem~\ref{wonderfulthmthmbis}, $ {\cal E}p ({\bf u} ^{\bot }):=
\mbox{\rm Capt}_{g} (  {\cal E}p ({\bf u}))$ is  the unique closed subset $
{\cal
{\bf v}}$,  and in particular, the unique epigraph ${\cal E}p ({\bf v})$ of
lower
semicontinuous function ${\bf v}$, which satisfies
\begin{displaymath} \left\{ \begin{array}{ll}
i)  &  {\cal E}p ({\bf u}) \; \subset  \; {\cal E}p ({\bf v})\; \subset  \;
{\bf  R}^{n} \times
{\bf  R}_{+} \\
ii) &  {\cal E}p ({\bf v}) \backslash {\cal E}p ({\bf u}) \;\; \mbox{\rm is
locally viable under $g$}\\
iii)  &  {\cal  E}p ({\bf v})  \;\;  \mbox{\rm is backward  invariant under
$g$}\\
\end{array} \right. \end{displaymath}
The first condition means that for any $x  \in {\bf   R}^{n}$, $0 \leq {\bf
v} (x) \leq {\bf u}
(x)$ and we observe that $(x,y) \in  {\cal E}p ({\bf v} ) \backslash {\cal
E}p ({\bf u})$ if and only if $y  \in [{\bf v} (x),{\bf u} (x)[$.
Hence, the first statement of the theorem ensues. 

By Theorem~\ref{viablcaptbascharthm}, $ {\cal E}p ({\bf u} ^{\bot }):=
\mbox{\rm Capt}_{g} ( {\cal E}p ({\bf u}))$ is also the smallest of the
nonnegative lower semicontinuous ${\bf v}$ satisfying 
     \begin{displaymath} \left\{ \begin{array}{ll}
i) &  {\cal E}p ({\bf u}) \; \subset  \; {\cal E}p ({\bf v}) \; \subset  \;
{\bf  R}^{n} \times
{\bf  R}_{+} \\
ii)  &  \forall \;   (x,y)  \in   {\cal E}p ({\bf v} ) \backslash {\cal E}p
({\bf u}), \; 
\;  (f (x),-l (x,f (x))-ay) \; \in \; T_{ {\cal E}p ({\bf v})} (x,y) 
\end{array} \right. \end{displaymath}

     When $y={\bf v} (x)$, the second condition can be written
     \begin{displaymath}
      \mbox{\rm if} \;  {\bf v} (x) <{\bf u} (x),  \; \; D _{\uparrow }{\bf
v} (x)+a{\bf v} (x)+ l
(x,f (x)) \; \leq  \; 0
     \end{displaymath}
      Conversely,  this  condition implies that  for any  $y   \in ]{\bf v}
(x),{\bf u}
(x)[$,  $(f (x),-l (x,f (x))-ay)  $ also belongs to  $  T_{ {\cal E}p ({\bf
v})}
(x,y)$ as in the proof of Theorem~\ref{optvalfctstoptimepbvkthm02}.
     
     Since $ {\cal E}p ({\bf u} ^{\bot }) = \mbox{\rm Capt}_{g} ( {\cal E}p
({\bf u}))$
is the unique closed subset satisfying the above properties and being
backward invariant under $g$,  this implies that
     \begin{displaymath}
      \mbox{\rm if} \; {\bf v} (x) \leq  {\bf u} (x) , \;  \;  D _{\uparrow
}{\bf v} (x) (-f
(x) -l (x),f (x)) - a {\bf v}(x)\; \leq  \;0
     \end{displaymath}
     The converse is true when $f$ and $l$ are Lipschitz, or whenever the
solution to the system $ (x',y')=g (x,y)$ is unique. In this case, the
function  ${\bf u}  ^{\bot  }$  is  the unique solution satisfying  the two
properties.
$\; \; \Box$ 
\section{Minimal Time and Minimal Length Solutions}

Let us consider a closed  subset $K \subset {\bf   R}^{n}$  and $ \psi_{K}$
its indicator and take $a=0$.

We observe that the hitting time (or minimal time) function $ \omega_{K}^{f
^{\flat }}$ is equal to \index{minimal time function}
\begin{displaymath}
 \omega_{K}^{f ^{\flat }} \; = \;  \alpha_{f,1}^{\bot} ( \psi_{K})
\end{displaymath}

In  the  same  way,  we  introduce  the  {\em  minimal  length  functional}
associating with $x ( \cdot )$
$$   \lambda_{K}(x  (  \cdot  ))  \;  :=  \;   \inf_{ \{t \;  | \;x (t) \in
K\}}\int_{0}^{t} \|x' (s)\| ds$$  
  the minimal length of the  curve $s \mapsto  x (s)$  from $0$ to $t$ such
that $x  (t)  \in  K$.  We  next  define  the  ``minimal  length'' function
\index{minimal length function} $ \lambda_{K}^{f ^{\flat }}$ by
\begin{displaymath}
 \lambda_{K}^{f ^{\flat }} (x) \; := \; \inf_{x ( \cdot ) \in  {\cal S}_{f}
(x)}\lambda_{K}(x ( \cdot ))
\end{displaymath}

We note that
\begin{displaymath}
 \lambda _{K}^{f ^{\flat }} \; = \;  \alpha_{f, \|f\|}^{\bot} ( \psi_{K})
\end{displaymath}

So, these two functions enjoy the properties proved above. For instance:
\begin{enumerate}
\item The minimal time function $ \omega_{K}^{f ^{\flat }}$ is the smallest
nonnegative lower semicontinuous ${\bf  v}$  function vanishing on $K$ such
that,
\begin{displaymath}
 \forall \; x \notin K, \;\; D _{\uparrow } {\bf v}  (x) (f (x)) +1 \; \leq
\;0
\end{displaymath}
If $f$  is assumed furthermore   to  be  Lipschitz,  it is the {\bf unique}
nonnegative lower semicontinuous solution  vanishing on  $K$ satisfying the
above inequalities and
\begin{displaymath}
 \forall \; x \in {\bf  R}^{n}, \;  \; D _{\uparrow } {\bf v}  (x) (-f (x))
- 1 \; \leq  \;0
\end{displaymath}
\item Assume that 
\begin{displaymath}
  \forall \;  x \in {\bf  R}^{n}, \;  \; \inf_{x ( \cdot ) \in {\cal S}_{f}
(x)} \int_{0}^{+ \infty } \|x' (s)\|ds \; = \; + \infty 
\end{displaymath}
Then the minimal  length  function  $  \lambda_{K}^{f  ^{\flat  }}$  is the
smallest nonnegative lower semicontinuous  ${\bf v}$  function vanishing on
$K$ such that,
\begin{displaymath}
  \forall \;  x \notin K,  \;\;  D _{\uparrow } {\bf v}   (x) (f (x)) + \|f
(x)\| \; \leq  \;0
\end{displaymath}
If $f$  is assumed furthermore   to  be  Lipschitz,  it is the {\bf unique}
nonnegative lower semicontinuous solution  vanishing on  $K$ satisfying the
above inequalities and
\begin{displaymath}
 \forall \; x \in {\bf  R}^{n}, \;  \; D _{\uparrow } {\bf v}  (x) (-f (x))
- \|f (x)\| \; \leq  \;0
\end{displaymath}

\end{enumerate}

Since  the minimal time  and minimal  length  functions  coincide  with the
indicator  of $  \psi_{K}$  on $K$,  the above conditions imply that $K$ is
backward invariant under $f$ whenever  $f$ is Lipschitz.

\section{Viscosity Type Solutions}

We now use the characterizations in terms of normal cones for deriving the
formulations in terms of subgradients instead of contingent epiderivatives.

\begin{Theorem} \label{optvalfctstoptimepbvkthm02gg}
We posit the assumptions of Theorem~\ref{optvalfctstoptimepbvkthm}.
Then the value function $ {\bf u}^{\top }$ is the solution to 
\begin{displaymath} \left\{ \begin{array}{ll}
        i) &  {\bf u} (x) \; \leq  \; {\bf u} ^{\top } (x)\\
        ii) &  \displaystyle{ \forall \; p \in  \partial _{-}{\bf u}^{\top}
(x),
\;  \; \left\langle   p,f (x)\right\rangle +l (x,f (x))  +a{\bf u} ^{\top }
(x)\;
\leq  \; 0}\\ & \mbox{\rm and} \\
        & \displaystyle{\forall \; p  \in \mbox{\rm Dom}(D _{\uparrow
}{\bf u}^{\top} (x))^{-}, \;  \;  \langle p,f (x) \rangle \; \leq 0 \;}\\
         iii) & \displaystyle{\forall \; p \in  \partial _{-}{\bf u}^{\top}
(x),
\;   \;({\bf  u}  (x)-{\bf  u}  ^{\top  }  (x))  \left(  \left\langle   p,f
(x)\right\rangle +l
(x,f (x))  +a{\bf u} ^{\top } (x) \right) \; = \; 0}\\
        & \mbox{\rm and} \\
        & \displaystyle{\forall \; p  \in \mbox{\rm Dom}(D _{\uparrow
}{\bf u}^{\top} (x))^{-}, \;  \;  ({\bf u} (x)-{\bf u} ^{\top } (x))\langle
p,f (x) \rangle
\; = 0 \;}\\
        \end{array} \right. \end{displaymath} 
\end{Theorem} 
Such a solution, recently discovered independently by Frankowska and
Barron \& Jensen,  are sometime called ``bilateral solutions'' to
Hamilton-Jacobi equation 
\begin{displaymath}
 \left\langle  \frac{ \partial }{ \partial x}{\bf u} (x),f (x)\right\rangle
+l
(x,f (x))  +a{\bf u}  (x)  \; = \; 0
\end{displaymath}
The method we present here is due to Frankowska.

     {\bf Proof} --- \hspace{ 2 mm}
      By  Theorem~\ref{optvalfctstoptimepbvkthm02},  we  know that whenever
${\bf u}
(x) \leq {\bf u}^{\top} (x)$,
     \begin{displaymath}
      (f (x),-a{\bf u}^{\top} (x)-l (x,f (x))) \; \in  \; T_{ {\cal E}p
({\bf u}^{\top})} (x,{\bf u}^{\top} (x))
     \end{displaymath}
     and that
     whenever ${\bf u} (x) < {\bf u}^{\top} (x)$,
     \begin{displaymath}
      (f (x),-a{\bf u}^{\top} (x)-l (x,f (x))) \; \in  \; T_{ {\cal E}p
({\bf  u}^{\top})} (x,{\bf  u}^{\top}  (x))  \cap  -  T_{  {\cal  E}p ({\bf
u}^{\top})} (x,{\bf u}^{\top}
(x))
     \end{displaymath}
     Theorem~\ref{3russianfrank} implies that these conditions are
equivalent to
     \begin{displaymath}
      \forall \; (p, \lambda ) \in  N_{ {\cal E}p ({\bf u}^{\top})} (x,{\bf
u}^{\top}
(x)), \;  \;  \langle p,f (x) \rangle - \lambda (a{\bf u}^{\top} (x)+l (x,f
(x)))
\; \leq  \; 0
     \end{displaymath}
     whenever ${\bf u} (x) \leq {\bf u}^{\top} (x)$, and to
     \begin{displaymath}
      \forall \; (p, \lambda ) \in  N_{ {\cal E}p ({\bf u}^{\top})} (x,{\bf
u}^{\top}
(x)), \;  \;  \langle p,f (x) \rangle - \lambda (a{\bf u}^{\top} (x)+l (x,f
(x)))
\; =  \; 0
     \end{displaymath}
     whenever ${\bf u} (x) < {\bf u}^{\top} (x)$.

     It remains now to recall that
     \begin{displaymath} \left\{ \begin{array}{ll}
      i)  &  (p,-1)  \in  N_{ {\cal E}p ({\bf u}^{\top})} (x,{\bf u}^{\top}
(x))
\;\mbox{\rm if and only if}\;p \in  \partial_{-}{\bf u}^{\top} (x)\\
      ii)  &  (p,0)  \in  N_{ {\cal E}p ({\bf u}^{\top})} (x,{\bf u}^{\top}
(x))
\;\mbox{\rm  if  and only if}\;   p   \in \mbox{\rm Dom}(D _{\uparrow }{\bf
u}^{\top}
(x))^{-}
     \end{array} \right. \end{displaymath}
     Taking $ \lambda =-1$, we obtain
     \begin{displaymath}
      \forall \;  p \in  \partial _{-}{\bf u}^{\top} (x), \;  \;\langle p,f
(x)
\rangle - \lambda (a{\bf u}^{\top} (x)+l (x,f (x))) \; \leq   \; 0
     \end{displaymath}
     and 
     whenever ${\bf u} (x) < {\bf u}^{\top} (x)$,
     \begin{displaymath}
     \forall \; p \in  \partial _{-}{\bf u}^{\top} (x), \;  \;  \langle p,f
(x)
\rangle - \lambda (a{\bf u}^{\top} (x)+l (x,f (x))) \; =  \; 0
     \end{displaymath}
     Taking $ \lambda =0$ yields that for all $p \in \mbox{\rm Dom}(D
_{\uparrow }{\bf u}^{\top} (x))^{-}$, $  \langle p,f (x) \rangle \leq 0$ if
${\bf u}
(x)  \leq   {\bf u}^{\top} (x)$ and $  \langle p,f (x) \rangle =0$ if ${\bf
u} (x) <
{\bf u}^{\top} (x)$. This means that $f (x)$ belongs to  the closure of
$\mbox{\rm Dom}(D _{\uparrow }{\bf u}^{\top} (x))$  in the general case and
that
$x$ belongs to  the vector space spanned by $\mbox{\rm Dom}(D _{\uparrow
}{\bf u}^{\top} (x))$ when ${\bf u} (x) < {\bf u}^{\top} (x)$. $\; \; \Box$

\mbox{}

We obtain an analogous statement for the function ${\bf u}^{\bot}$:

\begin{Theorem} \label{optvalfctstoptimepbthm28gg}
We posit the assumptions of Theorem~\ref{optvalfctstoptimepbthm28}.
Then the value function $ {\bf u}^{\bot }$ is the solution to 
\begin{displaymath} \left\{ \begin{array}{ll}
        i) &  0 \; \leq  \; {\bf u}^{\bot}(x) \; \leq  \; {\bf u} (x)\\
        ii) &  \displaystyle{ \forall \; p \in  \partial _{-}{\bf u}^{\bot}
(x),
\;  \; \left\langle   p,f (x)\right\rangle +l (x,f (x))  +a{\bf u} ^{\bot }
(x)\;
\geq  \; 0}\\ & \mbox{\rm and} \\
          &   \displaystyle{\forall  \;   p   \in  \mbox{\rm  Dom}(D{\bf u}
_{\downarrow }
^{\top} (x))^{-}, \;  \;  \langle p,f (x) \rangle \; \geq 0 \;}\\
         iii) & \displaystyle{\forall \; p \in  \partial _{-}{\bf u}^{\bot}
(x),
\;   \;({\bf  u}  (x)-{\bf  u}  ^{\bot  }  (x))  \left(  \left\langle   p,f
(x)\right\rangle +l
(x,f (x))  +a{\bf u} ^{\bot } (x) \right) \; = \; 0}\\
        & \mbox{\rm and} \\
        & \displaystyle{\forall \; p  \in \mbox{\rm Dom}(D _{\downarrow }
{\bf u}^{\top} (x))^{-},  \;  \;  ({\bf u} (x)-{\bf u} ^{\top } (x))\langle
p,f (x) \rangle \;
= 0 \;}\\
        \end{array} \right. \end{displaymath} 
\end{Theorem} 
     
     {\bf Proof} --- \hspace{ 2 mm} 
      By  Theorem~\ref{optvalfctstoptimepbthm28},  we  know  that  whenever
${\bf u}
(x) \geq {\bf u}^{\bot} (x)$,
     \begin{displaymath}
      (-f (x),a{\bf u}^{\bot} (x) + l (x,f (x))) \; \in  \; T_{ {\cal E}p
({\bf u}^{\bot})} (x,{\bf u}^{\bot} (x))
     \end{displaymath}
     and that
     whenever ${\bf u} (x) > {\bf u}^{\bot} (x)$,
     \begin{displaymath}
      (f (x),-a{\bf u}^{\bot} (x)-l (x,f (x))) \; \in  \; T_{ {\cal E}p
({\bf u}^{\bot})}  (x,{\bf  u}^{\bot}  (x))  \cap  -  T_{  {\cal  E}p ({\bf
u}^{\bot})} (x,{\bf u}^{\bot}
(x))
     \end{displaymath}
     Theorem~\ref{3russianfrank} implies that these conditions are
equivalent to
     whenever ${\bf u} (x) \geq {\bf u}^{\bot} (x)$,
     \begin{displaymath}
      \forall \; (p, \lambda ) \in  N_{ {\cal E}p ({\bf u}^{\bot})} (x,{\bf
u}^{\bot}
(x)), \;  \;  \langle p,f (x) \rangle - \lambda (a{\bf u}^{\bot} (x)+l (x,f
(x)))
\; \geq  \; 0
     \end{displaymath}
     and that
     whenever ${\bf u} (x) > {\bf u}^{\bot} (x)$,
     \begin{displaymath}
      \forall \; (p, \lambda ) \in  N_{ {\cal E}p ({\bf u}^{\bot})} (x,{\bf
u}^{\bot}
(x)), \;  \;  \langle p,f (x) \rangle - \lambda (a{\bf u}^{\bot} (x)+l (x,f
(x)))
\; =  \; 0
     \end{displaymath}
     It remains now to translate these statements in terms of subgradients.
$\; \; \Box$ 

\mbox{}

{\bf Remark: Viscosity Solutions} --- \hspace{ 2 mm}
When we  know  {\em  a priori\/}  that  the  solution  ${\bf  u}^{\bot}$ is
continuous,  we can prove that it is also a ``viscosity solution''  the the
Hamilton-Jacobi variational inequalities:

\begin{Theorem} \label{optvalfctstoptimepbthm28visc}
We posit the assumptions of  Theorem~\ref{optvalfctstoptimepbthm28}, and we
assume that $f$  and $l$ are Lipschitz and the function ${\bf u}^{\bot}$ is
continuous.
Then the value function $ {\bf u}^{\bot }$ is the solution to 
\begin{displaymath} \left\{ \begin{array}{ll}
        i) &  0 \; \leq  \; {\bf u}^{\bot}(x) \; \leq  \; {\bf u} (x)\\
        ii) &  \displaystyle{ \forall \; p \in  \partial _{+}{\bf u}^{\bot}
(x),
\;   \; \left\langle   -p,f (x)\right\rangle +l (x,f (x))  +a{\bf u} ^{\bot
}
(x)\; \geq  \; 0}\\
        \\ & \mbox{\rm and} \\
        & \displaystyle{\forall \; p  \in \mbox{\rm Dom}(D _{\downarrow }
{\bf u}^{\bot})^{-}, \;  \;  \langle p,f (x) \rangle \; \leq 0 \;}\\
         iii) & \displaystyle{\forall \; p \in  \partial _{-}{\bf u}^{\bot}
(x),
\;   \; \left\langle  p,f (x)\right\rangle +l (x,f (x))  +a{\bf u} ^{\bot }
(x) 
\; \leq  \; 0}\\& \mbox{\rm and} \\
        & \displaystyle{\forall \; p  \in \mbox{\rm Dom}(D _{\uparrow }
{\bf u}^{\bot} (x))^{-}, \;  \;  \langle p,f (x) \rangle \; \leq  0 \;}\\
        \end{array} \right. \end{displaymath} 
\end{Theorem} 
A solution to such a system of inequalities is called a `viscosity
solution' \index{viscosity solution} to Hamilton-Jacobi equation 
\begin{displaymath}
 \left\langle  \frac{ \partial }{ \partial x}{\bf u} (x),f (x)\right\rangle
+l
(x,f (x))  +a{\bf u} ^{\bot } (x)  \; = \; 0
\end{displaymath}
by Michael Crandall and Pierre-Louis Lions. 

\mbox{}

{\bf Proof} --- \hspace{ 2 mm}
We know that epigraph of ${\bf u}^{\bot}$  is backward invariant under $g$.
By
Lemmas~\ref{bacwinvcomp} and \ref{backinvab}, this implies that its
complement is forward invariant under $g$, and, since $f$ and $l$ are
Lipschitz, that its closure is also forward invariant. Since we assumed
that ${\bf u}^{\bot}$ is continuous, the closure of the complement of the
epigraph  of  ${\bf  u}^{\bot}$  is  the  hypograph  of  ${\bf  u}^{\bot}$.
Therefore, 
\begin{displaymath}
 \forall \;  (x,{\bf u} ^{\bot }(x)) \in \mbox{\rm Graph}({\bf u} ^{\bot })
, \;  \; (f
(x),  -a{\bf u} ^{\bot } (x) -l (x,f (x)))  \; \in  \; T_{ {\cal H}yp ({\bf
u} ^{\bot
})} (x,{\bf u} ^{\bot } (x))
\end{displaymath}
and thus,
     \begin{displaymath}
     \forall \; (p, \lambda ) \in  N_{ {\cal H}yp ({\bf u}^{\bot})} (x,{\bf
u}^{\bot}
(x)), \;  \;  \langle p,f (x) \rangle - \lambda (a{\bf u}^{\bot} (x)+l (x,f
(x)))
\; \leq  \; 0
     \end{displaymath}
     It remains now to recall that
     \begin{displaymath} \left\{ \begin{array}{ll}
      i)  &  (p,1)  \in  N_{ {\cal H}yp ({\bf u}^{\bot})} (x,{\bf u}^{\bot}
(x))
\;\mbox{\rm if and only if}\; p \in  \partial_{+}{\bf u}^{\bot} (x)\\
      ii)  &  (p,0) \in  N_{ {\cal H}yp ({\bf u}^{\bot})} (x,{\bf u}^{\bot}
(x))
\;\mbox{\rm if and only if}\; \forall \; p \in \mbox{\rm Dom}(D
_{\downarrow }{\bf u}^{\bot} (x))^{-}
     \end{array} \right. \end{displaymath}
     for deducing that the above condition is equivalent to
     \begin{displaymath}
       \forall  \;  p \in   \partial_{+}{\bf u}^{\bot} (x), \;  \;  \langle
-p,f (x)
\rangle +a{\bf u}^{\bot} (x)+ l (x,f (x)) \; \geq 0 \;
     \end{displaymath}
     and that $f (x)$ belongs to the closure of the domain of $D
_{\downarrow } {\bf u}^{\bot} (x)$ for achieving the proof. $\; \; \Box$

\chapter{Systems of First-Order Partial Differential Equations}

\vspace{ 22 mm}
{\Huge \bf Introduction}

\vspace{ 13 mm}

        We study here Dirichlet boundary value problems for  systems of
first-order partial differential equations of the form
\begin{displaymath}
\forall \;  j=1, \ldots ,p, \;  \;  \frac{ \partial }{ \partial t}u (t,x) +
\sum_{i=1}^{n}\frac{ \partial }{ \partial x_{i}}u_{i}
(t,x) f_{i}(t,x,u (t,x)) - g_{j}(t,x,u (t,x)) \; = \; 0
\end{displaymath}
on $  {\bf   R}_{+} \times  \Omega $, where $ \Omega \subset {\bf  R}^{n} $
is an open subset,  $  \Gamma $  its boundary and $K:=\overline{ \Omega } $
its closure.

 It is known that the solution $U$ to the above system can be set-valued,
describing ``shocks''. 
This is considered as a pathology  whenever the solution is regarded as a
map from  the   input space $  {\bf   R}_{+} \times K$  to the output space
${\bf  R}^{p}$,
but is quite natural when the solution $U$ is considered as a graph, i.e.,
a subset of ${\bf   R}_{+} \times K \times {\bf  R}^{p}$ and when the tools
of set-valued analysis  are used.
Since we are looking for set-valued map solutions, we begin by introducing
\begin{enumerate}
\item the graph of the {\em upper graphical limit} of a sequence of maps
$U_{n}:X\leadsto Y$ (single-valued or set-valued) is the upper limit of the
graphs of $U_{n}$,
\item {\em the contingent derivative} $DU (x,y)$ at  point $ (x,y)$ of the
graph of $U$ is the upper graphical limit of the difference quotients $
\nabla _{h} U (x,y)$, so that the graph of the contingent derivative is the
contingent cone to the graph of $U$:
\begin{displaymath}
 \mbox{\rm Graph}(DU (x,y)) \; = \; T_{ \mbox{\rm Graph}(U) (x,y)}
\end{displaymath}
\end{enumerate}
        
Introducing 
\begin{enumerate}
\item an {\em initial data} $u_{0}: K \mapsto {\bf  R}^{p}$,
\item a {\em boundary data} $v_{ \Gamma}: {\bf  R}_{+} \times  \partial
K \mapsto {\bf  R}^{p}$.
\end{enumerate}
we shall prove the {\em  existence and the  uniqueness\/}  of  a set-valued
solution $:{\bf 
R}_{+}  \times  K   \leadsto  {\bf   R}^{p} $  to the system of first-order
partial
differential equations satisfying the initial/boundary-value conditions

\begin{equation} \left\{ \begin{array}{ll} \label{initboucd001}
i) & \forall \; x \in K, \;  \; u_{0} (x) \; \in \; U(0,x)\\
ii)  &   \forall \;  t>0,  \; \forall \;  x \in  \Gamma, \;  \; v_{ \Gamma}
(t,x) \; \in  \; U (t,x)
\end{array} \right. \end{equation}

Again,   the  strategy  is  the   same  than  the   one   we  followed  for
Hamilton-Jacobi variational inequalities. It is enough to revive the method
of characteristics by  observing that  the graph  of  the  solution  is the
capture  basin of  the graph of  the initial/boundary-value  data  under an
auxiliary system (the ``characteristic system''), use the characterizations
derived from the Nagumo  Theorem and the fact  that the  contingent cone to
the graph is the graph of the contingent derivative of a set-valued map.

However,  the solution becomes single-valued when the  maps $f_{i}$  depend
only on the variables $x$. In this case, we even obtain explicit formulas.

        \section{Contingent Derivatives of Set-Valued Maps}
        \subsection{Set-Valued Maps}

        \begin{Definition} 
        Let  $X $  and  $Y $  be two spaces. A set-valued map $F $  from 
$X $ to  $Y $ is characterized by its \index{graph} {\em graph\/} $Graph(F)
$, the subset of the product space $X \times Y $  defined by
        \begin{displaymath}
        \mbox{\rm Graph}(F)  \;  := \;  \{ (x,y) \in  X \times Y \; |\;\; 
y \in F(x) \}
        \end{displaymath}

        We shall say that $F(x) $  is the \index{image} {\em image\/}  or
the \index{value} {\em value\/}  of $F $  at $x $. 
        
        A set-valued map is said to be \index{nontrivial} {\em
nontrivial\/}  if its graph is not empty, i.e., if there exists at least an
element $x \in X$ such that  $F(x)$ is not empty. 
        
        We say that $F$
        is \index{strict} {\em strict\/}  if all images $F(x)$ are not
empty. The \index{domain} {\em domain\/} of F is the subset  of elements $x
\in   X$  such that $F(x)$ is not empty:
        $$ \mbox{\rm Dom}(F) \;  :=  \; \{x \in X  \;\; |\;\;F(x) \ne
\emptyset \}$$
         The  \index{image} {\em image\/}  of $F$ is the union of the
images (or values) $F(x)$, when $x$ ranges over $X$:
        \begin{displaymath}
        \mbox{ \rm Im}(F) \; := \; \bigcup_{x \in X}F(x) 
        \end{displaymath}  
        The \index{inverse map} {\em inverse\/}    $F^{-1 }$   of $ F$ is
the set-valued map from $Y $  to $X $  defined by
        \begin{displaymath}
         x \in F^{-1}(y) \;   \Longleftrightarrow     \;  y \in F(x) \; 
\Longleftrightarrow  \;   (x,y) \in  \mbox{\rm Graph}(F)
        \end{displaymath}   
         \end{Definition}

We shall emphasize the characterization of a set-valued map (as well as a
single-valued map) by its graph. This point of view has been coined the
{\em graphical approach} by R.T. Rockafellar.

        The domain of $F$ is thus the image of $F^{-1 }$ and  coincides
with the projection of the graph onto the space $X $  and, in a  symmetric
way, the image of $F $  is equal to the domain of $F^{-1 }  $   and to the
projection of the graph of $F $  onto the space $Y $.
        
        Sequences of subsets can be regarded as set-valued maps defined on
the set ${\bf N}$ of integers.

        \subsection{Graphical Convergence of Maps}
        
        Since the graphical approach consists in regarding closed
set-valued maps as graphs, i.e., as closed subsets of the product space,
ranging over the space $ {\cal F} (X \times Y)$, one can supply this space
with upper and lower limits, providing the concept of upper and lower
graphical convergence:
        
        \begin{Definition}  \index{graphical convergence} \label{01A392} 
        Let us consider metric spac\-es $X, \; Y$ and a sequence of
set-valued maps $F_{n} : X \leadsto Y$. The set-valued maps $\mbox{\rm
Lim}^{\sharp}\mbox{}_{n \rightarrow \infty }F_{n}$ and $\mbox{\rm
Lim}^{\flat}\mbox{}_{n \rightarrow \infty }F_{n}$ from $X$ to $Y$ defined
by 
        \begin{displaymath} \left\{ \begin{array}{llll}
        i) & {\rm Graph}( \mbox{\rm Lim}^{\sharp}\mbox{}_{n \rightarrow
\infty }F_{n})  & := & \mbox{\rm Limsup}_{n \rightarrow \infty}{\rm
Graph}(F_{n})  \\
         & & & \\
        ii) & {\rm Graph}(\mbox{\rm Lim}^{\flat}\mbox{}_{n \rightarrow
\infty }F_{n})  & := & \mbox{\rm Liminf}_{n \rightarrow \infty}{\rm
Graph}(F_{n})
        \end{array} \right. \end{displaymath}
        are called the {\em (graphical) upper and lower limits}  of the
set-valued maps $F_{n}$ respectively. \index{graphical upper and lower
limits of  set-valued maps} \label{02A392}
         \end{Definition}
        
        Even for single-valued maps, this is a weaker convergence than the
pointwise convergence:
        \begin{Proposition}
        \mbox{}
        \begin{enumerate}
        \item If $f_{n}:X \mapsto Y$ converges pointwise to $f$, then, for
every $x \in X$, $f (x) \in f ^{\sharp } (x)$. If the sequence is
equicontinuous, then $f ^{\sharp } (x) = \{f (x)\}$.
        \item Let $ \Omega \subset  {\bf  R}^{n} $ be an open subset. If a
sequence $f_{n} \in L^{p} ( \Omega) $ converges to $f$ in $L^{p} ( \Omega
)$, then
        \begin{displaymath}
          \mbox{\rm for almost all}  \; x \in  \Omega, \;  \; f (x) \in f
^{\sharp } (x)
        \end{displaymath}
        \end{enumerate}
        \end{Proposition}

\subsection{Contingent Derivatives}
        
         Let $F:X \leadsto Y$ be a set-valued map.
        We introduce the {\em differential quotients}
        \begin{displaymath}
        u \; \leadsto \; \nabla _{h}F(x,y)(u) \; := \; \frac{F(x+hu)-y}{h}
        \end{displaymath}
        of a set-valued map $F:X \leadsto Y$ at $ (x,y) \in  \mbox{\rm
Graph}(F)$.
        \begin{Definition}
        The contingent derivative $DF(x,y)$ of $F$ at $ (x,y) \in  
        \mbox{\rm Graph}(F)$ is the graphical upper limit of differential
quotients:
        \begin{displaymath}
         DF (x,y) \; := \;  \mbox{\rm Lim}^{\sharp}\mbox{}_{h \rightarrow
0+}\nabla _{h}F(x,y)
        \end{displaymath}
         \end{Definition}

In other words, $v $ belongs to $ DF(x,y)(u)$ if and only if there exist
sequences
$h_n\rightarrow 0^+$, $u_n\rightarrow u$ and $v_n\rightarrow v$ such
that $\forall n\geq 0, ~~y+h_nv_n\in F(x+h_nu_n)$.

In particular, if $f:X\mapsto Y$ is a single valued
function, we set $Df(x)=Df(x,f(x))$.

        \mbox{}
        
        We deduce the fundamental formula on the graph of the contingent
derivative:
        \begin{Proposition}
        The graph of the contingent   derivative of a set-valued map is the
contingent  cone to its graph: for all $(x,y) \in  \mbox{\rm Graph}(F)$,
        \begin{displaymath}
         \mbox{\rm Graph}(DF (x,y))  =   T_{ \mbox{\rm Graph}(F)} (x,y) 
        \end{displaymath}
        \end{Proposition} 
        {\bf Proof} --- \hspace{ 2 mm}
        Indeed, we know that the contingent cone $$T_{ \mbox{\rm
Graph}(F)}(x,y) \; = \; \mbox{\rm Limsup}_{h \rightarrow 0+} \frac{
\mbox{\rm Graph}(F) -(x,y)}{h}$$ is the upper limit of the differential
quotients $ \frac{\mbox{\rm Graph}(F) -(x,y)}{h}$ when $h \rightarrow 0+$.
It is enough to observe that 
        \begin{displaymath}
         \mbox{\rm Graph}(\nabla _{h}F(x,y)) =  \frac{\mbox{\rm Graph}(F)
-(x,y)}{h}
        \end{displaymath}
        and to take the upper limit to conclude. $\; \; \Box$ 
        
        \mbox{}
        
         We can easily compute the derivative of  the inverse of a
set-valued map $F$ (or even of a noninjective single-valued map):  {\em
The  contingent derivative of the inverse of a set-valued map $F$  is the
inverse of the contingent derivative}:
        \begin{displaymath}
        D(F^{-1}) (y,x) \;\; = \;\; DF(x,y)^{-1}
        \end{displaymath}
        
        If $K$ is  a subset of $X$ and $f$ is a single-valued map which is
Fr\'{e}chet differentiable around a point $x \in K$, then {\em the
contingent derivative of the restriction of $f$ to $K$ is the restriction
of the derivative to the contingent cone}:
        \begin{displaymath}
        D(f|_{K})(x) = D(f|_{K})(x,f(x)) = f'(x)|_{T_{K}(x)}
        \end{displaymath}

\section{Frankowska Solutions to First-Order Partial Differential
Equations}

We consider two finite dimensional  vector spaces ${\bf   R}^{n}$ and ${\bf
R}^{p}$,  an open subset $  \Omega \subset {\bf   R}^{n} $,  its closure $K
:=\overline{ \Omega } $  closed,  its boundary $ \Gamma := \partial  \Omega
= \Gamma$, two time-dependent maps $f: {\bf R}_{+} \times K
\times  {\bf   R}^{p} \mapsto {\bf   R}^{n}$  and $g:  {\bf R}_{+} \times K
\times {\bf  R}^{p} \mapsto {\bf  R}^{p}$.

We shall study the system of first-order partial
differential equations 
        \begin{equation} \label{frankpde01}
 \frac{ \partial }{ \partial t}u (t,x) +  \frac{ \partial }{ \partial x}u
(t,x) f (t,x,u (t,x)) - g (t,x,u (t,x)) \; = \; 0
\end{equation}
on $ {\bf  R}_{+} \times K$.

 It is known that the solution $U$ to the above system can be set-valued,
describing ``shocks''. 
This is considered as a pathology  whenever the solution is regarded as a
map from  the   input space $  {\bf   R}_{+} \times K$  to the output space
${\bf  R}^{p}$,
but is quite natural when the solution $U$ is considered as a graph, i.e.,
a subset of ${\bf   R}_{+} \times K \times {\bf  R}^{p}$ and when the tools
of set-valued analysis are used.

Introducing 
\begin{enumerate}
\item an {\sf initial data} $u_{0}: K \mapsto {\bf  R}^{p}$,
\item  a {\sf  boundary  data} $v_{  \Gamma}:  {\bf   R}_{+} \times  \Gamma
\mapsto {\bf  R}^{p}$.
\end{enumerate}
we shall prove the existence and the uniqueness of a solution $:{\bf 
R}_{+} \times K  \leadsto $ to the system of first-order partial
differential equations (\ref{frankpde01})
satisfying the initial/boundary-value conditions

\begin{equation} \left\{ \begin{array}{ll} \label{initboucd003}
i) & \forall \; x \in K, \;  \; u_{0} (x) \; \in \; U(0,x)\\
ii)  &   \forall \;  t>0,  \; \forall \;  x \in  \Gamma, \;  \; v_{ \Gamma}
(t,x) \; \in  \; U (t,x)
\end{array} \right. \end{equation}

Actually,  we  associate with  the    initial  data $u_{0}:  K \mapsto {\bf
R}^{p}$ and the boundary data $v_{ \Gamma}: {\bf  R}_{+} \times  \partial
K \mapsto {\bf   R}^{p}$ the ``extended''  boundary data $ \Psi  (u_{0},v_{
\Gamma}): {\bf  R}_{+} \times K \leadsto {\bf  R}^{p}$ defined by
\begin{displaymath}
 \Psi (u_{0}, v_{ \Gamma})(s, x )   :=\left\{ \begin{array}{cllll}
 u_{0}( x) & \mbox{\rm if} & s=0   & \& & x\in K      \\
 v_{ \Gamma} (s , x ) & \mbox{\rm if} & s \geq 0 & \& & x \in \partial
K\\
\emptyset  & \mbox{\rm if} & s>0 & \& & x \in \mbox{\rm Int}(K) \\
 \end{array} \right. \end{displaymath}
which is a set-valued map since it takes (empty)  set values, the domain of
which is  $\mbox{\rm Dom}( \Psi (u_{0}, v_{ \Gamma})):= \partial (
{\bf  R}_{+} \times K) = (\{0\} \times K) \cup ({\bf  R}_{+} \times 
\Gamma)$.

The set-valued map $  \Psi $ encapsulates or replaces initial/boundary-value
data.
Hence initial and boundary  conditions (\ref{initboucd003})  can be written
in the form 
\begin{displaymath}
  \forall \;   (t,x) \; \in  \; {\bf  R}_{+} \times K, \;  \;  \Psi (u_{0},
v_{ \Gamma})(t, x ) \; \subset  \; U (t,x)
\end{displaymath}

By the way,  we  can study as  well  the case  when  $  \Psi : {\bf  R}_{+}
\times  {\bf   R}^{n} \leadsto {\bf   R}^{p}$  is any set-valued map, which
allows  to  study  other  problems   than  initial/boundary-value  problems
associated  with the  system of first-order partial  differential equations
(\ref{frankpde01}).

So,  in  the general case,  we introduce two set-valued  maps  $ \Psi :{\bf
R}_{+} \times K   \leadsto {\bf  R}^{p}$ and $ \Phi : {\bf  R}_{+} \times K
\leadsto {\bf  R}^{p}$ satisfying
\begin{displaymath}
 \forall \; (t,x) \; \in  \; {\bf  R}_{+} \times K, \;  \; \Psi (t,x) \;
\subset  \; \Phi  (t,x)
\end{displaymath}

We shall prove the existence and the uniqueness of a solution $:{\bf 
R}_{+} \times K  \leadsto $ to the system of first-order partial
differential equations (\ref{frankpde01})
satisfying the conditions
\begin{equation} \label{frankbdvc01}
 \forall \; (t,x) \in {\bf  R}_{+} \times K, \;  \; \Psi  (t,x) \; \subset 
\; U (t,x) \; \subset  \; \Phi  (t,x)
\end{equation}

The set-valued map $ \Phi $ describes ``viability constraints'' on the
solution the solution $U$ to the above system. The particular case without
constraints is naturally obtained when $ \Phi (t,x):={\bf  R}^{p}$.

{\bf Example: Impulse Boundary Value Problems}  \hspace{ 2 mm}
This is the case when we provide boundary condition $v_{
\Gamma}^{i}$ only at impulse times $t_{i}$ of an increasing sequence of
impulse times $t_{0}=0 < t_{1} \cdots < t_{n} < \cdots$.

We associate with them the map 
         $ \Psi : {\bf  R}_{+} \times K \leadsto {\bf  R}^{p}$ only derived
from the
initial  data $u_{0}:  K \mapsto {\bf   R}^{p}$  and the boundary data $v_{
\partial
K}^{i}$ by
\begin{displaymath}
 \Psi (u_{0}, \{v_{ \Gamma}\}^{i})(s, x )   :=\left\{
\begin{array}{cllll}
 u_{0}( x) & \mbox{\rm if} & s=0   & \& & x\in K      \\
 v_{ \Gamma}^{i}( x) & \mbox{\rm if} & s=t_{i}, \; i >0, & \& & x\in
\Gamma      \\
\emptyset  & \mbox{\rm if} & s \in ]t_{i},t_{i+1}[, i > 0, & \& & x \in K
\\
 \end{array} \right. \end{displaymath}
        defined on $\mbox{\rm Dom}( \Psi (u_{0}, \{v_{ \partial
K}\}^{i})):= (\{0\} \times K) \cup  \bigcup_{i}^{} ( \{t_{i}\} \times 
\Gamma)$. $\; \; \Box$ 

\mbox{}

We denote by  $  h:  {\bf R}_{+} \times K \times {\bf   R}^{p} \mapsto {\bf
R}^{n} \mapsto {\bf
R}_{+} \times K \times {\bf   R}^{p} \mapsto {\bf  R}^{n} $ the map defined
by 
\begin{displaymath}
 h ( \tau ,x,y) \; := \; (1,f ( \tau ,x,y),g ( \tau ,x,y))
\end{displaymath}
and the associated system $ ( \tau ',x',y')=h ( \tau ,x,y)$  of
differential equations 
\begin{equation} \left\{ \begin{array}{ll} \label{syschardeqeq}
i) & \tau ' (t) \; = \; 1 \\
ii) & x' (t) \; = \;  f ( \tau  (t),x (t),y (t))
\\
iii) & y' (t) \; = \; g ( \tau  (t),x (t),y (t))
\end{array} \right. \end{equation}
often called the associated ``characteristic system''.

\begin{Definition}
Given two time-dependent maps $f:  {\bf R}_{+} \times K \times {\bf  R}^{p}
\mapsto {\bf  R}^{n}$
and $g: {\bf R}_{+} \times K \times {\bf  R}^{p} \mapsto {\bf  R}^{p}$ and
two set-valued maps $  \Psi :{\bf   R}_{+} \times K  \leadsto {\bf  R}^{p}$
and $ \Phi :{\bf   R}_{+} \times K \leadsto {\bf  R}^{p}$ satisfying
\begin{displaymath}
 \forall \; (t,x) \; \in  \; {\bf  R}_{+} \times K, \;  \; \Psi (t,x) \;
\subset  \; \Phi  (t,x)
\end{displaymath}
we shall denote by $U := {\cal A}_{ (f,g)}^{ \Phi } ( \Psi ):  {\bf  R}_{+}
\times K \leadsto {\bf  R}^{p}$ the set-valued map  defined by 
\begin{equation} \label{defsolmapueq01}
\mbox{\rm Graph}(U) \; := \;  \mbox{\rm Capt}_{-h}^{ \mbox{\rm Graph}( \Phi
)}( \mbox{\rm Graph}( \Psi ))
\end{equation}
the graph of which is the viable-capture basin under $-h$ of the graph of $
\Psi $ in the graph of $ \Phi $.
 \end{Definition}

Even when $ \Psi $ is single-valued on its domain, this map $U$ can take
several values, defined as ``shocks'' in the language of physicists.
        
        Indeed, even though the trajectories of the solutions $ (x (t),y
(t))$ to the {\sf system of differential equation} 
        \begin{displaymath} \left\{ \begin{array}{ll}
        i) & x' (t) \; = \; f (t,x (t),y (t)) \\
        ii) & y' (t) \; = \; g(t,x (t),y (t)) \\
        \end{array} \right. \end{displaymath}
         never intersect in ${\bf   R}^{n} \times {\bf  R}^{p}$, their {\sf
projections} $x (t)$
onto ${\bf   R}^{n}$  --- privileged in his role of input space  --- may do
so. In other
words, if $x_{1} \ne x_{2}$, the solutions $(x_{i} (t),y_{i} (t))$ starting
from the initial conditions $ (x_{i}, u_{0}(x_{i}))$ ($i=1,2$) never
intersect, but one cannot exclude the case when for some $t$, we may have
$x_{1} (t)=x_{2} (t)$. At this time, the solution $U$  takes (at least) the
values $y_{1} (t) \ne y_{2} (t)$ associated with the common input $x :=
x_{1} (t)=x_{2} (t)$.

However, $U$ is single-valued whenever $f (t,x,y) \equiv f (t,x)$ is
independent of $y$ and the above system has a unique solution for any
initial condition.

Theorems~\ref{viablcaptbascharthm} and \ref{wonderfulthmthmbis}  can be
translated in terms of invariant manifold: 
\begin{Definition}
        We shall say that a set-valued map $ V: {\bf  R}_{+}
\times K \leadsto {\bf  R}^{p}$ defines an {\sf invariant manifold} under
the pair $ (f,g)$ \index{invariant manifold} if
for every $t_{0} \geq 0$, $ x_{0} \in K $ and $y_{0} \in V (t_{0},x_{0})$,
every solution $ (x ( \cdot ),y ( \cdot ))$ to the system of differential
equations
\begin{displaymath} \left\{ \begin{array}{ll}
i) & x' (t) \; = \; f (t,x (t),y (t)) \\
ii) & y' (t) \; = \; g (t,x (t),y (t))
\end{array} \right. \end{displaymath}
starting at $ (x_{0},y_{0})$ at time $t_{0}$ satisfies
\begin{displaymath}
\forall \;  t \geq t_{0},\;  \; y (t) \; \in \; V (t, x (t))
\end{displaymath}
         \end{Definition}
        We observe that $V$ defines  an invariant manifold under
$ (f,g)$ if and only if the graph of $ V $ is invariant under
$h$.

\begin{Theorem} \label{wontrackbistd}
Let us assume that the maps $f$ and $g$ are continuous with linear growth
and that the graphs of the set-valued maps $ \Psi  \subset  \Phi $ are
closed. The associated set-valued map $U := {\cal A}_{ (f,g)}^{ \Phi } (
\Psi ): {\bf  R}_{+}  \times K \leadsto {\bf  R}^{p}$ defined by
(\ref{defsolmapueq01}) is the {\bf largest} closed set-valued map  enjoying
the following properties:
\begin{enumerate}
\item $ \forall \; t \geq 0, \;  \forall \; x \in K, \; \; \Psi (t,x) \;
\subset  \;U(t,x) \; \subset   \; \Phi (t,x)$ 
\item  for every $t,x$ and $y \in U (t,x) \backslash \Psi (t,x)$, there
exist $s \in [0,t[$, $x_{s} \in K$  and $y_{s} \in U (s,x_{s})$ such that
the solution to the above system of differential equations starting at $
(x_{s},u(s,x (s)))$ at time $s$ satisfies 
\begin{displaymath}
\forall \;  \tau  \in [s,t],\;  \; y ( \tau ) \; \in \; U ( \tau , x ( \tau
)), \;  \; x (t) \; = \; x \; \& \; y (t)=y
\end{displaymath}
\end{enumerate}
If we  assume  furthermore that  the set-valued  map $  \Phi  $  defines an
invariant manifold under $  (f,g)$, then $U$ is the {\bf unique} set-valued
map satisfying the
above conditions, which also defines an invariant manifold under $ (f,g)$.
\end{Theorem}

     {\bf Proof} --- \hspace{ 2 mm}
     By assumption, $ \mbox{\rm Graph}( \Psi )$ is contained in $ \mbox{\rm
Graph}( \Phi )$. The graphs of the maps $ \Phi $ and $ \Psi $ are closed
subsets and repellers since they are contained in $ {\bf  R}_{+} \times K
\times {\bf  R}^{p}$, which is obviously a repeller  under the map $-h$.
     
     Therefore, by Theorem~\ref{viablcaptbascharthm}, $ \mbox{\rm
Graph}(U):= \mbox{\rm Capt}_{-h}^{ \mbox{\rm Graph}( \Phi )} ( \Psi )$ is
the largest closed subset $D:= \mbox{\rm Graph}(U)$ satisfying 
     \begin{displaymath} \left\{ \begin{array}{ll}
i) & \mbox{\rm Graph}( \Psi ) \; \subset  \; \mbox{\rm Graph}(U)\; \subset 
\; \mbox{\rm Graph}( \Phi )\\
ii) & \mbox{\rm Graph}(U) \backslash \mbox{\rm Graph}( \Psi ) \;\;\mbox{\rm
is locally viable under } \; -h \\
\end{array} \right. \end{displaymath}
     which are translated into the first properties of the above theorem.
      If we assume furthermore that $ \Phi $ defines an invariant manifold
under $ (f,g)$,
Theorem~\ref{wonderfulthmthmbis} implies that $ \mbox{\rm Graph}(U) =
\mbox{\rm Capt}_{-h} ( \Psi )$ is the unique closed subset satisfying the
above properties and being backward invariant under $-h$, i.e., invariant
under $h$.
$\; \; \Box$

Theorem~\ref{viablcaptbascharthm} and \ref{wonderfulthmthmbis}  also allow
us to derive the existence and the uniqueness of a set-valued map solution
to the system of first-order partial differential equations
(\ref{frankpde01})
        \begin{displaymath}
 \frac{ \partial }{ \partial t}u (t,x) +  \frac{ \partial }{ \partial x}u
(t,x) f (t,x,u (t,x)) - g (t,x,u (t,x)) \; = \; 0
\end{displaymath}
satisfying the conditions (\ref{frankbdvc01})
\begin{displaymath}
 \forall \; (t,x) \in {\bf  R}_{+} \times K, \;  \; \Psi  (t,x) \; \subset 
\; U (t,x) \; \subset  \; \Phi  (t,x)
\end{displaymath}

\begin{Definition}
We say that  a set-valued  map  $U:  {\bf   R}_{+}  \times  K  \mapsto {\bf
R}^{p}$  is a
Frankowska solution to the problem ((\ref{frankpde01}),(\ref{frankbdvc01}))
if the graph of $U$ is closed and if
 \begin{displaymath} \left\{ \begin{array}{ll}
i)& \forall \; (t,x) \in {\bf  R}_{+} \times K, \; \forall \; y \in U (t,x)
\backslash  \Psi  (t,x), \;  \; \\
 & 0 \; \in  \; DU (t,x,y) ( -1, -f (t,x,y)) +g (t,x,y) \\
 & \mbox{\rm and} \\
 ii) & \forall \; (t,x) \in {\bf  R}_{+} \times K, \; \forall \; y \in U
(t,x), \; 
0 \; \in  \; DU (t,x,y) ( 1, f (t,x,y)) -g (t,x,y)\\  \end{array} \right.
\end{displaymath} 
 \end{Definition}

Naturally, if $ \Psi := \psi $ and $U:=u$ are single-valued on their
domains, then a Frankowska solution can be written
 \begin{displaymath} \left\{ \begin{array}{l}
\forall \, (t,x) \in \mbox{\rm Dom}(u) \backslash \mbox{\rm Dom}( \psi ),
 \; 0 \in  Du (t,x) ( -1, -f (t,x,u (t,x))) +g (t,x,u (t,x))
\\\mbox{\rm and} \\
\forall \, t \geq 0, \,     \forall \, x \in K, 
0 \, \in  \, Du (t,x) ( 1,f (t,x,u (t,x)) -g (t,x,u (t,x))\\ 
 \end{array} \right. \end{displaymath}

If $Du (t,x) (-1,- \xi )=-Du (t,x) (1, \xi )$, these two equations boil
down to only one of them on $\mbox{\rm Dom}(u) \backslash \mbox{\rm Dom}(
\psi )$. If $u$ is differentiable in the usual sense, it satisfies the
above first order partial differential equation in the usual sense outside
the domain of $ \psi $.
        
        \mbox{}
\begin{Theorem} \label{frankwssolstatstrpbbistd}
Let us assume that the maps $f$ and $g$ are continuous with linear growth
and that the graphs of the set-valued maps $ \Psi  \subset  \Phi $ are
closed.

Then the set-valued map $U$ defined by (\ref{defsolmapueq01}) is the {\bf
largest} set-valued map with closed graph satisfying the condition
\begin{displaymath}
 \forall \; (t,x) \in {\bf  R}_{+} \times K, \;  \; \Psi  (t,x) \; \subset 
\; U (t,x) \; \subset  \; \Phi  (t,x)
\end{displaymath}
and solution to
\begin{displaymath}
\forall \, (t,x) \in {\bf  R}_{+} \times K, \; \forall \, y \in U (t,x)
\backslash  \Psi  (t,x), \;  0 \, \in  \, DU (t,x,y) ( -1, -f (t,x,y))
+g (t,x,y) 
\end{displaymath}
If we assume furthermore that $f$ and $g$ are uniformly Lipschitz with
respect to  $x$  and $y$  and that $  \Phi $  defines an invariant manifold
under $ (f,g)$, then
$U$ is the {\bf unique} Frankowska solution to the problem
((\ref{frankpde01}),(\ref{frankbdvc01})).
\end{Theorem}

     {\bf Proof} --- \hspace{ 2 mm}
     By Theorem~\ref{viablcaptbascharthm}, $ \mbox{\rm Graph}(U):=
\mbox{\rm Capt}_{-h}^{ \mbox{\rm Graph}( \Phi )} ( \mbox{\rm Graph}( \Psi
))$ is the largest closed subset $D:= \mbox{\rm Graph}(U)$ satisfying 
     \begin{displaymath} \left\{ \begin{array}{ll}
i) & \mbox{\rm Graph}( \Psi ) \; \subset  \; \mbox{\rm Graph}(U)\; \subset 
\; \mbox{\rm Graph}( \Phi )\\
ii) & \forall \; (t,x,y) \in \mbox{\rm Graph}(U) \backslash \mbox{\rm
Graph}( \Psi ),\\
 & -(1,f (t,x,t), g(t,x,t)) \; \in \; T_{ \mbox{\rm Graph}(U)} (t,x,y) \; =
\; \mbox{\rm Graph}(DU) (t,x,y)
\end{array} \right. \end{displaymath}
     which can be translated 
     \begin{displaymath}
     \forall \; (t,x), \;  \; \forall \; y \in U (t,x) \backslash \Psi 
(t,x), \; \; -g (t,x,y) \; \in  \; DU (t,x,y) (-1,-f (t,x,y))
     \end{displaymath}
      If we assume furthermore that $  \Phi $ defines an invariant manifold
under $ (f,g)$,
Theorem~\ref{wonderfulthmthmbis} implies that $ \mbox{\rm Graph}(U) =
\mbox{\rm Capt}_{-h} ( \mbox{\rm Graph}( \Psi ))$ is the unique closed
subset satisfying the above properties and being backward invariant under
$-h$, i.e., invariant under $h$. This can be translated by stating that
      \begin{displaymath} \left\{ \begin{array}{l}
          \forall \; (t,x,y) \in \mbox{\rm Graph}(U),  \\
 (1,f (t,x,t), g(t,x,t)) \; \in \; T_{ \mbox{\rm Graph}(U)} (t,x,y) \; = \;
\mbox{\rm Graph}(DU) (t,x,y)\\
      \end{array} \right. \end{displaymath} 
     i.e.,
     \begin{displaymath}
      \forall \; (t,x) \in {\bf  R}_{+} \times K, \; \forall \; y \in U
(t,x), \;  \;
     g (t,x,y)\; \in  \; DU (t,x,y) (1,f (t,x,y))
     \end{displaymath}

\section{Single-Valued Frankowska Solutions}

We already mentioned that even when $ \Psi $ is single-valued on its
domain, the solution $U$ can take several values, defined as ``shocks'' in
the language of physicists.

However, single-valuedness is naturally preserved whenever
 \begin{displaymath}
 h ( \tau ,x,y) \; := \; ( 1,\varphi (x), g( \tau ,x,y))
\end{displaymath}
when the second component of the map $h$ does not depend upon the second
variable $y$ and the differential equation $ ( \tau ',x',y')=h ( \tau ,x,y)
$ has a unique solution for any initial condition.

\mbox{}

Therefore, we proceed with the specific case when $f (t,x,y) \equiv \varphi
(x)$ depends only upon the variable $x$.

When $ (t,x) \in {\bf  R}_{+}  \times K$ is chosen, we introduce the
function $x ( \cdot ) := \vartheta_{ \varphi }( \cdot -t,x)$ the solution
to the differential equation $x' = \varphi  (x)$ starting at time $0$ at $
\vartheta_{ \varphi } (-t,x)$, or arriving at $x$ at time $t$. We associate
with it the map $g_{ (t,x)} : {\bf  R}_{+} \times {\bf  R}^{p} \mapsto {\bf
R}^{p}$ defined by 
\begin{displaymath}
\forall \;   \tau \geq 0,  \;  y \in {\bf  R}^{p}, \;  \; g_{ (t,x)} ( \tau
,y) \; := \;
g ( \tau ,\vartheta _{ \varphi } ( \tau -t,x),y)
\end{displaymath}

We denote by $ \vartheta _{g_{ (t,x)} } (t,s,y (s))$ the value at $t$ of
the solution to the  differential equation 
\begin{displaymath}
 y' ( \tau ) \; = \; g_{ (t,x)} ( \tau ,y ( \tau )) \; := \; g ( \tau ,
\vartheta_{ \varphi }( \tau -t,x), y ( \tau ))
\end{displaymath}
starting at $y (s)$ associated with the evolution $x ( \tau ):= \vartheta
_{ \varphi } ( \tau -t,x )$ starting at $x (s)=\vartheta _{ \varphi } (
s-t,x )$  at initial time $s$.

\mbox{}

We associate with the backward exit function the map $ \Theta_{K}^{ -
\varphi }$ defined by

\begin{displaymath}
 \forall \; x \in K, \;  \;
\Theta_{K}^{- \varphi } (x) \; := \; \vartheta _{- \varphi } ( \tau _{K}^{-
\varphi } (x),x)
\end{displaymath}
and we say that $ \Theta_{K}^{- \varphi }$ is the ``exitor'' (for exit
projector) of $K$. It maps $K$ to  its boundary $ \Gamma$ and
satisfies  $ \Theta_{K}^{- \varphi } (x)=x$ for every $x \in  \Theta
_{K}^{- \varphi } (K)$.

It will be very convenient to extend the function $ \tau_{K}^{- \varphi }$
defined on $K$  to the function (again denoted by) $\tau_{K}^{- \varphi }$
defined on $ {\bf  R}_{+} \times K$ by $ \tau _{K}^{- \varphi }(t,x) :=
\min (t, \tau _{K}^{- \varphi }(x))$:
\begin{displaymath}
 \tau_{K}^{- \varphi } (t,x) \; :=   \left\{ \begin{array}{cll}
t & \mbox{\rm if} & t \in [0, \tau _{K}^{- \varphi } (x)   ]  \\
 \tau _{K}^{- \varphi } (x) & \mbox{\rm if}  & t \in ] \tau _{K}^{- \varphi
} (x), \infty [
 \end{array} \right. \end{displaymath} 
so that we can also extend the exitor map  by setting
\begin{displaymath}
 \Theta_{K}^{ -\varphi } (t,x) \; := \; \vartheta _{ -\varphi } (
\tau_{K}^{- \varphi } (t,x),x)  \; \in  \; K
\end{displaymath}
because we observe that $\Theta_{K}^{ -\varphi }(t,x) $ is equal to

 \begin{displaymath} 
\left\{ \begin{array}{cll}
\vartheta _{ - \varphi } (t,x) & \mbox{\rm if} & t \in [0, \tau _{K}^{-
\varphi } (x)   ]  \\
\Theta_{K}^{ - \varphi } (x)  & \mbox{\rm if}  & t \in ] \tau _{K}^{-
\varphi } (x), \infty [
 \end{array} \right. \end{displaymath} 

\mbox{}
\begin{Proposition} \label{propgrencpropbistd}
We  posit assumption 
\begin{equation} \label{assumKonceandforall}
 K \;\;\mbox{\rm is closed and (forward) invariant under} \; \; \varphi 
\end{equation}
and 
\begin{equation} \left\{ \begin{array}{ll} \label{assumKonceandforall27}
i) & \varphi  \;\;\mbox{\rm is Lipschitz and that $g$ is continuous} \\
ii) & \mbox{\rm  $ ( \tau ',x',y')=h ( \tau ,x,y) $ has a unique solution
for any initial condition}\\
\end{array} \right. \end{equation}
Let us introduce
\begin{enumerate}
\item an {\sf initial data} $u_{0}: K \mapsto {\bf  R}^{p}$,
\item a {\sf boundary data} $v_{ \Gamma}: {\bf  R}_{+} \times  \partial
K \mapsto {\bf  R}^{p}$
\end{enumerate}
The solution $ u := {\bf  {\cal A}}_{ ( \varphi ,g)} ( \Psi  (u_{0},v_{
\Gamma}))$ is the single-valued map with closed graph   defined by
\begin{equation} \label{magmath077}
 u (t,x) \; = \; \vartheta_{ g_{ (t,x)} }( t, t- \tau _{K}^{ -\varphi }
(t,x), (\Psi (u_{0}, v_{ \Gamma}) ( t- \tau _{K}^{ -\varphi } (t,x),
\Theta_{K}^{- \varphi } (t,x)))
\end{equation}
or, more explicitly, by
\begin{displaymath} \left\{ \begin{array}{cll}
\vartheta_{ g_{ (t,x)} }( t,0, u_{0}( \vartheta _{- \varphi  } (t,x))) &
\mbox{\rm if} & t \in [0, \tau _{K}^{- \varphi } (x)   ]    \\ & & \\
\vartheta_{ g_{ (t,x)} }( t, t-  \tau _{K}^{ -\varphi } (x), v_{ \partial
K} (t- \tau _{K}^{ -\varphi } (x), \Theta_{K}^{ -\varphi  } (x))) &
\mbox{\rm if}  & t \in ] \tau _{K}^{- \varphi } (x), \infty [
 \end{array} \right. \end{displaymath} 

 Furthermore, if we assume the following viability assumptions on $\Phi $
\begin{equation} \left\{ \begin{array}{ll} \label{vainonautconstrpr}
i) &   \forall \; x \in K, \;  \forall \; t \geq 0, \; \; \forall \; y \in
\Phi  (t,x), \;  \; g (t,x,y) \; \in  \; D\Phi  (t,x,y) (1, \varphi (x))\\
ii) & \forall \;  \xi \in  \Gamma, \; \forall \; t>0, \;  \;  v_{
\Gamma} (t, \xi ) \; \in  \; \Phi  ( t,\xi )\\
iii)& \forall \; x \in K, \; \; u_{0} (x) \; \in  \; \Phi  (0,x)\\
\end{array} \right. \end{equation}
then
\begin{displaymath}
 \forall \;  (t,x) \in {\bf  R}_{+} \times K, \;  \; u (t,x) \; \in  \;
\Phi  (t,x)
\end{displaymath}

\end{Proposition} 

\mbox{}

{\bf Remark:} --- \hspace{ 2 mm}
We stress the fact  that the solution $u (t,x)$ {\sf depends only}
\begin{enumerate}
\item upon the initial condition $u_{0} (x)$ when  $t \leq \tau _{K}^{-
\varphi } (x)$,
\item upon the boundary condition $v_{ \Gamma}$ when  $t > \tau _{K}^{-
\varphi } (x)$.
\end{enumerate}
The second property proves a general principle concerning demographic
evolution  stating the state of the system {\sf eventually forgets its
initial condition $u_{0} ( \cdot )$}.  

\mbox{}

{\bf Proof} --- \hspace{ 2 mm} 
We  first take $  \Phi   (t,x)  \equiv {\bf  R}^{p}$, which is invariant by
assumption
(\ref{assumKonceandforall}).
Then the graph of $u$ is defined by  (\ref{defsolmapueq01}) is equal to
\begin{displaymath} 
\mbox{\rm Graph}(U) \; := \;  \mbox{\rm Capt}_{-h}( \mbox{\rm Graph}(
\Psi(u_{0},v_{ \Gamma}) ))
\end{displaymath}
An element $ (t,x,y)$ of the graph of $U$ is the value at  some $ h \geq 0$
of the solution $ ( \tau  ( \cdot ),x ( \cdot ),y ( \cdot ))$ to the system
of differential equations
\begin{displaymath} \left\{ \begin{array}{ll}
i) & \tau ' (t) \; = \; 1 \\
ii) & x' (t) \; = \; \varphi  (x (t)) \\
iii) & y' (t) \; = \; g ( \tau  (t),x (t),y (t))
\end{array} \right. \end{displaymath}
starting from $ (s,c, \Psi (u_{0},v_{ \Gamma})(s,c)))$, assumed to be
unique.

This implies that  $t= s + h$ and $x = \vartheta _{ \varphi } ( t,c)$.
If $s =0$, then $t=h$, $x= \vartheta_{ \varphi } (t,c)$ and
\begin{displaymath}
y \; = \; \vartheta_{ g_{ (t,x)} }( t,0, u (0, \vartheta _{- \varphi  }
(t,x)))
\end{displaymath}
If $s>0$, then $c \in  \Gamma$, so that $h= \tau _{K}^{- \varphi }
(x)$, $s=t-\tau _{K}^{- \varphi } (x)$ and $c= \Theta_{K}^{- \varphi } (x)$
and
\begin{displaymath}
y  \; = \; \vartheta_{ g_{ (t,x)} }( t, t-  \tau _{K}^{ -\varphi } (x), v_{
\Gamma} (t- \tau _{K}^{ -\varphi } (x), \Theta_{K}^{ -\varphi  } (x))) 
\end{displaymath}
Therefore, $y=:u (t,x)$ is uniquely determined by $t$ and $x$ so that  $U
=:u$ is single-valued. 

 Assumptions (\ref{vainonautconstrpr}) imply that the graph of the map $
(t,x) \leadsto \Phi  (t,x)$ is invariant under $ ( 1, \varphi ,g)$ and that
$ \mbox{\rm Graph} (\Psi _{u_{0},v_{ \Gamma}}) $ is contained in the
graph of $\Phi $. Hence the graph of $u$ is contained in the graph of
$\Phi $.
$\; \; \Box$ 

\mbox{}

{\bf Example}  \hspace{ 2 mm}  
    Let  us  consider  a $x$-dependent matrix  $A  (x)  \in  {\cal L} ({\bf
R}^{p},{\bf  R}^{p})$. We
associate the differential equation
\begin{displaymath}
 y' ( \tau ) \; = \; -A (x ( \tau )) y ( \tau ) 
\end{displaymath}
Then $u(t,x)$ is equal to
    
 \begin{displaymath} \left\{ \begin{array}{cll}
e^{- \int_{0}^{t} A ( \vartheta _{ \varphi } ( \tau -t,x ) d \tau }u_{0}(
\vartheta _{- \varphi  } (t,x))
 & \mbox{\rm if} & t \in [0, \tau _{K}^{- \varphi } (x)] \\
 e^{- \int_{t - \tau _{K}^{- \varphi } (x)}^{t} A ( \vartheta _{ \varphi }
( \tau - t,x ) d \tau } v_{ \Gamma} (t- \tau _{K}^{ -\varphi } (x),
\Theta_{K}^{- \varphi  } (x)) & \mbox{\rm if}  & t \in ] \tau _{K}^{-
\varphi } (x), \infty [
 \end{array} \right. \end{displaymath} 

\mbox{}

{\bf Example}  \hspace{ 2 mm}
Let us consider the case when $ K \; := \; K_{1} \times  \cdots \times
K_{n}$  here the $K_{i} \subset {\bf  R}^{n}_{i}$ are close subsets and set
${\bf  R}^{n} :=
\prod_{j=1}^{n}{\bf  R}^{n}_{j}$.

\begin{Lemma} \label{maurinlemma001}
Assume that $K := \prod_{j=1}^{n}K_{j}$ is the product of $n$ closed
subsets $K_{j} \subset {\bf  R}^{n}_{j}$.
We posit assumptions (\ref{assumKonceandforall}) on $K$.
        Denoting by
\begin{displaymath}
 \tau _{K_{j}}^{- \varphi } (x) \; := \; \tau _{K_{j}} ( ( \vartheta_{ -
\varphi } ( \cdot ,x))_{j})
\end{displaymath}
the {\sf partial backward exit time functions}, \index{partial backward
exit time functions}  the backward exit time function can then be written
\begin{displaymath}
 \tau _{K}^{- \varphi } (x) \; := \; \min_{j=1, \ldots ,n} \tau _{K_{j}}^{-
\varphi } (x) 
\end{displaymath}
\end{Lemma}
{\bf Proof} --- \hspace{ 2 mm}
We observe that
\begin{displaymath}
   \prod_{j=1}^{n}   {\bf    R}^{n}_{j}  \backslash  \left( \prod_{j=1}^{n}
K_{j}\right) \; =
\;  \bigcup_{j=1}^{n} \left( \prod_{i=1}^{j-1}{\bf  R}^{n}_{i} \times ({\bf
R}^{n}_{j} \backslash
K_{j}) \times  \prod_{l=j+1}^{n}{\bf  R}^{n}_{l} \right) 
\end{displaymath}
and thus, that for any function $t \mapsto x (t)= (x_{1} (t), \ldots ,x_{n}
(t))$, 
\begin{displaymath}
  \tau_{K} (x (  \cdot ))  \; := \; \inf_{x (t) \in {\bf  R}^{n} \backslash
K}t \; = \;
\min_{j=1, \ldots ,n} \left( \inf_{x_{j}(t) \in {\bf  R}^{n}_{j} \backslash
K_{j}}t
\right) \; = \; \min_{j=1, \ldots ,n} \tau _{K_{j}} (x_{j} ( \cdot ))
\end{displaymath}
since the infimum on a finite union of subsets is the minimum of the infima
on each subsets. $\; \; \Box$ 

 In this case, 
\begin{displaymath}
 \partial  \left( \prod_{j=1}^{n} K_{j}\right) \; = \; \bigcup_{j=1}^{n}
\left( \prod_{i=1}^{j-1}K_{i} \times \Gamma_{j} \times 
\prod_{l=j+1}^{n}K_{l} \right) 
\end{displaymath}
so that the  boundary data defined on $ \Gamma$ are defined by $n$ maps
\begin{displaymath}
 v_{ \Gamma}^{j} : \prod_{i=1}^{j-1}K_{i} \times \Gamma_{j} \times 
\prod_{l=j+1}^{n}K_{l}  \mapsto {\bf  R}^{p}
\end{displaymath}

\mbox{}

{\bf Example}  \hspace{ 2 mm}
     For instance, let us consider the case when the four-dimensional
causal variable $x:= (x_{1},x_{2},x_{3},x_{4})$ ranges over the product $
K:= \prod_{i=1}^{4} K_{i} $ with $ K_{1}:={\bf  R}_{+} $, $
K_{2}:=[0,r_{2}] $, $K_{3}:={\bf  R}_{+} $ and $ K_{4}:=[0,b] $. 
        
        We are looking for solutions to the system of first-order partial
differential equations
         \begin{displaymath} \left\{ \begin{array}{l}
        \displaystyle{     \frac{ \partial }{ \partial t}u (t,x) + \frac{
\partial }{ \partial x_{1}}u(t,x) - \rho  \frac{ \partial }{ \partial
x_{2}}u (t,x)+ \sigma   \frac{ \partial }{ \partial x_{3}}u (t,x)+ \beta 
(b-x_{4})x_{4} \frac{ \partial }{ \partial x_{4}}u(t,x)
        }\\ = \; g (t,x,u(t,x))
         \end{array} \right. \end{displaymath} 
(where the scalars  functions $ \rho , \; \sigma  $ and $ \beta $ are
positive) satisfying the initial and boundary conditions
        \begin{displaymath} \left\{ \begin{array}{ll}
        i) & \forall \; x_{1}\geq 0, \; \forall \; x_{2}\in [0,r_{2}], \;
x_{3} \geq 0, \; x_{4} \in [0,b],\\ &   u (0,x_{1},x_{2},x_{3},x_{4}) \; =
\; u_{0}(x_{1},x_{2},x_{3},x_{4}) \\
        ii) & \forall \; t \geq 0, \; \forall \; x_{2}\in [0,r_{2}], \;
x_{3} \geq 0, \; x_{4} \in [0,b], \\ &  u (t,0,x_{2},x_{3},x_{4}) \; = \;
v_{1}(t,x_{2},x_{3},x_{4}) \\
        i) & \forall \; t \geq 0, \; \forall \; x_{1}\geq 0, \;  x_{3} \geq
0, \; x_{4} \in [0,b], \\ &  u (t,x_{1},r_{2},x_{3},x_{4}) \; = \;
v_{r_{2}}(t,x_{1},x_{3},x_{4}) \\
        \end{array} \right. \end{displaymath}

Hence, we derive the existence and the uniqueness of  the Frankowska
solution to the above system of partial differential equations satisfying
an initial condition, a boundary condition for $x_{1}=0$ (births) and a
boundary condition for $x_{2}=r_{2}$.

        For computing it, we need to know the backward exit time and the
exitor for the associated characteristic system given by
        \begin{equation} \left\{ \begin{array}{ll}
        i) &  x'_{1} (t) \; = \; 1 \\
        ii) & x'_{2} (t) \; = \; - \rho  x_{2} (t)\\
        iii)& x'_{3} (t) \; = \; \sigma x_{3} (t)\\
        iv)& x'_{4} (t) \; = \; \beta (b -x_{4} (t))x_{4} (t)
        \end{array} \right. \end{equation}
where $ \varphi  (x) := ( \varphi _{i} (x))_{i=1, \ldots ,4}$ with $
\varphi_{1} (x_{1}):=1$, $ \varphi _{2} (x_{2}):=- \rho x_{2}$,
$\varphi_{3} (x_{3}):= \sigma x_{3}$ and $ \varphi _{4} (x_{4}):= \beta 
(b-x_{4})x_{4}$.

We recall that the solution to the purely logistic equation $y' (t)= \beta
(t) (b-y (t))y (t)$ starting at $y_{s}$ at time $s$ is given by
\begin{displaymath}
 y (t) \; := \;  \frac{b}{1 + \left( \frac{b}{y_{s}}-1 \right)e^{-b
\int_{s}^{t} \beta  ( \tau ) d \tau } }
\end{displaymath}

The closed subset $K$ is obviously (forward) invariant under $ \varphi $
defined by $ \varphi  $ and one can observe easily that
\begin{displaymath}
 \tau _{K_{3}}^{- \varphi } (x) \; = \; + \infty \; \& \;  \tau _{K_{4}}^{-
\varphi } (x) \; = \; + \infty 
\end{displaymath} 
so that
\begin{displaymath}
 \tau _{K}^{- \varphi } (t,x_{1},x_{2},x_{3},x_{4}) \; = \; \min \left( t,
x_{1}, \frac{1}{ \rho }\log \left( \frac{r_{2}}{x_{2}} \right) \right) 
\end{displaymath}
Therefore, 
\begin{enumerate}
\item if $t \leq  \min \left(  x_{1}, \frac{1}{ \rho }\log \left(
\frac{r_{2}}{x_{2}} \right) \right)$, then 
\begin{displaymath}
 \Theta_{K}^{- \varphi } (t,x) \; = \; \left( 0,x_{1}-t,e^{ \rho
t}x_{2},e^{- \sigma t}x_{3}, \frac{b}{1+ \left( \frac{b}{x_{4}}-1
\right)e^{ \beta bt}} \right)
\end{displaymath}
\item if $x_{1} \leq \min \left( t, \frac{1}{ \rho }\log \left(
\frac{r_{2}}{x_{2}} \right) \right)$, then 
\begin{displaymath}
 \Theta_{K}^{- \varphi } (t,x) \; = \; \left( t-x_{1},0,e^{ \rho
x_{1}}x_{2},e^{- \sigma x_{1}}x_{3}, \frac{b}{1+ \left( \frac{b}{x_{4}}-1
\right)e^{ \beta bx_{1}}} \right)
\end{displaymath}

\item  if $ \frac{1}{ \rho }\log \left( \frac{r_{2}}{x_{2}} \right)\leq
\min \left( t, x_{1}  \right)$, then 
 \begin{displaymath} \left\{ \begin{array}{l}
 \Theta_{K}^{- \varphi } (t,x) \; = \\ \displaystyle{ \left( t-\frac{1}{
\rho }\log \left( \frac{r_{2}}{x_{2}} \right),x_{1}-\frac{1}{ \rho }\log
\left( \frac{r_{2}}{x_{2}} \right),r_{2}, e^{ \left( \frac{x_{2}}{r_{2}}
\right)^{ \frac{ \sigma }{ \rho }}}x_{3}, \frac{b}{1+ \left(
\frac{b}{x_{4}}-1 \right) \left( \frac{r_{2}}{x_{2}} \right)^{ \frac{ \beta
b}{ \rho }}} \right)}     \\
 \end{array} \right. \end{displaymath} 
\end{enumerate}
If the right-hand side $g (t,x,y) :=-A(t,x)y$ where $A (t,x) \in  {\cal L}
({\bf  R}^{p},{\bf  R}^{p})$ is linear, then we set
\begin{displaymath}
 A (t,x; \tau ) \; := \; A \left( \tau , x_{1}-t+ \tau , e^{- \rho (t- \tau
)}x_{2}, e^{ \sigma (t- \tau )}x_{3}, \frac{b}{1+ \left( \frac{b}{x_{4}}-1
\right) e^{ \beta b (t- \tau )}}
 \right)
\end{displaymath}
Then the solution is given by
\begin{enumerate}

\item if $t \leq  \min \left(  x_{1}, \frac{1}{ \rho }\log \left(
\frac{r_{2}}{x_{2}} \right) \right)$, then 
\begin{displaymath}
u (t,x) \; = \; e^{- \int_{0}^{t}A (t,x; \tau )d \tau } u_{0} 
 \left( x_{1}-t,e^{ \rho t}x_{2},e^{- \sigma t}x_{3}, \frac{b}{1+ \left(
\frac{b}{x_{4}}-1 \right)e^{ \beta bt}} \right)
\end{displaymath}
\item if $x_{1} \leq \min \left( t, \frac{1}{ \rho }\log \left(
\frac{r_{2}}{x_{2}} \right) \right)$, then 
\begin{displaymath}
u (t,x) \; = \; e^{- \int_{t-x_{1}}^{t}A (t,x; \tau )d \tau } v_{1}  \left(
t-x_{1},e^{ \rho x_{1}}x_{2},e^{- \sigma x_{1}}x_{3}, \frac{b}{1+ \left(
\frac{b}{x_{4}}-1 \right)e^{ \beta bx_{1}}} \right)
\end{displaymath}

\item  if $ \frac{1}{ \rho }\log \left( \frac{r_{2}}{x_{2}} \right)\leq
\min \left( t, x_{1}  \right)$, then 
 \begin{displaymath} \left\{ \begin{array}{l}
\displaystyle{u (t,x) \; = \; e^{- \int_{t- \frac{1}{ \rho }\log \left(
\frac{r_{2}}{x_{2}} \right)}^{t}A (t,x; \tau )d \tau } }\\
\displaystyle{v_{r_{2}} 
 \left( t-\frac{1}{ \rho }\log \left( \frac{r_{2}}{x_{2}}
\right),x_{1}-\frac{1}{ \rho }\log \left( \frac{r_{2}}{x_{2}} \right), e^{
\left( \frac{x_{2}}{r_{2}} \right)^{ \frac{ \sigma }{ \rho }}}x_{3},
\frac{b}{1+ \left( \frac{b}{x_{4}}-1 \right) \left( \frac{r_{2}}{x_{2}}
\right)^{ \frac{ \beta b}{ \rho }}} \right)}     \\
 \end{array} \right. \end{displaymath} 
\end{enumerate}

\mbox{}

\section{Regularity Properties}

We begin by proving that the operator $ {\bf {\cal A}}_{ ( \varphi ,g)}$
preserves continuity  and boundedness:

\begin{Proposition} \label{propgrencpropbis2td}
We posit assumptions (\ref{assumKonceandforall}) and
(\ref{assumKonceandforall27}) and the regularity property 
\begin{equation} \label{assumKonceandforallbis}
\tau _{K}^{- \varphi }: K \mapsto  \Gamma \;\;\mbox{\rm is
continuous}\;\;
\end{equation}
Assume also that $g$ enjoys uniform linear growth with respect to $y$ in
the sense that there exists a positive constant $c$ such that 
\begin{equation} \label{lingrassg}
  \forall \;  x \in K,  \;  \forall \;  y \in {\bf   R}^{p},  \;   \; \|g (
t,x,y)\| \; \leq 
\;c (1+ \|y\|)
\end{equation}
Then the ${\bf  {\cal A}}_{ (\varphi ,g)} ( \Psi (u_{0},v_{ \Gamma}))$ 
is continuous and bounded whenever $u_{0}$ and $v_{ \Gamma}$ are
continuous and bounded.
\end{Proposition} 

{\bf Proof} --- \hspace{ 2 mm}
By assumption (\ref{assumKonceandforallbis}), $ \tau _{K} (t,x):= \min (t,
\tau _{K}^{- \varphi  } (x))$ is also continuous, so that $ \Theta_{K}^{-
\varphi }(t,x) := \vartheta_{- \varphi } ( \tau _{K}^{- \varphi  }
(t,x),x)$ is also continuous, and thus
\begin{displaymath}
 u (t,x) \; = \;  \vartheta _{g_{ (t,x)}} (t,t- \tau _{K} (t,x), \Psi 
(u_{0},v_{ \Gamma}) (t- \tau _{K} (t,x),\Theta_{K}^{- \varphi }(t,x)
))
\end{displaymath}
is also continuous.
 When $ \|g (t, x,y)\| \; \leq  \;c (1+ \|y\|)$, we also infer that
\begin{displaymath}
 \|u (t,x)\| \; \leq  \; e^{c \tau _{K}^{- \varphi } (t,x)} \max ( \|u_{0}
(x)\|, \|v_{ \Gamma} ( \Theta_{K}^{- \varphi} (x))\|)
\end{displaymath}
and thus, that
\begin{displaymath}
 \|u (t, \cdot )\|_{ \infty } \; \leq  \; e^{ct} \max ( \|u_{0}\|_{ \infty
}, \|v_{ \Gamma}\|_{ \infty })
\end{displaymath}

\mbox{}

Next, we prove that under monotonicity assumptions, the operator 
$ {\bf {\cal A}}_{ ( \varphi ,g)}$ is Lipschitz:

\begin{Definition}
        We say that $g$ is {\sf uniformly monotone}
\index{uniformly monotone} with respect to $t$ and $y$ if there
exists $ \mu  \in  {\bf R}$ such that 
\begin{equation} \label{monassg37}
  \langle g (t,x, y_{1}) - g (t,x,y_{2}) ,y_{1}-y_{2} \rangle \; \leq  \;
-\mu \|y_{1}-y_{2}\|^{2}
\end{equation}
 \end{Definition}

 The interesting case is obtained when $   \mu >0$. When $g$ is uniformly
Lipschitz with respect to $t$ and $y$ with constant $ \lambda $, then it is
uniformly monotone with $ \mu =- \lambda $.

\begin{Proposition} \label{propgrencpropasybistdmon}
We posit assumption  (\ref{assumKonceandforall}) and assume that $ \varphi
$ is Lipschitz and that $g$ is continuous and that $g$ is {\sf uniformly
monotone} \index{uniformly monotone} with respect to $y$.

 Then, for each $t>0$, 
 \begin{displaymath} \left\{ \begin{array}{l}
\sup_{ x \in K}\| ({\bf {\cal A}}_{ ( \varphi ,g)}\Psi (u^{1}_{0}, v^{1}_{
\Gamma})) ( (t,x)- ({\bf {\cal A}}_{ ( \varphi ,g)} \Psi (u^{2}_{0},
v^{2}_{ \Gamma}) (t,x)\| \\  \\ 
\leq  \; e^{- \mu t} \max ( \|u^{1}_{0}-u^{2}_{2}\|_{ \infty } , \|v^{1}_{
\Gamma} (t, \cdot )-v^{2}_{ \Gamma}) (t, \cdot )\|_{ \infty }   )\\
 \end{array} \right. \end{displaymath}

Consequently, if we posit assumptions  (\ref{assumKonceandforallbis}) and
(\ref{lingrassg}), we deduce that for any $T>0$
$$  {\bf   {\cal A}}_{ (  \varphi ,g)} :{\cal C}_{ \infty }(K ,{\bf  R}^{p}
)\times {\cal
C}_{ \infty }([0,T] \times   \Gamma,{\bf  R}^{p}) \mapsto {\cal C}_{ \infty
}([0,T]
\times K,{\bf  R}^{p})$$ 
is  a Lipschitz operator  from  the space $  {\cal  C}_{  \infty  }(K ,{\bf
R}^{p} )\times
{\cal  C}_{  \infty  }([0,T]  \times   \Gamma,{\bf   R}^{p})$  of  pairs of
continuous and
bounded initial and boundary data to the space $ {\cal C}_{ \infty }([0,T]
\times K,{\bf   R}^{p})$  of continuous and bounded maps from $[0,T] \times
K$ to ${\bf  R}^{p}$.
\end{Proposition} 

{\bf Proof} --- \hspace{ 2 mm}
 Indeed, setting $u^{i} (t,x):= ({\bf {\cal A}}_{ (\varphi ,g)} ( \Psi
(u_{0}^{i}, v_{ \Gamma}^{i})))(x) $, $i=1,2$, we know that
\begin{displaymath}
 \forall \; x \in K , \;  \;u^{i} (t,x)\; := \; \vartheta _{g_{t,x}} (t,t-
\tau _{K}^{- \varphi } (t,x), \Psi ( u_{0}^{i}, v_{ \Gamma}^{i}) (
\Theta_{K}^{- \varphi } (t,x)) )
\end{displaymath}
is the value $ y^{i}(t):=u^{i} (t,x)$ at time $ t$ of the solution $ y (
\cdot )$ to the differential equation $y' ( \tau )=g ( \tau ,x ( \tau ), y
( \tau ))$ starting from  $ \Psi (u_{0}^{i},v_{ \Gamma}^{i}) ( \Theta
_{K}^{ - \varphi } (t,x))$ at time $t-\tau _{K}^{- \varphi } (t,x)$, where
$x ( \cdot )$ is the solution to the differential equation $x' = \varphi 
(x)$ starting from $ \Theta_{K}^{- \varphi } (t,x)$ at time $t-\tau
_{K}^{- \varphi } (t,x)$.

Recalling that $ \tau _{K}^{- \varphi } (t,x)= \min (t, \tau _{K}^{-
\varphi } (x)) \leq t$, we infer from assumption (\ref{monassg37})
that\footnote{Indeed, integrating the two sides of  inequality
\begin{displaymath}
 \frac{d}{dt} \|y_{1} (t)-y_{2} (t)\|^{2} \; = \; 2  \langle g (t,x
(t),y_{1} (t))-g (t,x (t), y_{2} (t)),y_{1} (t)-y_{2} (t) \rangle  \; \leq 
\; - 2 \mu \|y_{1} (t)-y_{2} (t)\|^{2}
\end{displaymath}
yields 
\begin{displaymath}
 \|y_{1} (t)-y_{2} (t)\|^{2} \; \leq  \; e^{-2 \mu t} \|y_{1} (0)-y_{2}
(0)\|^{2}
\end{displaymath}}
 \begin{displaymath} \left\{ \begin{array}{l}
     \|y^{1}(t)- y^{2}(t)\| \\ \leq  \; e^{ - \mu  \tau _{K}^{- \varphi }
(t,x)} \| \Psi (u_{0}^{1},v_{ \Gamma}^{1})( \Theta_{K}^{ - \varphi }
(t,x))- \Psi (u_{0}^{2},v_{ \Gamma}^{2})( \Theta_{K}^{ - \varphi }
(t,x))\|\\
 \leq \; e^{ - \mu  t} \max (\sup_{x \in K} \|u_{0}^{1} (x)-u_{0}^{2}
(x)\|,
  \sup_{ \xi \in \Gamma}\|v_{ \Gamma}^{1} ( \xi  )-v_{ \partial
K}^{2} ( \xi )\| ) \\ = \; e^{ - \mu  t } \max ( \|u_{0}^{1}-u_{0}^{2}\|_{
\infty }+\|v^{1}_{ \Gamma}-v^{2}_{ \Gamma}\|_{ \infty })
 \end{array} \right. \end{displaymath}

        \end{document}